\documentclass[11pt,twoside]{article}
\usepackage{times}
\usepackage{amsmath,amssymb,amsthm}
\usepackage{enumerate}
\usepackage{cite}
\usepackage{mathrsfs,graphicx}
\usepackage{float}

\allowdisplaybreaks

\newcommand{\no}{\nonumber}

\newcommand{\R}{\mathbb R}
\newcommand{\N}{\mathbb N}

\newcommand{\p}{\partial}
\newcommand{\ve}{\varepsilon}
\newcommand{\f}{\frac}
\newcommand{\la}{\lambda}

\newcommand{\al}{\alpha}

\newcommand{\vp}{\varphi}

\newcommand{\dl}{\delta}
\newcommand{\ds}{\displaystyle}
\newcommand{\RN}[1]{\textup{\uppercase\expandafter{\romannumeral#1}}}

\newcommand{\mcD}{\mathcal{D}}
\newcommand{\mcS}{\mathcal{S}}
\newcommand{\bs}{\boldsymbol}

\theoremstyle{plain}

\theoremstyle{definition}

\theoremstyle{remark}

\numberwithin{equation}{section}

\title{On the steady non-Newtonian fluids in domains with noncompact boundaries}

\author{Yang Jiaqi$^{1, *}$ \qquad Yin Huicheng$^{2, }$\footnote{Yang Jiaqi (\texttt{yjqmath@163.com
}) and
    Yin Huicheng (\texttt{huicheng$@$nju.edu.cn}, \texttt{05407$@$njnu.edu.cn}) are supported by
    the NSFC (No.~11571177) and the Priority Academic Program
    Development of Jiangsu Higher Education Institutions.}\\
    [12pt] {\small 1. Department of Mathematics and IMS,
  Nanjing University, Nanjing 210093, China}\\
  {\small 2. School of Mathematical Sciences, Jiangsu Provincial Key Laboratory for
  Numerical}\\
  {\small  Simulation of Large Scale Complex Systems,
  Nanjing Normal University, Nanjing 210023, China}}

\begin{document}

\date{}
\maketitle
\thispagestyle{empty}

\begin{abstract}
In this paper, we study the steady non-Newtonian fluids in a class of unbounded domains
with noncompact boundaries.  With respect to the resulting mathematical problems, we establish
the global existence of solutions with arbitrary large flux under some suitable conditions,
and meanwhile, show the uniqueness of the solutions when the flux is sufficiently small.
Our results are an extension or an improvement of those obtained in some previous references.


\end{abstract}

\section{Introduction}
Although the steady Navier-Stokes equations have been investigated extensively
(see \cite{Amick}-\cite{G-2},  \cite{Gil-mar}, \cite{La-1}-\cite{La-5}, \cite{Mar-pa}-\cite{Str}
and the references therein),
the global well-posedness  of a flow in a domain $\Omega\subset\Bbb R^d$ $(d=2,3)$ with noncompact boundaries is
still an interesting question for arbitrary fluxes. A special case is that the domain  $\Omega$ is a 
distorted infinite cylinder or channel
(see \cite{G-2} and so on), namely, $\Omega$ can be described as follows (see Figure 1 and Figure 2 below):
\begin{equation}
 \Omega=\bigcup_{i=0}^2\Omega_i,
\end{equation}
where $\Omega_0$ is a smooth bounded subset of $\Omega$, while $\Omega_1$ and $\Omega_2$
are disjoint regions which may be expressed in possibly different coordinate systems
$(x_1^1, ..., x_d^1)$ and $(x_1^2, ..., x_d^2)$  by
\begin{equation}
\Omega_i=\{(x_1^i, ..., x_d^i)\in\mathbb{R}^d: x^i_1>0, (x_2^i,..,x_d^i)\in\Sigma_i(x_1^i)\},\qquad i=1, 2,
\end{equation}
here $\Sigma_i(x_1^i)$ represents the bounded cross section of $\Omega_i$ for fixed $x_1^i$.
\vskip 0.2 true cm \includegraphics[width=12cm,height=9cm]{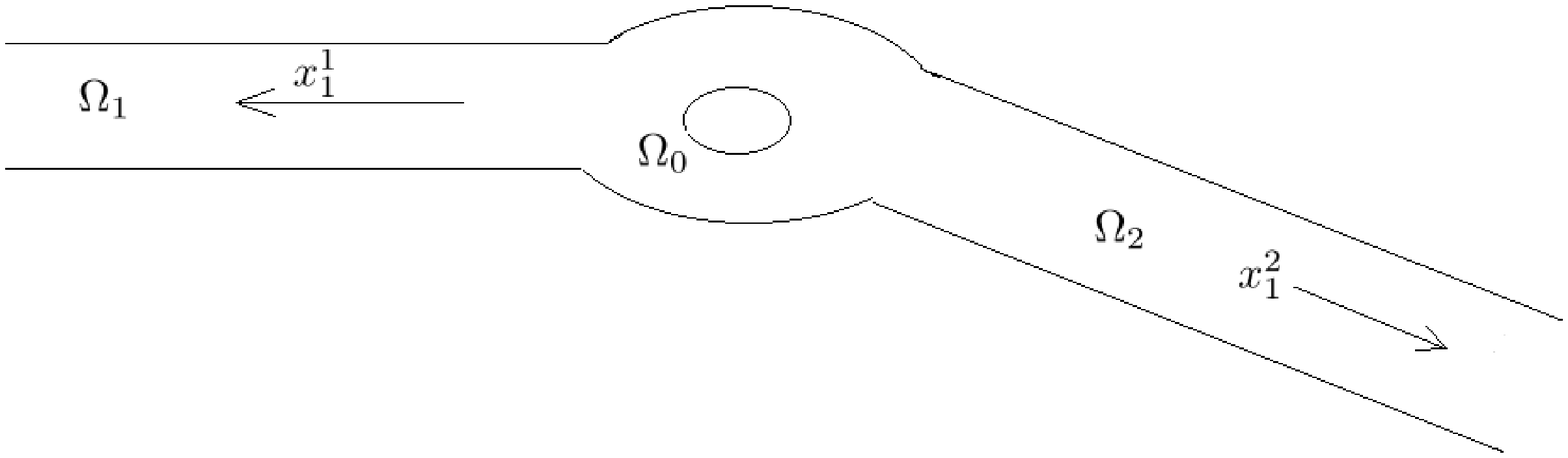}

\centerline{\bf Figure 1. Domain $\Omega$ for $d=2$}
\vskip 0.2 true cm

\includegraphics[width=14cm,height=10cm]{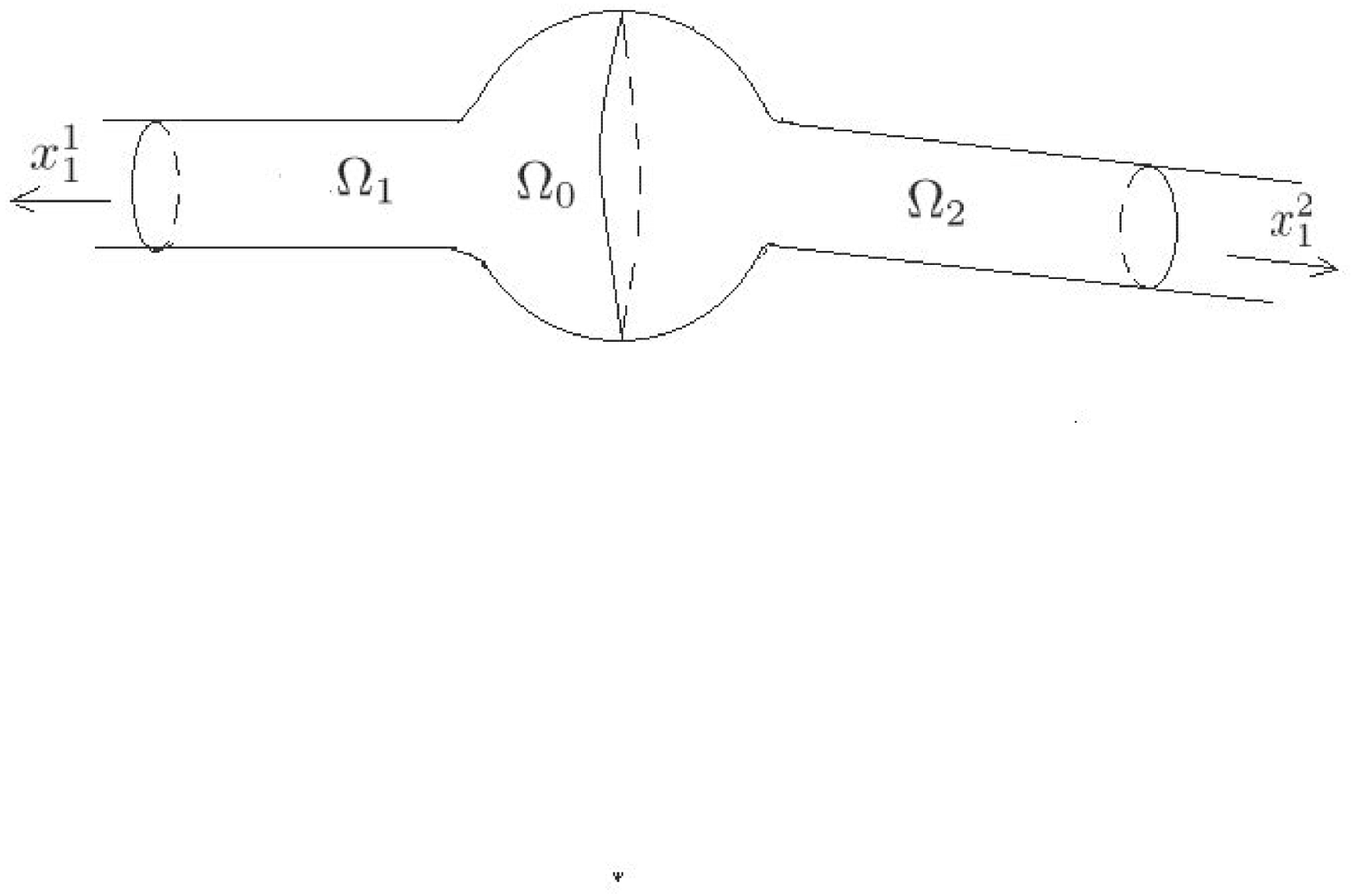}

\centerline{\bf Figure 2. Domain $\Omega$ for $d=3$}
\vskip 0.3 true cm

Owing to the incompressibility of the fluids and the vanishing property of the current velocity $\mathbf{v}=(v_1, ..., v_d)$
on the boundary  $\p \Omega$,
we deduce that the flux $\alpha_i\equiv\int_{\Sigma_i(x_1^i)}\mathbf{v}\cdot\mathbf{n}\ dS$ of velocity  $\mathbf{v}$
($\mathbf{n}$ stands for the unit outer normal direction of $\Sigma_i(x_1^i)$)  through
$\Sigma_i(x_1^i)$ is a constant independent of the variable $x_1^i$, and $\alpha_i'$s ($i=1, 2$) satisfy
\begin{equation}
\alpha_1+\alpha_2=0.
\end{equation}

When the cross section $\Sigma_i(x_1^i)$ is independent of $x_1^i$,
which means that each outlet $\Omega_i$ is a semi-infinite strip for $d=2$
or a semi-infinite straight cylinder for $d=3$, respectively, $\Sigma_i(x_1^i)$ will be
simply denoted by $\Sigma_i$.
In this case, the classical Leray's problem (see \cite{Leary}) is to study the well-posedness of the following steady flows:
\begin{equation}
\begin{cases}
-\mu\Delta \mathbf{v}+\mathbf{v}\cdot\nabla\mathbf{v}+\nabla \pi=0& \text{in $\Omega$}, \\
$div$ \ \mathbf{v}=0 & \text{in $\Omega$},\\
\mathbf{v}=0,& \text{on $\partial{\Omega}$},\\
\int_{\Sigma_i}\mathbf{v}\cdot \mathbf{n} \ dS=\al_i& i=1,2,\\
\mathbf{v}\rightarrow \mathbf{v}_0^{i}&\text {as $x_1^i\to\infty$ in $\Omega_i$ for $i=1,2$},
\end{cases}
\end{equation}
where $\mathbf{v}_0^{i}$  stands for the velocity of the Poiseuille flow corresponding to the given constant
$\alpha_i$, which is determined by
\begin{equation}
\begin{cases}
\mathbf{v}_0^i=v_0^i(z^i)\mathbf{e}_1,\\
\ds\mu\sum_{j=2}^d\p_{x_j^i}^2v_0^i(z^i)=-C_i &\text{in  $\Sigma_i$},\\
v_0^i(z^i)=0 & \text{on $\partial\Sigma_i$}
\end{cases}
\end{equation}
with $z^i=(x_2^i, ..., x_d^i)$, $\mathbf{e}_1=(1,0,..., 0)$ and $C_i$ being a constant uniquely determined by $\alpha_i$.

Leray's problem (1.4)-(1.5) has been extensively studied.
In \cite{La-1}, under the smallness assumption of the flux $\al_i$ $(i=1,2)$, O.A.Ladyzhenskaya proved the  existence of solution $\mathbf{v}$
but the uniqueness argument was not given.  In  \cite{Amick}, C.J.Amick completed the proof of both existence and uniqueness
when the flux $\al_i$ is sufficient small. Alternately, O.A.Ladyzhenskaya and V.A.Solonnikov in
\cite{La-5} considered problem (1.4) together with (1.5) under the weaker assumption that the section $\Sigma_i(x_1^i)$ is uniformly bounded
with respect to the variable $x_1^i$ instead of the straight outlet $\Sigma_i$ in \cite{Amick}.
In this case, one cannot pose the condition of $\mathbf{v}$
at infinity by Poiseuille flow since section $\Sigma_i(x_1^i)$ changes for different $x_1^i$.
Consequently, the authors in \cite{La-5} considered problem  (1.4) in another way, which is called
Ladyzhenskaya-Solonnikov Problem I (by prescribing a growth condition
of $|\mathbf{v}|$ with respect to the distance along the direction of each outlet instead of
condition (1.5)),
and they established the global existence of $\mathbf{v}$ for arbitrary large flux by utilizing a variant of Saint-Venant's principle.
Furthermore, if the flux is sufficient small, they got the uniqueness of solution  $\mathbf{v}$.
In particular, if flux is small and both exits $\Omega_1$ and $\Omega_2$ are straight,
then it has been shown that the solution
$\mathbf{v}$ to Ladyzhenskaya-Solonnikov Problem I  tends to the corresponding Poiseuille solution of (1.5).

In \cite{La-5}, the authors also studied another problem for (1.4), i.e.,
Ladyzhenskaya -Solonnikov Problem II.
At this time, the sections $\Omega_1$  and $\Omega_2$ of $\Omega$ are not uniformly bounded and admit
some certain rates of ``growth", i.e.,
\begin{equation}
\Omega_i=\{x=(x_1^i, y^i)\in\mathbb{R}^d: x^i_1>0, |y^i|\equiv\sqrt{(x_2^i)^2+...+(x_d^i)^2}<g_i(x^i_1)\},
\end{equation}
where $g_i(x_1^i)$ is a global Lipschitz function.
Later, in a series of papers \cite{P-1,P-2,P-3,P-4}, K.Pileckas shows that the Ladyzhenskaya -Solonnikov Problem II
is uniquely solvable if flux is small.
Simultaneously, it is shown in  \cite{P-1,P-2,P-3,P-4} that the decay rate of solution $\mathbf{v}$
at infinity is related to the inverse power of the functions $g_1$
and  $g_2$.

The Navier-Stokes  model of incompressible fluids is  based on the Stokes-hypothesis£¬ which  simplifies the relation between
the stress tensor and the velocity. However, a number of experiments show that many other incompressible fluids, including bloods,
cannot be described by this model. In the late 1960s, see \cite{La-2,La-3}, O.A.Ladyzhenskaya started a systematic investigation on
the well-possedness of the boundary value problems associated to certain generalized Newtonian models.
In contrast to Newtonian flows, for non-Newtonian flows, the viscosity coefficient $\mu$ is no longer constant,
it depends on the magnitude of $\mathcal{D}(\mathbf{v})$, i.e.,
\begin{equation}
\mu \big{(}\mathcal{D}(\mathbf{v})\big{)}=\mu_0+\mu_1|\mathcal{D}(\mathbf{v})|^{p-2},
\end{equation}
where $\mu_0\geqslant0$, $\mu_1>0$, $p>1$, and $\mathcal{D}(\mathbf{v})=(D_{ij}(\mathbf{v}))_{i,j=1}^d$
with $D_{ij}(\mathbf{v})=\f12(\p_iv_j+\p_jv_i)$. In this case,
the corresponding Leray's problem in the unbounded pipe domain $\Omega$ is described as follows
\begin{equation}
\begin{cases}
-$div$ \ \big{(}\mu_0\mathcal{D}(\mathbf{v})+\mu_1|\mathcal{D}(\mathbf{v})|^{p-2}\mathcal{D}(\mathbf{v})\big{)}
+\mathbf{v}\cdot\nabla\mathbf{v}+\nabla\pi=0& \text{in $\Omega$},\\
$div$ \ \mathbf{v}=0 & \text{in $\Omega$},\\
\mathbf{v}=0 &\text{on $\partial{\Omega}$},\\
\int_{\Sigma_i}\mathbf{v}\cdot \mathbf{n} \ dS=\al_i& i=1,2,\\
\mathbf{v}\rightarrow \mathbf{v}_{P_i} &\text {as $x_1^i\to\infty$ in $\Omega_i$ for $i=1,2$},
\end{cases}
\end{equation}
where $\mathbf{v}_{P_i}$ is the Hagen-Poiseuille flow, which satisfies
\begin{equation}
\begin{cases}
\mathbf{v}_{P_i}=v_{P_i}(z^i)\mathbf{e}_1,\\
\mu_0\Delta_i'v_{P_i}+\nabla_i'\cdot\big{(}\mu_1|\mathcal{D}(v_{P_i})|^{p-2}\mathcal{D}(v_{P_i})\big{)}=-C_i,& \text{in $\Sigma_i$}\\
v_{P_i}(z^i)=0, & \text{on $\partial\Sigma_i$}
\end{cases}
\end{equation}
here $v_{P_i}$ is a scalar function, $z^i=(x_2^i, ..., x_d^i)$, $\mathbf{e}_1=(1, 0, ..., 0)\in\Bbb R^d$, $\Delta_i'=\p_{x_2^i}^2+...+\p_{x_d^i}^2$
and $\nabla_i'=(\p_{x_2^i}, ..., \p_{x_d^i})$.

For the nonlinear equation systems in (1.8), O.A.Ladyzhenskaya \cite{La-3} and J.L.Lions \cite{Lions} proved
the existence of the solution $\mathbf{v}$ by the
monotone operator theory in a bounded domain when $p\geqslant\frac{3d}{d+2}$. This result has been improved
by some authors, in particular, in \cite{Fre-ma-st}, the same result is established for $p>\frac{2d}{d+2}$.
For noncompact boundaries, particularly, for piping-system,
G.P.Galdi \cite{G-0} proved that if $\mu_0>0$, $p>2$, and flux is small, then problem (1.8) together with (1.9)
has a unique weak solution $\mathbf{v}$. If $\mu_0 = 0$, by deriving some ¡°weighted¡± energy estimates,
E. Maru\v{s}i\'{c}-Paloka  in \cite{Mar-pa} established the existence and uniqueness of the weak solution
$\mathbf{v}$ to problem (1.8) with (1.9) when $p> 2$ and the flux is small. Since the approach in \cite{Mar-pa}
requires  a detailed information about the dependence of $v_{P_i}$ on the cross-sectional coordinates,
where an explicit background solution is known, it seems that the resulting proof
in \cite{Mar-pa} is only suitable for the case of a circular cross section.
For arbitrary large flux, motivated by Ladyzhenskaya and Solonnikov's results
in \cite{La-5}, the authors in \cite{Gil-mar}
prove the existence and uniqueness of solution to the Ladyzhenskaya-Solonnikov Problem I for the non-Newtonian fluids when $p> 2$ and $\mu_0 = 0$.
In this paper, we shall consider both Ladyzhenskaya-Solonnikov Problem I and Ladyzhenskaya-Solonnikov Problem II for the
non-Newtonian fluids,  and intend to establish some systematic results. Here we point out that the
restriction of  $p>2$ when $\mu_0=0$ is essentially required in the proof of \cite{Gil-mar}
(one can see the statements of lines 8-9 from below on pages 3874 in \cite{Gil-mar}: ``As far as we know, the Leray problem for $p<2$ (with small
fluxes) is an open problem"), meanwhile only the corresponding Ladyzhenskaya-Solonnikov problem I is considered in \cite{Gil-mar}.
We shall  study problem (1.8) together with (1.9) for $p>1$ and $\mu_0\geqslant0$
(when $1<p<2$, the condition $\mu_0>0$ will be needed). On the other hand,
for the corresponding Ladyzhenskaya-Solonnikov Problem II of (1.8)
(i.e., the outlets of $\Omega$ may be permitted to be unbounded), we
shall establish both the existence and uniqueness of the solution $\mathbf{v}$
for $\mu_0>0$ and $p>1$  or $\mu_0=0$ and $2< p \le 3-\f{2}{d}$ (hence $d=3$), especially,
when the sections of $\Omega$ are uniformly bounded,
the resulting conclusions also hold for $\mu_0=0$ and $p>2$ (here we point out that
this case has been solved in \cite{Gil-mar}).

Let us comment on the proofs of our results. For the case of $\mu_0=0$  in (1.8),
if one wants to directly deal with the nonlinear
term $div(|D(\mathbf{v})|^{p-2}D(\mathbf{v}))$ for $p>2$ and apply the integration by parts
for equation (1.8) multiplying the solution $\mathbf{v}$ to obtain a priori estimate of  $\mathbf{v}$,
then the regularities of $\mathbf{v}\in W^{2,l}$ and $\pi\in W^{1,l}$
for some positive number $l$ in bounded domains are required as pointed out in \cite{Gil-mar}.
However this regularity is not expected for the weak solutions
of (1.8) if $p\not=2$ as stated in \cite{Gil-mar} (see lines 15-16
of pages 3875). To overcome this kind of difficulty,
the authors in \cite{Gil-mar} studied the following truncated modified problem
\begin{equation}\label{0.1}
\begin{cases}
-$div$ \ \big{(}\frac{1}{T}\mathcal{D}(\mathbf{v}^{T})+\mu_1|\mathcal{D}(\mathbf{v}^{T})|^{p-2}\mathcal{D}(\mathbf{v}^{T})\big{)}
+\mathbf{v}^{T}\cdot\nabla\mathbf{v}^{T}+\nabla \pi^{T}=0 &\text{in $\Omega(T)$},\\
$div$ \ \mathbf{v}^{T}=0 &\text{in $\Omega(T)$},\\
\mathbf{v}^{T}=0 &\text{on $\partial{\Omega(T)}$},
\end{cases}
\end{equation}
where $\Omega(T)=\Omega_0\cup\{x\in \Omega: 0<x_1^1< T, 0<x_1^2<T\}$.
By deriving the uniform estimates of $\mathbf{v}^{T}$ under the key assumption of $p>2$
and applying a local version of the Minty trick, the authors  in \cite{Gil-mar}
proved the existence and uniqueness of solution $\mathbf{v}$ to the Ladyzhenskaya-Solonnikov Problem I of (1.8)
when $\mu_0=0$ and $p>2$.
We now state our ingredients for treating problem (1.8) in this paper.
At first, we consider the following truncated modified problem
instead of (1.8)
\begin{equation}
\label{0.a}
\begin{cases}
-$div$ \ (\mu_0\mcD(\mathbf{v}^{T})+\mu_1|\mcD(\mathbf{v}^{T})|^{p-2}\mcD(\mathbf{v}^{T}))
+\mathbf{v}^{T}\cdot\nabla\mathbf{v}^{T}+\nabla \pi^T=0&\text{in $\Omega(T)$},\\
$div$ \ \mathbf{v}^{T}=0 &\text{in $\Omega(T)$},\\
\mathbf{v}^{T}=0&\text{on $\partial{\Omega(T)}$}.
\end{cases}
\end{equation}
As in \cite{Gil-mar} and \cite{La-5}, we assume that the velocity $\mathbf{v}^{T}$ of \eqref{0.a}
has the form $\mathbf{u}^T+\mathbf{a}$,
where $\mathbf{u}^T$ is the new unknown with zero flux, and $\mathbf{a}$ is a specially constructed solenoidal field satisfying
$\int_{\Sigma_i(x^i_1)}\mathbf{a}\cdot\mathbf{n}\ dS=\al_i$ and admitting some other ``good" properties.
To obtain a priori estimates of $\mathbf{u}^T$, we have to control the nonlinear term $\mathbf{u}^T\cdot\nabla\mathbf{u}^T\cdot\mathbf{a}$.
If one only assumes that $\mathbf{a}$ is bounded as in \cite{Gil-mar}, then it follows from  Young inequality and Poincar\'{e} inequality
that only the following estimate for $p>2$ is obtained
\begin{align}\label{0.2}
&|\int_{\Omega_i(t)}\mathbf{u}^T\cdot\nabla\mathbf{u}^T\cdot\mathbf{a}dx|\leqslant \ve\int_{\Omega_i(t)}|\nabla\mathbf{u}^T|^pdx+c(\ve)\int_{\Omega_i(t)}|\mathbf{u}^T|^{p'}dx\no\\
&\leqslant \ve\int_{\Omega_i(t)}|\nabla\mathbf{u}^T|^pdx+c(\ve)\int_{\Omega_i(t)}|\nabla\mathbf{u}^T|^{p'}dx\no\\
&\leqslant  \ve\int_{\Omega_i(t)}|\nabla\mathbf{u}^T|^pdx+\int_{\Omega_i(t)}(c(\ve)+\ve|\nabla\mathbf{u}^T|^p)dx\no\\
&\leqslant  2\ve\int_{\Omega_i(t)}|\nabla\mathbf{u}^T|^pdx+c(\ve)t,
\end{align}
where $\Omega_i(t)=\{x\in\Omega_i: 0<x^i_1<t\}$ and $c(\ve)>0$ stands for a generic constant depending on $\ve>0$.
From \eqref{0.2}, the authors in \cite{Gil-mar} obtained the crucial uniform estimate of $\int_{\Omega_i(t)}|\nabla\mathbf{u}^T|^pdx$
for the solution $\mathbf{v}^T$ to problem (1.10).
To relax the restriction of power $p$ and get the uniform control for solution $\mathbf{v}^T$ of problem (1.11),
we need more properties of $\mathbf{a}$ and
other interesting observations.
Note that, for $\mu_0>0$,
the leading term is $\int_{\Omega_i(t)}|\mcD(\mathbf{u}^T)|^2dx$ in the energy estimate of $\mathbf{u}^T$
(see (4.4) in $\S 4$), if one can  find a field $\mathbf{a}$ such that
\begin{equation}\label{0.3}
|\int_{\Omega_i(t)}\mathbf{u}^T\cdot\nabla\mathbf{u}^T\cdot\mathbf{a}dx|\leqslant \ve\int_{\Omega_i(t)}|\nabla\mathbf{u}^T|^2dx,
\end{equation}
then $\int_{\Omega_i(t)}|\nabla\mathbf{u}^T|^2dx$ instead of $\int_{\Omega_i(t)}|\nabla\mathbf{u}^T|^pdx$ can be estimated;
while, for $\mu_0=0$ and $p>2$, the leading term is $\int_{\Omega_i(t)}|\mcD(\mathbf{u})|^pdx$
in the estimate of $\mathbf{u}^T$ (see (4.4) in $\S 4$), if one can construct a vector field $\mathbf{a}$ such that $|\mathbf{a}|\leqslant c|\Sigma_i(x^i_1)|^{-1}$ in $\Omega_i(t)$, then it follows from the Young inequality and
Poincar\'{e} inequality that
\begin{equation}
\begin{split}
&|\int_{\Omega_i(t)}\mathbf{u}\cdot\nabla\mathbf{u}\cdot\mathbf{a}dx|\\
&\leqslant \ve\int_{\Omega_i(t)}|\mathbf{u}|^p|\Sigma_i(x^i_1)|^{\f{p}{d-1}}dx
+\ve\int_{\Omega_i(t)}|\nabla\mathbf{u}|^pdx+c(\ve)\int_{\Omega_i(t)}|\Sigma_i(x^i_1)|^{\f{p(d-2)}{(p-2)(d-1)}}dx\\
&\leqslant 2\ve\int_{\Omega_i(t)}|\nabla\mathbf{u}|^pdx+c(\ve)\int_0^t|\Sigma_i(s)|^{1-\f{p(d-2)}{(p-2)(d-1)}}ds,
\end{split}
\end{equation}
which derives the uniform estimate of $\int_{\Omega_i(t)}|\nabla\mathbf{u}^T|^pdx$
if $\int_0^t|\Sigma_i(s)|^{1-\f{p(d-2)}{(p-2)(d-1)}}ds<\infty$.
Thanks to \cite{G-2} Lemma III.4.3 and \cite{Pileckas} Lemma 2-Lemma 3, the aforementioned
$\mathbf{a}$ in (1.13) and (1.14) can be found.
On the other hand, for Ladyzhenskaya-Solonnikov Problem I and II,
the condition $\int_{\Omega_i(t)}|\mcD(\mathbf{u})|^pdx\leqslant c\int_0^t|\Sigma_i(s)|^{1-\f{dp}{d-1}}ds$
should be required (see problem (2.4) and problem (2.5) in $\S 2$). Hence, by (1.14) we require
such an inequality
\begin{align}
\int_0^t|\Sigma_i(s)|^{1-\f{p(d-2)}{(p-2)(d-1)}}ds\leqslant c\int_0^t|\Sigma_i(s)|^{1-\f{dp}{d-1}}ds.
\end{align}
In the case of $\mu_0=0$, for Ladyzhenskaya-Solonnikov Problem I, (1.15) is automatically satisfied for
any $p>2$ since $\Sigma_i(x^i_1)$ is bounded, while for Ladyzhenskaya-Solonnikov Problem II, (1.15)
is satisfied only for $2<p\le 3-\f{2}{d}$ $(d=3)$. Based on the uniform estimates of $\mathbf{u}^T$,
inspired by \cite{G-2} and \cite{Wolf},
through choosing some suitable test functions and taking some delicate
analysis on the resulting nonlinear terms, we can
show $\mathbf{v}^T\rightarrow\mathbf{v}$ a.e. in any compact subset of $\Omega$ by establishing
the uniform interior estimates of solution $\mathbf{v}^T$ to (1.11). From this, together with
some methods introduced in \cite{La-5} for treating the Newtonian fluids
and involved analysis on the resulting nonlinear terms in non-Newtonian fluids, we eventually
complete the proofs on the existence and uniqueness of solution $\mathbf{v}$ to the related
Ladyzhenskaya-Solonnikov Problem I
and Ladyzhenskaya-Solonnikov Problem II of (1.8) under some suitable conditions.

Our paper is organized as follows. In $\S 2$, the detailed descriptions on the resulting Ladyzhenskaya-Solonnikov Problems
for the non-Newtonian flows are given. In $\S 3$, we present some preliminary conclusions
which will be applied to prove our main results in subsequent sections. In $\S 4$, we establish
the existence of the solutions $\mathbf{v}^{T}$ to the bounded truncated problem corresponding to
(1.8).
In $\S 5$, we study the interior regularity of solutions $\mathbf{v}^{T}$ obtained in $\S 4$.
Based on $\S 4$ and $\S 5$, we shall complete the proofs on  Ladyzhenskaya-Solonnikov Problem I
and Ladyzhenskaya-Solonnikov Problem II of (1.8) in  $\S 6$ and $\S 7$ respectively.

\section{Descriptions of Ladyzhenskaya-Solonnikov Problems for non-Newtonian fluids}
We focus on the following non-Newtonian fluid problem in the domain $\Omega$
with noncompact boundaries (see Figure 3 and  Figure 4 below):

\begin{equation}
\begin{cases}
-$div$ \ \big{(}\mu_0\mathcal{D}(\mathbf{v})+\mu_1|\mathcal{D}(\mathbf{v})|^{p-2}\mathcal{D}(\mathbf{v})\big{)}
+\mathbf{v}\cdot\nabla\mathbf{v}+\nabla \pi=0&\text{in $\Omega$},\\
$div$ \ \mathbf{v}=0&\text{in $\Omega$},\\
\mathbf{v}=0&\text{on $\partial{\Omega}$},\\
\int_{\Sigma_i(x_1^i)}\mathbf{v}\cdot \mathbf{n} \ dS=\alpha_i&\text{with $\sum\limits^N_{i=1}\alpha_i=0$,}
\end{cases}
\end{equation}
where
\begin{equation*}
\Omega=\Omega_0\cup(\bigcup\limits_{i=1}^{N}\Omega_i),
\end{equation*}
and
\begin{equation*}
\Omega_i=\{x\in\mathbb{R}^n: x^i_1>0, y^i=(x^i_2,..,x^i_d)\in\Sigma_i(x^i_1)\}.
\end{equation*}

\vskip 0.2 true cm \includegraphics[width=12cm,height=9cm]{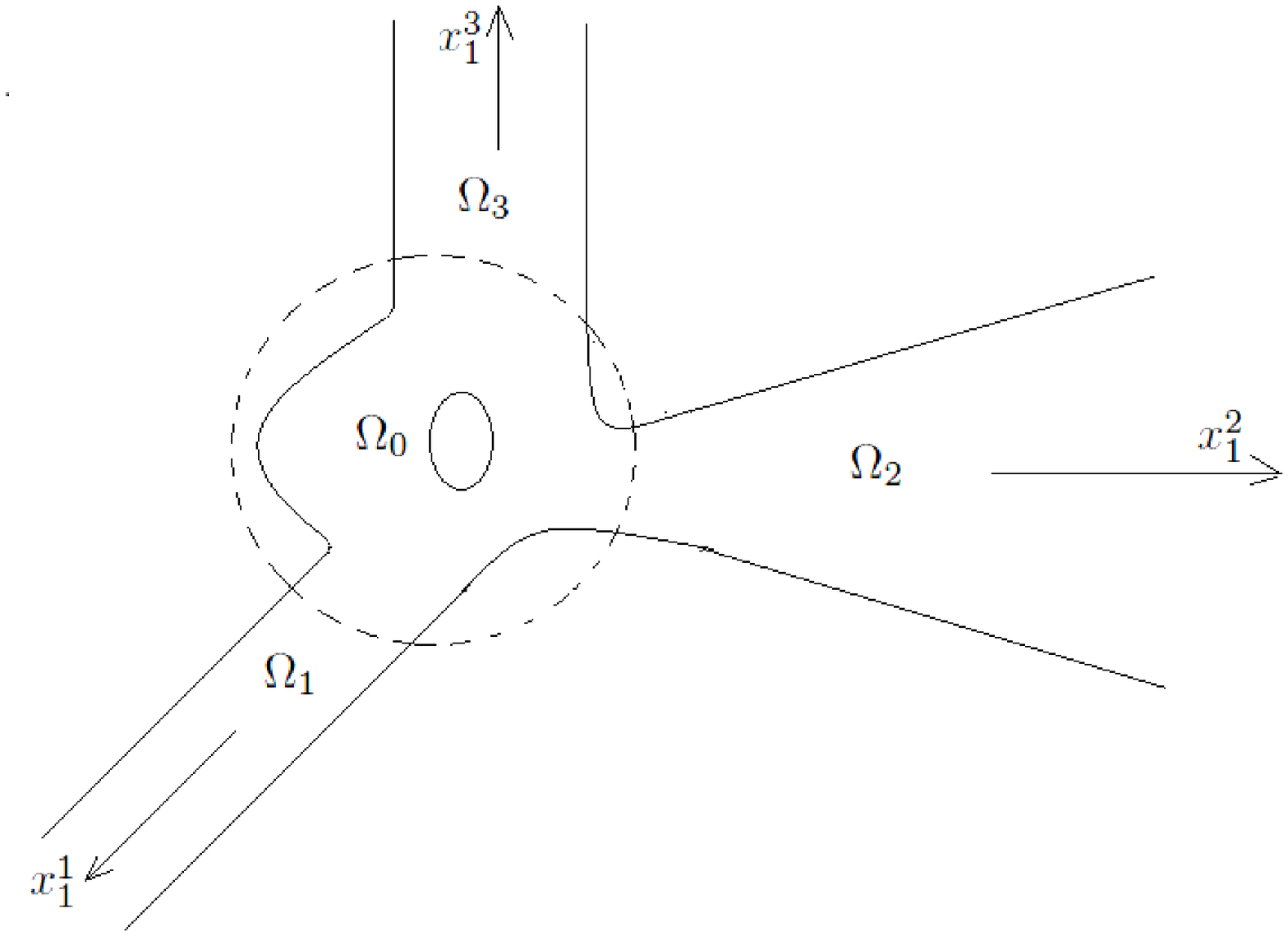}

\centerline{\bf Figure 3. Domain $\Omega$ for $d=2$.}

\vskip 0.2 true cm

Suppose that $\mathbf{v}$ is a solenoidal field and $\mathbf{v}=0$ holds on $\partial{\Omega_i}\backslash\{x^i=( x_1^i, y^i): x_1^i=0, y^i\in\Sigma_i(0)\}$.
Then
\begin{equation*}
|\alpha_i|^p=|\int_{\Sigma_i(t)}\mathbf{v}\cdot\mathbf{n} \ dS|^p\leqslant |\Sigma_i(t)|^{p-1}\int_{\Sigma_i(t)}
|\mathbf{v}|^p \ dS\leqslant c|\Sigma_i(t)|^{\frac{dp}{d-1}-1}\int_{\Sigma_i(t)}|\nabla'\mathbf{v}|^p \ dS,
\end{equation*}
where $c>0$ stands for a generic constant.
This means
\begin{equation*} |\al_i|^p\int_0^t|\Sigma_i(s)|^{1-\frac{dp}{d-1}}ds\leqslant c\int_{\Omega_i(t)}|\nabla\mathbf{v}|^pdx, \end{equation*}
where $\Omega_i(t)=\{x^i\in\Omega_i: 0<x_1^i<t\}$ for $1\le i\le N$.
Hence, if  $\al_i\neq0$  and
\begin{equation}
\text{$I_i(t)\equiv\int_0^t|\Sigma_i(s)|^{1-\frac{dp}{d-1}}ds\to +\infty$ as $t\to +\infty$},
\end{equation}
then
\begin{equation}
\text{$Q_i(t)\equiv\int_{\Omega_i(t)}|\nabla\mathbf{v}|^pdx\to +\infty$  as $t\to +\infty$}.
\end{equation}

\vskip 0.2 true cm
\includegraphics[width=15cm,height=14cm]{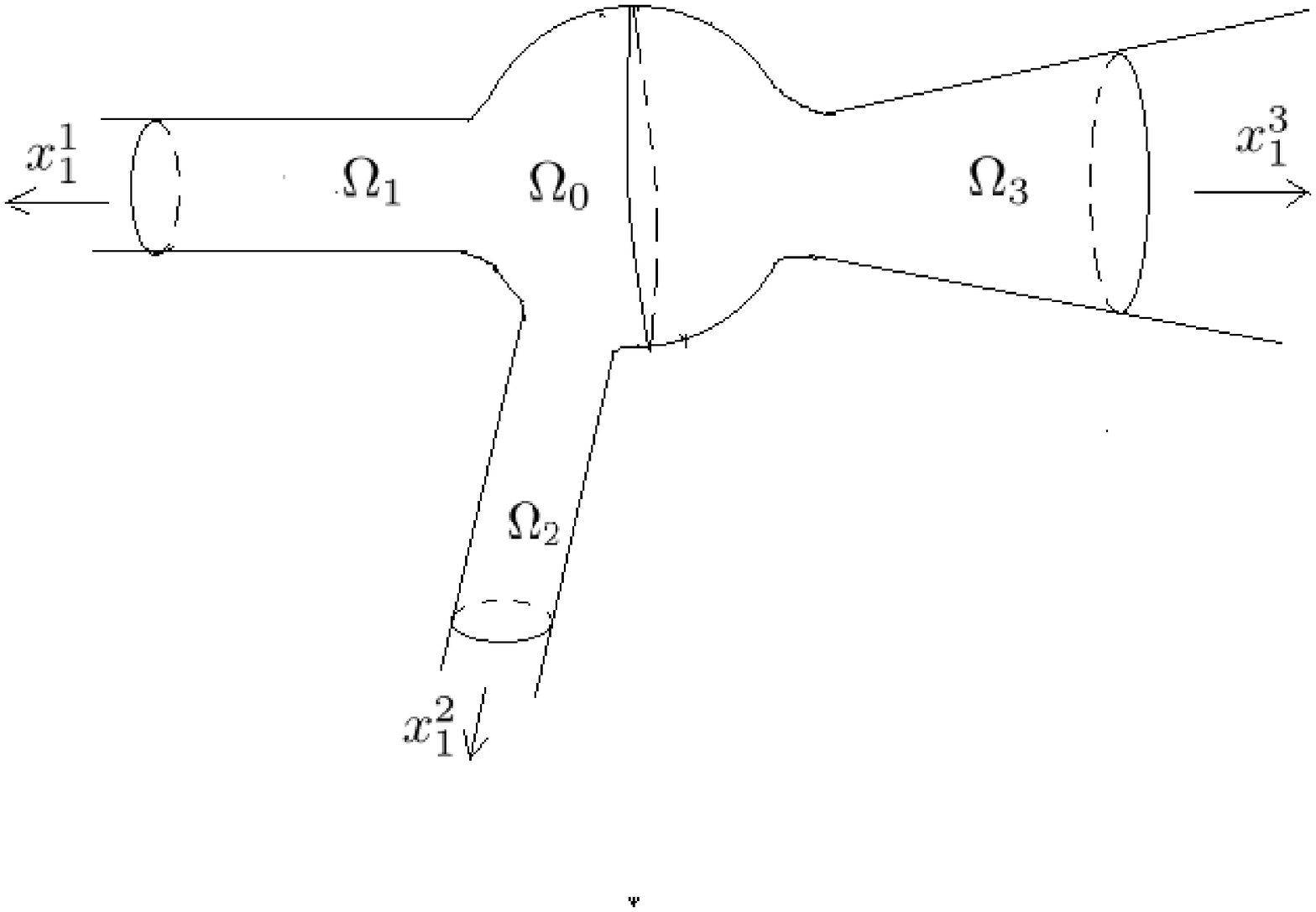}

\centerline{\bf Figure 4. Domain $\Omega$ for $d=3$.}
\vskip 0.2 true cm

From (2.2) and (2.3), it is natural to consider the following two problems

$\mathbf{Ladyzhenskaya-Solonnikov \ Problem \ I.}$  {\it Suppose that there are two positive constant $c_1$ and $c_2$
such that $c_1<|\Sigma_i(t)|<c_2$ \ for \ $1\le i\le N$.
We look for a pair vector field $(\mathbf{v},\pi)$ to fulfill
\begin{equation}
\label{problem 1}
\begin{cases}
-$div$ \ \big{(}\mu_0\mathcal{D}(\mathbf{v})+\mu_1|\mathcal{D}({\mathbf{v}})|^{p-2}\mathcal{D}(\mathbf{v})\big{)}
+\mathbf{v}\cdot\nabla\mathbf{v}+\nabla \pi=0&\text{in $\Omega$},\\
$div$ \ \mathbf{v}=0&\text{in $\Omega$}, \\
\mathbf{v}=0&\text{on $\partial{\Omega}$},\\
\int_{\Sigma_i(x_1^i)}\mathbf{v}\cdot \mathbf{n} \ dS=\al_i & \text{with $\ds\sum_{i=1}^N\al_i=0$},\\
\sup\limits_{t>0}t^{-1}Q_i(t)<\infty& \text{for $1\le i\le N$},
\end{cases}
\end{equation}
where $Q_i(t)$ is defined in (2.3) for $\mu_0=0$, and $Q_i(t)$ is defined as $Q_i(t)\equiv\int_{\Omega_i(t)}(|\nabla\mathbf{v}|^2+|\nabla\mathbf{v}|^p)dx$ for $\mu_0>0$.}

$\mathbf{Ladyzhenskaya-Solonnikov\ Problem \ II.}$ {\it Suppose that $I_i(\infty)=\infty$ for $1\le i\le m$,
while $I_i(\infty)<\infty$ for $m+1\le i\le N$, we look for a pair vector field $(\mathbf{v},\pi)$ such that
\begin{equation}
\label{problem 2}
\begin{cases}
-$div$ \ (\mu_0\mathcal{D}(\mathbf{v})+\mu_1|\mathcal{D}({\mathbf{v}})|^{p-2}\mathcal{D}(\mathbf{v}))
+\mathbf{v}\cdot\nabla\mathbf{v}+\nabla \pi=0&\text{in $\Omega$},\\
$div$ \ \mathbf{v}=0&\text{in $\Omega$},\\
\mathbf{v}=0&\text{on $\partial{\Omega}$},\\
\int_{\Sigma_i(x_1^i)}\mathbf{v}\cdot \mathbf{n} \ dS=\al_i &\text{with $\ds\sum_{i=1}^N\al_i=0$},\\
\sup\limits_{t>0}I_i^{-1}(t)Q_i(t)<\infty&\text{ for $1\le i\le N$},\end{cases}
\end{equation}
where $I_i(t)$ is defined in (2.2) for $\mu_0=0$, and $I_i(t)$ is defined as
$I_i(t)\equiv\int_0^t(|\Sigma_i(s)|^{-\frac{d+1}{d-1}}+|\Sigma_i(s)|^{1-\frac{dp}{d-1}})ds$ for $\mu_0>0$.}

In subsequent sections, we shall focus on the studies  on these two problems above.
The obtained results will be stated in Theorem 6.1-Theorem 6.2 and Theorem 7.3-Theorem 7.4 respectively.
In addition, for notational convenience, we introduce some function spaces as follows:
\begin{align*}
&D^{1,p}(\Omega)=\{\mathbf{u}\in L^1_{loc}(\Omega):\nabla\mathbf{u}\in L^p(\Omega)\},\\
&D^{1,p}_0(\Omega)= \{\text{\ completion\ of\ $C_0^{\infty}(\Omega)$\ in\ the\ semi-norm\ $\|\nabla\mathbf{u}\|_{p,\Omega}
\equiv(\int_{\Omega}|\nabla\mathbf{u}|^pdx)^{\f{1}{p}}$}\},\\
&D^{-1,p'}(\Omega)=(D^{1,p}(\Omega))',\quad  D_0^{-1,p'}(\Omega)=(D_0^{1,p}(\Omega))',\\
&\mathcal{D}(\Omega)=\{\mathbf{u}\in C_0^{\infty}(\Omega): \nabla\cdot\mathbf{u}=0\},\\
&\mathcal{D}^{1,p}_0(\Omega)=\{\text{completion\ of\ $\mathcal{D}(\Omega)$ \ in\ the\ semi-norm\ $\|\nabla\mathbf{u}\|_{p,\Omega}$}\}.
\end{align*}

\section{Preliminary results}

In this part, some preliminary results will be listed so that we can apply them
to study the described problems in $\S 2$. It follows from the proof of Appendix of \cite{La-5} that we have

{\bf Lemma 3.1.}
{\it Let $\omega=\{x=(x_1, y)\in\Bbb R^d: t_1<x_1<t_2, y\in\Sigma(x_1)\}$ and $\mathbf{u}|_{\{x: t_1<x_1<t_2, y\in\partial\Sigma(x_1)\}}=0$. Then
\begin{equation*}
\|\mathbf{u}\|_{r,\omega}\leqslant c(d,r,q)\max\{1,(t_2-t_1)^{-\f{1}{d}}\max_{x_1\in[t_1,t_2]}{|\Sigma(x_1)|^{\f{1}{d(d-1)}}}\}
|\omega|^{\f{1}{d}+\f{1}{r}-\f{1}{q}}\|\nabla\mathbf{u}\|_{q,\omega},
\end{equation*}
where $1\leqslant r\leqslant\f{d q}{d-q}$, if $q<d$; $1\leqslant r<\infty$, if $q=d$; $1\leqslant r\leqslant\infty$, if $q>d$.}

\vskip 0.1 true cm

The following result will play a crucial role in estimating pressure $\pi$ in problems (2.4)-(2.5).

\vskip 0.1 true cm

{\bf Lemma 3.2. (Theorem III.3.3 of \cite{G-2})} {\it For a bounded Lipschitzian domain $\omega\subset \Bbb R^d$, suppose that
\begin{equation}
\text{$f\in L^{p}(\omega)$ and $\int_{\omega}f dx=0$,}
\end{equation}
where $1<p<\infty$. Then one can find a vector field $\mathbf{w}$ such that
\begin{equation}
\begin{cases}
\nabla\cdot\mathbf{w}=0,\\
\mathbf{w}\in W^{1,p}_0(\omega),\\
\|\mathbf{w}\|_{1,p}\leqslant M(\omega)\|f\|_p,
\end{cases}
\end{equation}
where $M(\omega)>0$ is a constant depending only on the Lebesgue's measure $|\omega|$ of domain $\omega$.}

\vskip 0.1 true cm

{\bf Remark 3.1.} {\it Assume that $\omega$ is a star-shaped domain
with respect to a ball B with the radius $R_0$.
Then it follows from Theorem III.3.1 of \cite{G-2} that the positive constant $M(\omega)$ in
(3.2) satisfies $M(\omega) \leqslant c(d,p)(\frac{diam{(\omega)}}{R_0})^d(1+\frac{diam{(\omega)}}{R_0})$.
This property will be useful in order to solve problem (2.5).}

\vskip 0.1 true cm

{\bf Remark 3.2.} {\it  If $f\in L^{p}(\omega)\cap L^{r}(\omega)$ with $1<p,r<\infty$, then
one can find a vector field $\mathbf{w}\in W^{1,p}(\omega)\cap W^{1,r}(\omega)$ (see Remark III.3.12 of \cite{G-2})
such that
\begin{equation*}
\text{$\|\mathbf{w}\|_{1,p}\leqslant M(\omega)\|f\|_p$ \quad and \quad $\|\mathbf{w}\|_{1,r}\leqslant M(\omega)\|f\|_r$.}
\end{equation*}  }

Next, we list some results, whose proofs can be
found in Lemma 2.3 of \cite{La-5} or Lemma 3.1 of \cite{Gil-mar}.

\vskip 0.1 true cm

{\bf  Lemma 3.3.}
{\it   Let $\dl$ be a fixed constant with $\delta\in (0,1)$ and $t_0<T$. In addition,
we suppose that
$\Psi(\tau)$ is a monotonically increasing function, equal to zero for $\tau=0$
and equal to infinity for $\tau=\infty$.

(i) Assume that the nondecreasing, nonnegative smooth functions z(t) and $\varphi(t)$,
not identically equal to zero, satisfy the following inequalities for all $t\in [t_0,T]$,
\begin{equation}
\label{saint-venant 1}
z(t)\leqslant \Psi(z'(t))+(1-\delta)\varphi(t),
\end{equation}
and
\begin{equation}
\label{saint-venant 2}
\varphi(t)\geqslant\delta^{-1}\Psi(\varphi'(t)).
\end{equation}
If
\begin{equation}
\label{saint-venant 3}
z(T)\leqslant\varphi(T),
\end{equation}
then for all $t\in[t_0,T]$,
\begin{equation}
\label{saint-venant 4}
z(t)\leqslant\varphi(t).
\end{equation}

(ii) Assume that inequalities \eqref{saint-venant 1} and \eqref{saint-venant 2} are fulfilled for all $t\geqslant t_0$.
Then \eqref{saint-venant 4} holds for $t\geqslant t_0$ if
\begin{equation}
\liminf\limits_{t\rightarrow\infty}\frac{z(t)}{\varphi(t)}<1,
\end{equation}
or if z(t) has an order of growth for $t\rightarrow\infty$, less than the order of growth of the positive solutions to
the equation
\begin{equation}
\tilde{z}(t)=\delta^{-1}\Psi(\tilde{z}'(t)).
\end{equation}

(iii) Assume that the nonidentical zero nonnegative functions z(t), satisfying the homogenous inequality
\begin{equation}
\text{$z(t)\leqslant\delta^{-1}\Psi(z'(t))$ for $t\geqslant t_0$,}
\end{equation}
increases unboundedly for $t\rightarrow\infty$.
If $\delta^{-1}\Psi(\tau)\leqslant c_0\tau^m$ holds for $m>1$ and $\tau\geqslant\tau_1$, then \begin{equation}
\liminf\limits_{t\rightarrow\infty}t^{-\frac{m}{m-1}}z(t)>0;
\end{equation}
if, however, $\delta^{-1}\Psi(\tau)\leqslant c_0\tau$ holds for $\tau\geqslant\tau_1$, then \begin{equation}
\liminf\limits_{t\rightarrow\infty}z(t)exp(-\frac{t}{c_0})>0.
\end{equation}
}

The following Korn-type inequality can be referred  in
Theorem 3.2 of \cite{Kondratev-olenik} or Theorem 1 of \cite{G-1}.

{\bf Lemma 3.4.} {\it Let K be a cone in $\mathbb{R}^d$ and $p>1$. If $\int_{K}|\mathcal{D}(\mathbf{u})|^pdx<+\infty$,
then there is a skew-symmetric matrix A with constant coefficients such that
\begin{equation}
\int_K|\nabla(\mathbf{u}(x)-Ax)|^pdx\leqslant C\int_K |\mcD(\mathbf{u})(x)|^pdx, \end{equation}
where the positive constant C does not depend on the function $\mathbf{u}$ itself.
}

\vskip 0.1 true cm

{\bf Remark 3.3.} {\it If $\int_K|\nabla\mathbf{u}|^pdx<+\infty$, then $A=0$ holds in (3.12).}

\vskip 0.1 true cm

Finally, we state a conclusion as follows, whose
proof can be found in \cite{G-2} Lemma III.4.3, and \cite{Pileckas} Lemma 2-Lemma 3.

\vskip 0.1 true cm

{\bf Lemma 3.5.} {\it
Assume that the domain $\Omega$ and the numbers $\al_i$ $(1\le i\le N)$
are defined in problem (2.4) or problem  (2.5). Let $\al=
\max\limits_{1\le i\le N}|\al_i|$. Then for any fixed $\varepsilon>0$, there exists a smooth divergence-free
vector field $\mathbf{a}(x,\varepsilon)$ which vanish in a neighborhood of $\partial{\Omega}\cap\partial{\Omega_i}$
($1\le i\le N$), and which satisfies
\begin{align*}
 &(i) \ |\mathbf{a}|\leqslant c(\ve)\al|\Sigma_i(t)|^{-1} \ and \ |\nabla\mathbf{a}|\leqslant c(\varepsilon)\alpha|\Sigma_i(t)|^{-\frac{d}{d-1}}\quad\text{for $x\in\Omega_i(t)$ and  $1\le i\le N$}.\\
&(ii)\ \int_{\Sigma_i(t)}\mathbf{a}\cdot\mathbf{n}\ dS=\alpha_i\quad\text{for $1\le i\le N$}.\\
&(iii) \ \int_{\Omega_0}\mathbf{a}^2\mathbf{w}^2dx\leqslant\varepsilon\alpha^2\int_{\Omega_0}|\nabla\mathbf{w}|^2dx\quad\text{for any $\mathbf{w}\in\mathcal{D}(\Omega)$}.\\
&(iv) \ \int_{\Omega_i(t_2)\backslash\Omega_i(t_1)}\mathbf{a}^2\mathbf{w}^2dx
\leqslant\varepsilon\alpha^2\int_{\Omega_i(t_2)\backslash\Omega_i(t_1)}|\nabla\mathbf{w}|^2dx\quad
 \text{for any $\mathbf{w}\in\mathcal{D}(\Omega)$, $t_2>t_1>0$, and $1\le i\le N$}.
\end{align*}}

\section{Existence of solutions to problems (2.4) and  (2.5) in bounded truncated domains}
In this part, for the following problem in the bounded domain $\Omega(T)=\Omega_0\cup\{x\in\Omega: 0<x_1^1<T,..., 0<x_1^d<T\}$
\begin{equation}
\label{eq:truncated doamin}
\begin{cases}
-$div$ \ (\mu_0\mcD(\mathbf{v}^{T})+\mu_1|\mcD(\mathbf{v}^{T})|^{p-2}\mcD(\mathbf{v}^{T}))+\mathbf{v}^{T}\cdot\nabla\mathbf{v}^{T}+\nabla \pi^T=0&\text{in $\Omega(T)$},\\
$div$ \ \mathbf{v}^{T}=0&\text{in $\Omega(T)$}, \\
\mathbf{v}^{T}=0&\text{on $\partial{\Omega(T)}$},
\end{cases}
\end{equation}
we intend to find a weak solution $(\mathbf{v}^{T}, \pi^T)=(\mathbf{u}^{T}+\mathbf{a}, \pi^T)$ such that
\begin{multline}
\label{def weak sol:truncated domain}
\mu_0\big{(}\mathcal{D}(\mathbf{u}^T)+\mathcal{D}(\mathbf{a}),\mathcal{D}(\boldsymbol{\psi})\big{)}
+\mu_1\Big{(}|\mathcal{D}(\mathbf{u}^T)+\mathcal{D}(\mathbf{a})|^{p-2}\big{(}\mathcal{D}(\mathbf{u}^T)
+\mathcal{D}(\mathbf{a})\big{)},\mathcal{D}(\boldsymbol{\psi})\Big{)}\\
=(\mathbf{u}^T\cdot\nabla\boldsymbol{\psi},\mathbf{u}^T)+(\mathbf{u}^T\cdot\nabla\boldsymbol{\psi},\mathbf{a})
+(\mathbf{a}\cdot\nabla\boldsymbol{\psi},\mathbf{u}^T)+(\mathbf{a}\cdot\nabla\boldsymbol{\psi},\mathbf{a}),
\ \forall \ \boldsymbol{\psi}\in\mcD(\Omega(T)),
\end{multline}
where the vector value function $\mathbf{a}$ is given in Lemma 3.5.

{\bf Theorem 4.1.} {\it Let $\mu_0>0$ and $p>1$ or $\mu_0=0$ and $p>2$.
Then there is a vector field $\mathbf{u}^T$ such
that \eqref{def weak sol:truncated domain} holds, and $\mathbf{u}^T\in\mathcal{D}^{1,2}_0(\Omega(T))\bigcap\mathcal{D}^{1,p}_0(\Omega(T))$, if $\mu_0>0$;
$\mathbf{u}^T\in\mathcal{D}^{1,p}_0(\Omega(T))$, if $\mu_0=0$.
}

\vskip 0.2 true cm

{\bf Proof.}  Although the proof of Theorem 4.1 is standard as in \cite{G-0}-\cite{G-1} and \cite{Gil-mar},
where the authors treated problem (4.2) for different  vector value function $\mathbf{a}$,
we still give out the detailed proof for the sake of completeness.\\

\vskip 0.1 true cm

$\mathbf{Case\ I}$. $\mu_0>0$, $p>1$ \\

\vskip 0.1 true cm

Let \{$\boldsymbol{\psi}_k^T\}$ be a basis in $\mathcal{D}^{1,2}_0(\Omega(T))$.
We look for a series $\{c_{km}^T\}$ such that $\mathbf{u}_m^{T}=\ds\sum_{i=1}^m c^{T}_{km}\boldsymbol{\psi}_k^T$
satisfies
\begin{multline}
\mu_0\big{(}\mathcal{D}(\mathbf{u}_m^T)+\mathcal{D}(\mathbf{a}),\mathcal{D}(\boldsymbol{\psi}_k^T)\big{)}+\mu_1\Big{(}|\mathcal{D}(\mathbf{u}_m^T)
+\mathcal{D}(\mathbf{a})|^{p-2}\big{(}\mathcal{D}(\mathbf{u}_m^T)+\mathcal{D}(\mathbf{a})\big{)},\mathcal{D}(\boldsymbol{\psi}_k^T)\Big{)}\\
=(\mathbf{u}_m^T\cdot\nabla\boldsymbol{\psi}_k^T,\mathbf{u}_m^T)+(\mathbf{u}_m^T\cdot\nabla\boldsymbol{\psi}_k^T,\mathbf{a})
+(\mathbf{a}\cdot\nabla\boldsymbol{\psi}_k^T,\mathbf{u}_m^T)+(\mathbf{a}\cdot\nabla\boldsymbol{\psi}_k^T,\mathbf{a})
\qquad\text{ for $k=1,...,m$.}
\end{multline}
Multiplying both sides of (4.3) by $c_{km}^T$
and summing over $k$ yield
\begin{equation}
\begin{split}
\label{weak form:Galerkin solution}
&\mu_0\|\mathcal{D}\mathbf{u}_m^T\|^2_2+\mu_0\big{(}\mathcal{D}(\mathbf{a}),\mathcal{D}(\mathbf{u}_m^T)\big{)}+\mu_1\Big{(}|\mathcal{D}(\mathbf{u}_m^T)
+\mathcal{D}(\mathbf{a})|^{p-2}\big{(}\mathcal{D}(\mathbf{u}_m^T)+\mathcal{D}(\mathbf{a})\big{)},\mathcal{D}(\mathbf{u}_m^T)+\mathcal{D}(\mathbf{a})\Big{)}\\
&\quad -\mu_1\Big{(}|\mathcal{D}(\mathbf{u}_m^T)+\mathcal{D}(\mathbf{a})|^{p-2}\big{(}\mathcal{D}(\mathbf{u}_m^T)
+\mathcal{D}(\mathbf{a})\big{)},\mathcal{D}(\mathbf{a})\Big{)}\\
&=(\mathbf{u}_m^T\cdot\nabla\mathbf{u}_m^T,\mathbf{a})+(\mathbf{a}\cdot\nabla\mathbf{u}_m^T,\mathbf{a}).
\end{split}
\end{equation}
Using Schwarz inequality we get
\begin{equation}
\mu_0\big{(}\mathcal{D}(\mathbf{a}),\mathcal{D}(\mathbf{u}_m^T)\big{)}
\geqslant-\frac{\mu_0}{2}\|\mathcal{D}(\mathbf{u}_m^T)\|^2_2-\frac{\mu_0}{2}\|\mathcal{D}(\mathbf{a})\|^2_2.
\end{equation}
In addition,
\begin{equation}
\|\mathcal{D}(\mathbf{u}^T_m)+\mathcal{D}(\mathbf{a})\|^p_p\geqslant\frac{1}{2^{p-1}}\|\mathcal{D}(\mathbf{u}^T_m)\|^p_p
-\|\mathcal{D}(\mathbf{a})\|^p_p.
\end{equation}
We also notice that, by H\"{o}lder inequality and Young inequality,
\begin{equation}
\begin{split}
&\Big{|}\mu_1\Big{(}|\mathcal{D}(\mathbf{u}_m^T)+\mathcal{D}(\mathbf{a})|^{p-2}\big{(}\mathcal{D}(\mathbf{u}_m^T)
+\mathcal{D}(\mathbf{a})\big{)},\mathcal{D}(\mathbf{a})\Big{)}\Big{|}\\
&\leqslant \frac{\mu_1}{2^{p}}\|\mathcal{D}(\mathbf{u}^T_m)\|^p_p +c\|\mathcal{D}(\mathbf{a})\|^p_p,
\end{split}
\end{equation}
where and below $c>0$ denotes by a generic positive constant.
By Lemma 3.5 (iii) and (iv) we have that for any fixed $\ve>0$,
\begin{equation}
\Big{|}\int_{\Omega(T)}\mathbf{u}_m^T\cdot\nabla\mathbf{u}_m^T\cdot\mathbf{a}dx\Big{|}\leqslant
c\int_{\Omega(T)}\mathbf{a}^2|\mathbf{u}^T_m|^2dx+\frac{\varepsilon}{2}\int_{\Omega(T)}|\nabla\mathbf{u}^T_m|^2dx\leqslant
\varepsilon \int_{\Omega(T)}|\nabla\mathbf{u}_m^T|^2dx.
\end{equation}
On the other hand, by H\"{o}lder inequality and Korn inequality,
\begin{equation}
|(\mathbf{a}\cdot\nabla\mathbf{u}_m^T,\mathbf{a})|\leqslant
c\|\mathcal{D}(\mathbf{u}_m^T)\|_2 \|\mathbf{a}\|^2_{4} \leqslant \frac{\mu_0}{4} \|\mathcal{D}(\mathbf{u}_m^T)\|^2_2+c\|\mathbf{a}\|_{4}^{4}.
\end{equation}
Set $\mathbf{c}_m^T=(c^T_{1m},...,c^T_{mm})\in\R^m$.
To obtain a solution of (4.3),
we define a function $P: \R^m\rightarrow\R^m$, whose components are
\[
\begin{split}
P_k(\mathbf{c}_m^T)=\mu_0\big{(}\mathcal{D}(\mathbf{u}_m^T)+\mcD(\mathbf{a}),\mathcal{D}(\boldsymbol{\psi}_k^T)\big{)}+\mu_1\Big{(}|\mathcal{D}(\mathbf{u}_m^T)
+\mathcal{D}(\mathbf{a})|^{p-2}\big{(}\mathcal{D}(\mathbf{u}_m^T)+\mathcal{D}(\mathbf{a})\big{)},\mathcal{D}(\boldsymbol{\psi}_k^T)\Big{)}\\
\qquad -(\mathbf{u}_m^T\cdot\nabla\boldsymbol{\psi}_k^T,\mathbf{u}_m^T)+(\mathbf{u}_m^T\cdot\nabla\boldsymbol{\psi}_k^T,\mathbf{a})
+(\mathbf{a}\cdot\nabla\boldsymbol{\psi}_k^T,\mathbf{u}_m^T)+(\mathbf{a}\cdot\nabla\boldsymbol{\psi}_k^T,\mathbf{a}),
\qquad\text{$k=1,...,m$.}
\end{split}
\]
It follows from (4.5)-(4.9) and Korn inequality that
\[
\begin{split}
P(\mathbf{c}^T_m)\cdot \mathbf{c}^T_m&\geqslant\frac{\mu_0}{2}\|\mcD(\mathbf{u}^T_m)\|^2_2
-\frac{\mu_0}{2}\|\mcD(\mathbf{a})\|^2_2+\frac{\mu_1}{2^{p-1}}\|\mcD(\mathbf{u}^T_m)\|^p_p-\mu_1\|\mathcal{D}(\mathbf{a})\|^p_p\\
&\quad -\frac{\mu_1}{2^{p}}\|\mcD(\mathbf{u}^T_m)\|^p_p-c\|\mcD(\mathbf{a})\|^p_p-\ve\|\nabla{u}^T_m\|^2_2
-\frac{\mu_0}{4}\|\mcD(\mathbf{u}^T_m)\|^2_2-c\|\mathbf{a}\|^{4}_{4}\\
&\geqslant\frac{\mu_0}{8}\|\mcD(\mathbf{u}^T_m)\|^2_2+\frac{\mu_1}{2^p}\|\mcD(\mathbf{u}^T_m)\|^p_p
-c\|\mcD(\mathbf{a})\|^2_2-c\|\mcD(\mathbf{a})\|^p_p-c\|\mathbf{a}\|^{4}_{4}>0
\end{split}
\]
provided that $\|\mcD(\mathbf{u}^T_m)\|_2$ is large enough. This,
together with Lemma I.4.3 of \cite{Lions}, yields that there exists $\bar{\mathbf{c}}^T_m\in\R^m$ such that $P(\bar{\mathbf{c}}^T_m)=0$.
Hence we find a solution of Equation (4.3) for any fixed $m\in\N$.
Moreover, by (4.4)-(4.9)  we get
\[
\begin{split}
&\frac{\mu_0}{2}\|\mcD(\mathbf{u}^T_m)\|^2_2+\frac{\mu_1}{2^{p-1}}\|\mcD(\mathbf{u}^T_m)\|^p_p\\
&\leqslant\frac{\mu_0}{2}\|\mcD(\mathbf{a})\|^2_2+\mu_1\|\mcD(\mathbf{a})\|^p_p
+\frac{\mu_1}{2^p}\|\mcD(\mathbf{u}^T_m)\|^p_p+c\|\mcD(\mathbf{a})\|^p_p+\ve\|\nabla\mathbf{u}^T_m\|^2_2\\
&\quad +\frac{\mu_0}{4}\|\mcD(\mathbf{u}^T_m)\|^2_2+c\|\mathbf{a}\|^{4}_{4},
\end{split}
\]
and by Korn inequality,
\begin{equation}
\|\nabla\mathbf{u}^T_m\|^2_2+\|\nabla\mathbf{u}^T_m\|^p_p\leqslant c(\|\nabla\mathbf{a}\|^2_2+\|\nabla\mathbf{a}\|^p_p+\|\mathbf{a}\|^{4}_{4}).
\end{equation}
From (4.10), we obtain that there is a vector field $\mathbf{u}^T$  and a subsequence of $\{\mathbf{u}^T_m\}$, which is still denoted by $\{\mathbf{u}^T_m\}$,
such that
\begin{equation}
\begin{split}
\text{$\mathbf{u}_m^T\rightharpoonup\mathbf{u}^T$ in  $\mcD^{1,2}_0(\Omega(T))$},\\
\text{$\mathbf{u}_m^T\rightharpoonup\mathbf{u}^T$ in  $\mcD^{1,p}_0(\Omega(T))$}
\end{split}
\end{equation}
and
\begin{equation}
\text{$\mathbf{u}_m^T\rightarrow\mathbf{u}^T$  in  $L^2(\Omega(T))$}.
\end{equation}
Meanwhile, by (4.10)
\[
\begin{split}
&\||\mcD(\mathbf{u}^T_m)+\mcD(\mathbf{a})|^{p-2}\big{(}\mcD(\mathbf{u}^T_m)+\mathcal{D}(\mathbf{a})\big{)}
\|^{p'}_{p'}\\
&= \|\mcD(\mathbf{u}^T_m)+\mcD(\mathbf{a})\|^p_p\\
&\leqslant 2^{p-1}(\|\mcD(\mathbf{u}^T_m)\|^p_p+\|\mcD(\mathbf{u}^T_m)\|^p_p)\\
&\leqslant c(\|\nabla\mathbf{a}\|^2_2+\|\nabla\mathbf{a}\|^p_p+\|\mathbf{a}\|^{4}_{4}),
\end{split}
\]
which means that one can find a vector function $\mathbf{G}^T\in L^{p'}(\Omega(T))$ such that
\begin{equation}
\text{$\mu_1|\mcD(\mathbf{u}^T_m)+\mcD(\mathbf{a})|^{p-2}\big{(}\mcD(\mathbf{u}^T_m)+\mcD(\mathbf{a})\big{)}
\rightharpoonup \mathbf{G}^T$ in  $L^{p'}(\Omega(T))$}.
\end{equation}
Thus, by (4.11)-(4.13) we arrive at
\begin{multline}
\label{truncated domain:equility 1}
\mu_0\big{(}\mcD(\mathbf{u}^T)+\mcD(\mathbf{a}),\mcD(\boldsymbol{\psi}_k^T)\big{)}+\big{(}\mathbf{G}^T,\mcD(\boldsymbol{\psi}_k^T)\big{)}\\
=(\mathbf{u}^T\cdot\nabla\boldsymbol{\psi}_k^T,\mathbf{u}^T)+(\mathbf{u}^T\cdot\nabla\boldsymbol{\psi}_k^T,
\mathbf{a})+(\mathbf{a}\cdot\nabla\boldsymbol{\psi}_k^T,\mathbf{u}^T)
+(\mathbf{a}\cdot\nabla\boldsymbol{\psi}_k^T,\bf{a}).
\end{multline}
To prove $\mathbf{u}^T$ is a weak solution of (4.2), one should establish that for any $\boldsymbol{\vp}\in\mcD(\Omega(T))$,
\begin{equation}
\big{(}\mathbf{G}^T,\mathcal{D}(\boldsymbol{\varphi})\big{)}=\Big{(}\mu_1|\mathcal{D}(\mathbf{u})
+\mcD(\mathbf{a})|^{p-2}\big{(}\mcD(\mathbf{u})+\mcD(\mathbf{a})\big{)},\mcD(\boldsymbol{\varphi})\Big{)}.
\end{equation}
In fact, multiplying both sides of \eqref{truncated domain:equility 1} with $c_{km}^T$ and summing over $k$ yield
\begin{multline}
\mu_0\big{(}\mcD(\mathbf{u}^T)+\mcD(\mathbf{a}),\mcD(\mathbf{u}_m^T)\big{)}+\big{(}\mathbf{G}^T,\mcD(\mathbf{u}_m^T)\big{)}\\
=(\mathbf{u}^T\cdot\nabla\mathbf{u}_m^T,\mathbf{u}^T)+(\mathbf{u}^T\cdot\nabla\mathbf{u}_m^T,\mathbf{a})
+(\mathbf{a}\cdot\nabla\mathbf{u}_m^T,\mathbf{u}^T)
+(\mathbf{a}\cdot\nabla\mathbf{u}_m^T,\mathbf{a}).
\end{multline}
Subtracting (4.16) by  (4.4) and then passing to limit as $m\rightarrow\infty$, we get
\begin{equation}
\label{truncated domain:equility 3}
\lim_{m\rightarrow\infty}\left(\mu_0\mcD(\mathbf{u}_m^T)+\mcS\left(\mcD\left(\mathbf{u}_m^T\right)
+\mcD\left(\mathbf{a}\right)\right),\mcD(\mathbf{u}_m^T)\right)
=\mu_0\|\mcD(\mathbf{u}^T)\|_2^2+\big{(}G^T,\mcD(\mathbf{u}^T)\big{)},
\end{equation}
where and below $\mcS(D)\equiv\mu_1|D|^{p-2}D$ for the tensor $D$.
Since for any pair of tensors $D$ and $C$, we have the monotonicity property
\begin{equation*}
\big{(}\mcS(D)-\mcS(C)\big{)}\cdot(D-C)\geqslant0.
\end{equation*}
This yields that for any $\mathbf{\Phi}\in\mcD_0^{1,p}(\Omega(T))\cap\mcD_0^{1,2}(\Omega(T))$,
\begin{equation}
\Big{(}\mu_0\big{(}\mcD(\mathbf{u}^{T}_m)-\mcD(\mathbf{\Phi})\big{)}+\mcS\big{(}\mcD(\mathbf{u}_m^T)+\mcD(\mathbf{a})\big{)}-\mcS\big{(}\mcD({\mathbf{\Phi}})+\mcD(\mathbf{a})\big{)}
,\mcD(\mathbf{u}^T_m)-\mcD(\mathbf{\Phi})\Big{)}\geqslant0.
\end{equation}
Together with (4.17), we have that for $m\rightarrow\infty$,
\begin{equation*}
\Big{(}\mu_0\big{(}\mcD(\mathbf{u}^T)-\mcD(\mathbf{\Phi})\big{)}+G^T-\mcS\big{(}\mcD(\mathbf{\Phi})+\mcD(\mathbf{a})\big{)},\mcD(\mathbf{u}^T)-\mcD(\mathbf{\Phi})\Big{)}\geqslant0.
\end{equation*}
Choosing $\mathbf{\Phi}=\mathbf{u}^T-\ve \boldsymbol{\vp}$ with $\ve>0$ and $\boldsymbol{\vp}\in\mcD(\Omega(T))$.
Then
\begin{equation}
\Big{(}\ve\mu_0\mcD(\boldsymbol{\vp})+G^T-\mu_1|\mcD(\mathbf{u}^T)+\mcD(\mathbf{a})
-\ve\mcD(\boldsymbol{\vp})|^{p-2}\big{(}\mcD(\mathbf{u}^T)+\mcD(\mathbf{a})
-\ve\mcD(\boldsymbol{\vp})\big{)},\mcD(\boldsymbol{\vp})\Big{)}\geqslant0.
\end{equation}
Let $\varepsilon\rightarrow0$, we arrive at
\begin{equation}
\Big{(}\mathbf{G}^T-\mu_1|\mcD(\mathbf{u}^T)+\mcD(\mathbf{a})|^{p-2}\big{(}\mcD(\mathbf{u}^T)+\mcD(\mathbf{a})\big{)}
,\mcD(\boldsymbol{\vp})\Big{)}\geqslant0.
\end{equation}
If  $\boldsymbol{\varphi}$ is replaced by $-\boldsymbol{\varphi}$ in (4.20), then
\begin{equation}
\Big{(}\mathbf{G}^T-\mu_1|\mcD(\mathbf{u}^T)+\mcD(\mathbf{a})|^{p-2}\big{(}\mcD(\mathbf{u}^T)+\mcD(\mathbf{a})\big{)}
,\mcD(\boldsymbol{\vp})\Big{)}\leqslant0.
\end{equation}
Combining (4.18) with (4.19) yields
\begin{equation}
\big{(}\mathbf{G}^T,\mcD(\boldsymbol{\vp})\big{)}=\big{(}\mu_1|\mcD({\mathbf{u}}^T)+\mcD(\mathbf{a})|^{p-2}(\mcD(\mathbf{u}^T)+\mcD(\mathbf{a}))
,\mcD(\boldsymbol{\vp})\big{)}.
\end{equation}
Thus, $\mathbf{u}^T$ is a weak solution of (4.2).\\
\vskip 0.2 true cm

$\mathbf{Case\ II}$. $\mu_0=0$, $p>2$

\vskip 0.1 true cm
As in Case I, let \{$\boldsymbol{\psi}_k^T\}$ be a basis in $\mathcal{D}^{1,p}_0(\Omega(T))$
and set $\mathbf{u}_m^{T}=\ds\sum_{i=1}^m c^{T}_{km}\boldsymbol{\psi}_k^T$. Then for $k=1,...,m$,
\begin{multline}
\mu_1\Big{(}|\mathcal{D}(\mathbf{u}_m^T)+\mathcal{D}(\mathbf{a})|^{p-2}\big{(}\mathcal{D}(\mathbf{u}_m^T)
+\mathcal{D}(\mathbf{a})\big{)},\mathcal{D}(\boldsymbol{\psi}_k^T)\Big{)}\\
=(\mathbf{u}_m^T\cdot\nabla\boldsymbol{\psi}_k^T,\mathbf{u}_m^T)+(\mathbf{u}_m^T\cdot\nabla\boldsymbol{\psi}_k^T,\mathbf{a})
+(\mathbf{a}\cdot\nabla\boldsymbol{\psi}_k^T,\mathbf{u}_m^T)+(\mathbf{a}\cdot\nabla\boldsymbol{\psi}_k^T,\mathbf{a}).
\end{multline}
Multiplying both sides of (4.23) by $c_{km}^T$
and summing over $k$ yield
\begin{equation}
\begin{split}
&\mu_1\Big{(}|\mathcal{D}(\mathbf{u}_m^T)
+\mathcal{D}(\mathbf{a})|^{p-2}\big{(}\mathcal{D}(\mathbf{u}_m^T)+\mathcal{D}(\mathbf{a})\big{)},\mathcal{D}(\mathbf{u}_m^T)+\mathcal{D}(\mathbf{a})\Big{)}\\
&-\mu_1\Big{(}|\mathcal{D}(\mathbf{u}_m^T)+\mathcal{D}(\mathbf{a})|^{p-2}\big{(}\mathcal{D}(\mathbf{u}_m^T)
+\mathcal{D}(\mathbf{a})\big{)},\mathcal{D}(\mathbf{a})\Big{)}
=(\mathbf{u}_m^T\cdot\nabla\mathbf{u}_m^T,\mathbf{a})+(\mathbf{a}\cdot\nabla\mathbf{u}_m^T,\mathbf{a}).
\end{split}
\end{equation}
Note that
\begin{equation}
\|\mathcal{D}(\mathbf{u}^T_m)+\mathcal{D}(\mathbf{a})\|^p_p
\geqslant\frac{1}{2^{p-1}}\|\mathcal{D}(\mathbf{u}^T_m)\|^p_p-\|\mathcal{D}(\mathbf{a})\|^p_p
\end{equation}
and
\begin{equation}
\begin{split}
&\Big{|}\mu_1\Big{(}|\mathcal{D}(\mathbf{u}_m^T)+\mathcal{D}(\mathbf{a})|^{p-2}\big{(}\mathcal{D}(\mathbf{u}_m^T)
+\mathcal{D}(\mathbf{a})\big{)},\mathcal{D}(\mathbf{a})\Big{)}\Big{|}\\
&\leqslant \frac{\mu_1}{2^{p+1}}\|\mathcal{D}(\mathbf{u}^T_m)\|^p_p +c\|\mathcal{D}(\mathbf{a})\|^p_p.
\end{split}
\end{equation}
In addition, it follows from Lemma 3.5 (i), Poincar\'{e} inequality and Young's inequality that
\begin{align}
\Big{|}\int_{\Omega(T)}\mathbf{u}_m^T\cdot\nabla\mathbf{u}_m^T\cdot\mathbf{a}dx\Big{|}
&\leqslant\ve\int_{\Omega(T)}|\mathbf{a}|^{\f{p}{d-1}}|\mathbf{u}_m^T|^p+\varepsilon\int_{\Omega(T)}|\nabla\mathbf{u}_m^T|^pdx
+c\int_{\Omega(T)}|\mathbf{a}|^{\frac{p(d-2)}{(p-2)(d-1)}}dx\no\\
&\leqslant c\ve\|\nabla\mathbf{u}_m^T\|^p_p+c\|\mathbf{a}\|^{\frac{p(d-2)}{(p-2)(d-1)}}_{\frac{p(d-2)}{(p-2)(d-1)}}.
\end{align}
As in (4.9), we have
\begin{equation}
|(\mathbf{a}\cdot\nabla\mathbf{u}_m^T,\mathbf{a})|\leqslant
\|\mathcal{D}(\mathbf{u}_m^T)\|_p \|\mathbf{a}\|^2_{2p'} \leqslant \frac{\mu_1}{2^p} \|\mathcal{D}(\mathbf{u}_m^T)\|^p_p+c\|\mathbf{a}\|_{2p'}^{2p'}.
\end{equation}
Similarly to Case I, set $\mathbf{c}_m^T=(c^T_{1,m},...,c^T_{m,m})\in\R^m$ and define a function $P: \R^m\rightarrow\R^m$ as follows
\[
\begin{split}
&P_k(\mathbf{c}_m^T)=\mu_1\Big{(}|\mathcal{D}(\mathbf{u}_m^T)
+\mathcal{D}(\mathbf{a})|^{p-2}\big{(}\mathcal{D}(\mathbf{u}_m^T)+\mathcal{D}(\mathbf{a})\big{)},\mathcal{D}(\boldsymbol{\psi}_k^T)\Big{)}\\
&-(\mathbf{u}_m^T\cdot\nabla\boldsymbol{\psi}_k^T,\mathbf{u}_m^T)+(\mathbf{u}_m^T\cdot\nabla\boldsymbol{\psi}_k^T,\mathbf{a})
+(\mathbf{a}\cdot\nabla\boldsymbol{\psi}_k^T,\mathbf{u}_m^T)+(\mathbf{a}\cdot\nabla\boldsymbol{\psi}_k^T,\mathbf{a})
\qquad\text{ for $k=1,...,m$.}
\end{split}
\]
Then, by (4.25)-(4.28) and Korn inequality, we arrive at
\[
\begin{split}
P(\mathbf{c}^T_m)\cdot \mathbf{c}^T_m&\geqslant\frac{\mu_1}{2^{p-1}}\|\mcD(\mathbf{u}^T_m)\|^p_p-\mu_1\|\mathcal{D}(\mathbf{a})\|^p_p
-\frac{\mu_1}{2^{p+1}}\|\mcD(\mathbf{u}^T_m)\|^p_p\\
&\quad -c\|\mcD(\mathbf{a})\|^p_p-c\ve\|\nabla{u}^T_m\|^p_p-c\|\mathbf{a}\|^{\frac{p(d-2)}{(p-2)(d-1)}}_{\frac{p(d-2)}{(p-2)(d-1)}}
-\frac{\mu_1}{2^p}\|\mcD(\mathbf{u}^T_m)\|^p_p
-c\|\mathbf{a}\|^{2p'}_{2p'}\\
&\geqslant\frac{\mu_1}{2^{p+1}}\|\mcD(\mathbf{u}^T_m)\|^p_p-c\|\mcD(\mathbf{a})\|^p_p-c\|\mathbf{a}\|^{2p'}_{2p'}
-c\|\mathbf{a}\|^{\frac{p(d-2)}{(p-2)(d-1)}}_{\frac{p(d-2)}{(p-2)(d-1)}}>0
\end{split}
\]
for sufficiently large $\|\mcD(\mathbf{u}^T_m)\|_p$. From this, we then obtain the existence of solution
to Equation (4.3) for any fixed $m\in\N$.
Moreover, by (4.24)-(4.28)  we get
\[
\begin{split}
&\f{\mu_1}{2^{p-1}}\|\mcD(\mathbf{u}^T_m)\|^p_p\\
&\leqslant\frac{\mu_1}{2^{p+1}}\|\mcD(\mathbf{u}^T_m)\|^p_p+c\|\mcD(\mathbf{a})\|^p_p+c\ve\|\nabla\mathbf{u}^T_m\|^p_p\\
&\quad +c\|\mathbf{a}\|^{\frac{p(d-2)}{(p-2)(d-1)}}_{\frac{p(d-2)}{(p-2)(d-1)}}
+\frac{\mu_1}{2^p}\|\mcD(\mathbf{u}^T_m)\|^p_p+c\|\mathbf{a}\|^{2p'}_{2p'}
\end{split}
\]
and
\begin{equation}
\|\nabla\mathbf{u}^T_m\|^p_p\leqslant c(\|\nabla\mathbf{a}\|^p_p+\|\mathbf{a}\|^{\frac{p(d-2)}{(p-2)(d-1)}}_{\frac{p(d-2)}{(p-2)(d-1)}}+\|\mathbf{a}\|^{2p'}_{2p'}).
\end{equation}
Based on (4.29), we know that there is a vector field $\mathbf{u}^T$  and a subsequence of $\{\mathbf{u}^T_m\}$, which we still denote by $\{\mathbf{u}^T_m\}$,
such that
\begin{equation}
\begin{split}
\text{$\mathbf{u}_m^T\rightharpoonup\mathbf{u}^T$ in  $\mcD^{1,p}_0(\Omega(T))$}
\end{split}
\end{equation}
and
\begin{equation}
\text{$\mathbf{u}_m^T\rightarrow\mathbf{u}^T$  in  $L^p(\Omega(T))$}.
\end{equation}
Meanwhile, by (4.29),
\[
\begin{split}
&\||\mcD(\mathbf{u}^T_m)+\mcD(\mathbf{a})|^{p-2}\big{(}\mcD(\mathbf{u}^T_m)+\mathcal{D}(\mathbf{a})\big{)}
\|^{p'}_{p'}\\
&= \|\mcD(\mathbf{u}^T_m)+\mcD(\mathbf{a})\|^p_p
\leqslant 2^p(\|\mcD(\mathbf{u}^T_m)\|^p_p+\|\mcD(\mathbf{u}^T_m)\|^p_p)\\
&\leqslant  c(\|\nabla\mathbf{a}\|^p_p+\|\mathbf{a}\|^{\frac{p(d-2)}{(p-2)(d-1)}}_{\frac{p(d-2)}{(p-2)(d-1)}}+\|\mathbf{a}\|^{2p'}_{2p'}),
\end{split}
\]
which means that one can find $\mathbf{G}^T\in L^{p'}(\Omega(T))$ such that
\begin{equation}
\text{$\mu_1|\mcD(\mathbf{u}^T_m)+\mcD(\mathbf{a})|^{p-2}\big{(}\mcD(\mathbf{u}^T_m)+\mcD(\mathbf{a})\big{)}
\rightharpoonup \mathbf{G}^T$ in  $L^{p'}(\Omega(T))$}.
\end{equation}
Thus, by (4.30)-(4.32) we can obtain
\begin{equation}
\big{(}\mathbf{G}^T,\mcD(\boldsymbol{\psi}_k^T)\big{)}
=(\mathbf{u}^T\cdot\nabla\boldsymbol{\psi}_k^T,\mathbf{u}^T)+(\mathbf{u}^T\cdot\nabla\boldsymbol{\psi}_k^T,
\mathbf{a})+(\mathbf{a}\cdot\nabla\boldsymbol{\psi}_k^T,\mathbf{u}^T)
+(\mathbf{a}\cdot\nabla\boldsymbol{\psi}_k^T,\bf{a}).
\end{equation}
Completely analogous to the proof in Case I, one can prove that for any
$\boldsymbol{\vp}\in\mcD(\Omega(T))$,
\begin{equation}
\big{(}\mathbf{G}^T,\mcD(\boldsymbol{\vp})\big{)}=\big{(}\mu_1|\mcD({\mathbf{u}}^T)
+\mcD(\mathbf{a})|^{p-2}(\mcD(\mathbf{u}^T)+\mcD(\mathbf{a}))
,\mcD(\boldsymbol{\vp})\big{)}.
\end{equation}
Namely, $\mathbf{u}^T$ is a weak solution of (4.2).
\qquad \qquad \qquad \qquad \qquad \qquad \qquad \qquad \qquad \qquad  
\qquad \qquad $\square$
\section{Interior regularity of weak solutions}
In this part, we will establish the uniform interior estimates of weak solution $\mathbf{u}$ to the steady
non-Newtonian fluid equations in (2.4) or (2.5).

{\bf Theorem 5.1.}  {\it
Let $\Omega$ be any domain in $\mathbb{R}^d$, and $\Omega'\subset\subset\Omega$. Suppose that $\mathbf{v}=\mathbf{u}+\mathbf{a}$ is a weak solution to steady non-Newtonian fluid equations in $\Omega$, which satisfies for any $\boldsymbol{\psi}\in\mcD(\Omega)$,
\begin{multline}
\mu_0\big{(}\mcD(\mathbf{u})+\mcD(\mathbf{a}),\mcD(\boldsymbol{\psi})\big{)}+\mu_1\big{(}|\mcD(\mathbf{u})
+\mcD(\mathbf{a})|^{p-2}(\mcD(\mathbf{u})+\mcD(\mathbf{a})),\mcD(\boldsymbol{\psi})\big{)}\\
=\big{(}\mathbf{u}\cdot\nabla\boldsymbol{\psi},\mathbf{u}\big{)}+\big{(}\mathbf{u}\cdot\nabla\boldsymbol{\psi},\mathbf{a}\big{)}
+\big{(}\mathbf{a}\cdot\nabla\boldsymbol{\psi},\mathbf{u}\big{)}+\big{(}\mathbf{a}\cdot\nabla\boldsymbol{\psi},\mathbf{a}\big{)},
\end{multline}
where the vector value function $\mathbf{a}$ is given in Lemma 3.5.
Let $\Omega'_{4r}\equiv\{x: d(x, \Omega')<4r\}\subset\subset \Omega$ for any fixed number $r>0$.
Then we have that:

If $p>2$ and $\mu_0\geqslant0$, then $\nabla\mathbf{u}\in W^{\kappa,p}(\Omega')$ and
\begin{equation*}
\|\nabla\mathbf{u}\|_{\kappa,p,\Omega'}\leqslant c(\kappa,r,|\Omega'_{4r}|,\|\nabla\mathbf{u}\|_{p,\Omega'_{4r}}),
\end{equation*}
where $\kappa\in[0,\frac{2\hat{{\theta}}}{p})$ with $\hat{{\theta}}=\min\{p',\theta\}$, and $\theta=\frac{p(d+2)-3d}{2p}$, if $p<d$;
$\theta$ is an arbitrary constant less than 1,
if $d\leqslant p<3$ (only for $d=2$); $\theta=1$, if $ p\geqslant3$.

If $1<p<2$ and $\mu_0>0$, then $\nabla\mathbf{u}\in W^{\frac{2\theta}{p}-\ve,p}(\Omega')\cap W^{\theta-\ve,2}(\Omega')$
for any $\varepsilon>0$, and
\begin{equation*}
\|\nabla\mathbf{u}\|_{\theta-\ve,p,\Omega'}+\|\nabla\mathbf{u}\|_{\frac{2\theta}{p}-\varepsilon,p,\Omega'}\leqslant c(\ve,r,|\Omega'_{4r}|,\|\nabla\mathbf{u}\|_{2,\Omega'_{4r}}),
\end{equation*}
where $\theta=\f{1}{4}$, if $d=3$; $\theta$ is an arbitrary constant less than 1, if $d=2$.}

{\bf Proof.}  By $\Omega'\subset\subset\Omega$, then $dist(\partial\Omega',\partial\Omega)>0$.
For any $\rho$ with $0<\rho< dist(\partial\Omega',\partial\Omega)$, set
\begin{equation}
\Omega'_{\rho}=\{x\in\Omega:dist(x,\Omega')<\rho\}.
\end{equation}
Let $0<r<\frac{1}{4}dist(\partial\Omega',\partial\Omega)$ and choose a cutoff function
$\eta$ such that $\eta=1$ in $\Omega'_r$, $\eta=0$ in $\mathbb{R}^3\backslash\Omega'_{2r}$,
$0\leqslant\eta\leqslant1$ and $|\nabla\eta|<\frac{1}{r}$, $|\nabla^2\eta|<\frac{1}{r^2}$ in $\Omega'_{2r}$.
Next we study the following functional $F$: $D^{1,p}_0(\Omega'_{4r})\rightarrow \mathbb{R}$, where
\begin{align}
F(\bs{\psi})=&
\mu_0\big{(}\mcD(\mathbf{u})+\mcD(\mathbf{a}),\mcD(\boldsymbol{\psi})\big{)}_{\Omega'_{4r}}+\mu_1\big{(}|\mcD(\mathbf{u})
+\mcD(\mathbf{a})|^{p-2}(\mcD(\mathbf{u})+\mcD(\mathbf{a})),\mcD(\boldsymbol{\psi})\big{)}_{\Omega'_{4r}}\no\\
&-\big{(}\mathbf{u}\cdot\nabla\boldsymbol{\psi},\mathbf{u}\big{)}_{\Omega'_{4r}}-\big{(}\mathbf{u}\cdot\nabla\boldsymbol{\psi},\mathbf{a}\big{)}_{\Omega'_{4r}}
-\big{(}\mathbf{a}\cdot\nabla\boldsymbol{\psi},\mathbf{u}\big{)}_{\Omega'_{4r}}-\big{(}\mathbf{a}\cdot\nabla\boldsymbol{\psi},\mathbf{a}\big{)}_{\Omega'_{4r}}.
\end{align}

\vskip 0.1 true cm

$\mathbf{Case\ I}$. {\bf $\mu_0\geqslant0$, $p>2$}\\

\vskip 0.1 true cm

In this case, we just only treat the case of $\mu_0=0$
since the treatment for $\mu_0>0$ is easier. Since $\mathbf{u}$ satisfies (5.1), it follows that $ker{F}=\mathcal{D}_0^{1,p}(\Omega'_{4r})$.
Then according to De Rham Theorem,
we know that there exists a function $\pi \in L^{p'}(\Omega'_{4r})$ such that  for any $\boldsymbol{\psi}\in D^{1,p}_0(\Omega'_{4r})$,
\begin{equation}
F(\boldsymbol{\psi})=(\pi,\nabla\cdot\boldsymbol{\psi}).
\end{equation}
Without loss of generality, $\int_{\Omega'_{4r}}\pi dx=0$ can be assumed. By Lemma 3.2
we can find a vector filed $\widehat{\boldsymbol{\psi}}\in D^{1,p}_0(\Omega'_{4r})$ such that
$\nabla\cdot\widehat{\boldsymbol{\psi}}=|\pi|^{p'-2}\pi-\f{1}{|\Omega'_{4r}|}\int_{\Omega'_{4r}}|\pi|^{p'-2}\pi dx$
and
$\|\nabla\widehat{\boldsymbol{\psi}}\|_{p,\Omega'_{4r}}\leqslant c(|\Omega'_{4r}|)\|\pi\|^{p'-1}_{p',\Omega'_{4r}}$.
Combining (5.3) with (5.4), in terms of $\int_{\Omega'_{4r}}\pi dx=0$, we have that
\begin{equation*}
\|\pi\|^{p'}_{p',\Omega'_{4r}}=F(\widehat{\boldsymbol{\psi}})\leqslant c(|\Omega'_{4r}|)\|\pi\|^{p'-1}_{p',\Omega'_{4r}}\|F\|_{D^{-1,p'}_0({\Omega'_{4r}})}.
\end{equation*}
In addition,
\[
\begin{split}
\|F\|_{D^{-1,p'}_0(\Omega'_{4r})}&\leqslant\mu_1\|\mcD(\mathbf{u})+\mcD({\mathbf{a}})\|^{p-1}_{p,\Omega'_{4r}}+\|\mathbf{u}\|^2_{2p',\Omega'_{4r}}
+2\|\mathbf{a}\|_{\infty}\|\mathbf{u}\|_{p',\Omega'_{4r}}+\|\mathbf{a}\|^2_{2p',\Omega'_{4r}}\\
&\leqslant c(|\Omega'_{4r}|)(1+\|\nabla\mathbf{u}\|^{p-1}_{p,{\Omega'_{4r}}}+\|\nabla\mathbf{u}\|^2_{p,{\Omega'_{4r}}}+\|\nabla\mathbf{u}\|_{p,{\Omega'_{4r}}}).
\end{split}
\]
Hence we have
\begin{equation}
\|\pi\|_{p',\Omega'_{4r}}\leqslant c(|\Omega'_{4r}|)(1+\|\nabla\mathbf{u}\|^{p-1}_{p,{\Omega'_{4r}}}+\|\nabla\mathbf{u}\|^2_{p,{\Omega'_{4r}}}+\|\nabla\mathbf{u}\|_{p,{\Omega'_{4r}}}).
\end{equation}
Set $\Delta_{\la,k}\mathbf{u}=\mathbf{u}(x+\la \mathbf{e}_k)-\mathbf{u}(x)$. Let
$\Delta_{\la}\mathbf{u}$ stand for $\Delta_{\la,k}\mathbf{u}$ for any $k$ ($1\le k\le d$).
Choosing $\Delta_{-\lambda}(\eta^p\Delta_{\la}\mathbf{u})$ as a test function, we then get
that from (5.4),
\begin{equation}
\begin{split}
&\mu_1\int_{\Omega'_{3r}}|\mcD(\mathbf{u})+\mcD(\mathbf{a})|^{p-2}(\mcD(\mathbf{u})+\mcD(\mathbf{a})):\mathcal{D}(\Delta_{-\lambda}(\eta^p\Delta_{\lambda}\mathbf{u}))dx\\
&=-\int_{\Omega'_{3r}}\mathbf{u}\cdot\nabla\mathbf{u}\cdot(\Delta_{-\lambda}(\eta^p\Delta_{\lambda}\mathbf{u}))dx
+\int_{\Omega'_{3r}}\pi\nabla\cdot(\Delta_{-\lambda}(\eta^p\Delta_{\lambda}\mathbf{u}))dx\\
&\quad -\int_{\Omega'_{3r}}\mathbf{u}\cdot\nabla\mathbf{a}\cdot(\Delta_{-\lambda}(\eta^p\Delta_{\lambda}\mathbf{u}))dx
-\int_{\Omega'_{3r}}\mathbf{a}\cdot\nabla\mathbf{u}\cdot(\Delta_{-\lambda}(\eta^p\Delta_{\lambda}\mathbf{u}))dx\\
&\quad -\int_{\Omega'_{3r}}\mathbf{a}\cdot\nabla\mathbf{a}\cdot(\Delta_{-\lambda}(\eta^p\Delta_{\lambda}\mathbf{u}))dx.
\end{split}
\end{equation}
A direct computation yields
\begin{equation}
\begin{split}
&\mu_1\int_{\Omega'_{3r}}|\mcD(\mathbf{u})+\mcD(\mathbf{a})|^{p-2}(\mcD(\mathbf{u})
+\mcD(\mathbf{a})):\mcD(\Delta_{-\la}(\eta^p\Delta_{\lambda}\mathbf{u}))dx\\
=&p\mu_1\int_{\Omega'_{3r}}|\mcD(\mathbf{u})+\mcD(\mathbf{a})|^{p-2}(\mcD(\mathbf{u})
+\mcD(\mathbf{a})):\Delta_{-\la}Sym(\Delta_{\la}\mathbf{u}\otimes\eta^{p-1}\nabla\eta)dx\\
&+\mu_1\int_{\Omega'_{3r}}\Delta_{\la}(|\mcD(\mathbf{u})+\mcD(\mathbf{a})|^{p-2}(\mcD(\mathbf{u})
+\mcD(\mathbf{a}))):(\eta^p\Delta_{\lambda}\mcD(\mathbf{u}))dx\\
=&\mu_1\int_{\Omega'_{3r}}|\mcD(\mathbf{u})+\mcD(\mathbf{a})|^{p-2}(\mcD(\mathbf{u})
+\mcD(\mathbf{a})):\Delta_{-\la}Sym(\Delta_{\la}\mathbf{u}\otimes\eta^{p-1}\nabla\eta)dx\\
&+\mu_1\int_{\Omega'_{3r}}\eta^p\Delta_{\lambda}(|\mcD(\mathbf{u})
+\mcD(\mathbf{a})|^{p-2}(\mcD(\mathbf{u})+\mcD(\mathbf{a}))):\Delta_{\lambda}(\mcD(\mathbf{u})+\mcD(\mathbf{a}))dx\\
&-\mu_1\int_{\Omega'_{3r}}\eta^p \Delta_{\lambda}(|\mcD(\mathbf{u})|^{p-2}\mcD(\mathbf{u})):\Delta_{\lambda}(\mcD(\mathbf{a}))dx,
\end{split}
\end{equation}
where $Sym(D)\equiv D+D^{t}$ for a given second order tensor $D$.
Since
\begin{equation}
(|D|^{p-2}D-|C|^{p-2}C)\cdot(D-C)\geqslant \delta(|D|+|C|)^{p-2}|D-C|^2
\end{equation}
holds for any pair of tensors D and C, where $\delta>0$ is some suitable constant,
we have
\begin{equation}
\begin{split}
&\mu_1\int_{\Omega'_{3r}}\eta^p\Delta_{\lambda}(|\mcD(\mathbf{u})+\mcD(\mathbf{a})|^{p-2}(\mcD(\mathbf{u})
+\mcD(\mathbf{a}))):\Delta_{\lambda}(\mcD(\mathbf{u})+\mcD(\mathbf{a}))dx\\
&\geqslant \f{\delta}{2^{p-1}}\mu_1\int_{\Omega'_{3r}}\eta^p|\Delta_{\la}(\mcD(\mathbf{u}))|^p dx-\delta\mu_1\int_{\Omega'_{3r}}\eta^p|\Delta_{\la}(\mcD(\mathbf{a}))|^p dx.\\
\end{split}
\end{equation}
Collecting (5.6)-(5.9) yields
\begin{equation}
\begin{split}
&\f{\delta}{2^{p-1}}\mu_1\int_{\Omega'_{3r}}\eta^p|\Delta_{\la}(\mcD(\mathbf{u}))|^p dx-\delta\mu_1\int_{\Omega'_{3r}}\eta^p|\Delta_{\la}(\mcD(\mathbf{a}))|^p dx.\\
&\leqslant-p\mu_1\int_{\Omega'_{3r}}|\mcD(\mathbf{u})+\mcD(\mathbf{a})|^{p-2}(\mcD(\mathbf{u})
+\mcD(\mathbf{a})):\Delta_{-\la}Sym(\Delta_{\la}\mathbf{u}\otimes\eta^{p-1}\nabla\eta)dx\\
&\quad +\mu_1\int_{\Omega'_{3r}}\eta^p \Delta_{\lambda}(|\mcD(\mathbf{u})|^{p-2}\mcD(\mathbf{u})):\Delta_{\lambda}(\mcD(\mathbf{a}))dx\\
&\quad +\int_{\Omega'_{3r}}\pi\nabla\cdot(\Delta_{-\lambda}(\eta^p\Delta_{\lambda}\mathbf{u}))dx\\
&\quad -\int_{\Omega'_{3r}}\boldsymbol{\Phi}\cdot(\Delta_{-\lambda}(\eta^p\Delta_{\lambda}\mathbf{u}))dx\\
&\quad -\int_{\Omega'_{3r}}\mathbf{u}\cdot\nabla\mathbf{u}\cdot(\Delta_{-\lambda}(\eta^p\Delta_{\lambda}\mathbf{u}))dx\\
&\equiv I_1+...+I_5,
\end{split}
\end{equation}
where $\boldsymbol{\Phi}=-\mathbf{u}\cdot\nabla\mathbf{a}-\mathbf{a}\cdot\nabla\mathbf{u}-\mathbf{a}\cdot\nabla\mathbf{a}$.
We now focus on the treatments on $I_i$ $(1\le i\le 5)$ in (5.10).
At first, it is well known that $\|\Delta_{\lambda,k}u\|_{p,\Omega'}\leqslant|\lambda|\|\partial_ku\|_{p,\Omega}$
for $0<|\lambda|<dist(\Omega',\partial\Omega)$ (see Chapter 7 of \cite{G-T}).
For the term $I_1$, we have
\begin{align*}
|I_1|\leqslant p\mu_1\|\mathcal{D}(\mathbf{u})+\mathcal{D}(\mathbf{a})\|_{p,\Omega'_{3r}}^{p-1}
\|\Delta_{-\lambda}Sym(\Delta_{\lambda}\mathbf{u}\otimes\eta^{p-1}\nabla\eta)\|_{p,\Omega'_{3r}}.
\end{align*}
Note that
\begin{align*}
&\|\Delta_{-\lambda}Sym(\Delta_{\la}\mathbf{u}\otimes\eta^{p-1}\nabla\eta)\|^p_{p,\Omega'_{3r}}\\
&=\int_{\Omega'_{3r}}\big{|}\Delta_{-\la}[\eta^{p-1}(\partial_j\eta\Delta_{\lambda}u_i+\partial_i\eta\Delta_{\lambda}u_j)]\big{|}^pdx\\
&\leqslant 2^{p-1}|\la|^p\bigg(\int_{\Omega'_{3r}}(p-1)|\partial_k\eta(\partial_j\eta\Delta_{\lambda}u_i+\partial_i\eta\Delta_{\lambda}u_j)|^pdx\\
&\quad +\int_{\Omega'_{3r}}|\eta|^p|\partial_{jk}\eta\Delta_{\lambda}u_i+\p_j\eta\Delta_{\lambda}\partial_ku_i
+\partial_{ik}\eta\Delta_{\lambda}u_j+\p_i\eta\Delta_{\lambda}\partial_ku_j|^pdx\bigg)\\
&\leqslant c|\lambda|^p\bigg(\frac{|\la|^p}{r^2}\int_{\Omega'_{4r}}|\nabla\mathbf{u}|^pdx
+\frac{1}{r}\int_{\Omega'_{3r}}|\eta\nabla(\Delta_{\lambda}\mathbf{u})|^pdx\bigg).
\end{align*}
In addition, it follows from $
\eta\nabla(\Delta_{\lambda}\mathbf{u})=\nabla(\eta\Delta_{\lambda}\mathbf{u})-(\nabla\eta)\cdot\Delta_{\lambda}\mathbf{u}$
and Korn inequality that
\begin{align*}
&\bigg(\int_{\Omega'_{3r}}|\eta\nabla(\Delta_{\lambda}\mathbf{u})|^pdx\bigg)^{1/p}\\
&\leqslant c\bigg(\int_{\Omega'_{3r}}|\mathcal{D}(\eta(\Delta_{\lambda}\mathbf{u}))|^pdx\bigg)^{1/p}
+\frac{c}{r}\bigg(\int_{\Omega'_{3r}}|\Delta_{\lambda}\mathbf{u}|^pdx\bigg)^{1/p}\\
&\leqslant c\frac{|\lambda|}{r}\bigg(\int_{\Omega'_{4r}}|\nabla\mathbf{u}|^pdx\bigg)^{1/p}
+c\bigg(\int_{\Omega'_{3r}}|\eta\mathcal{D}(\Delta_{\lambda}\mathbf{u})|^pdx\bigg)^{1/p}.
\end{align*}
Hence, we arrive at
\begin{equation*}
\|\Delta_{-\lambda}Sym(\Delta_{\lambda}\mathbf{u}\otimes\eta^{p-1}\nabla\eta)\|_{p,\Omega'_{3r}}
\leqslant c\frac{\lambda^2}{r^2}\|\nabla\mathbf{u}\|_{p,\Omega'_{4r}}
+c\frac{|\lambda|}{r}\bigg(\int_{\Omega'_{4r}}|\eta\mathcal{D}(\Delta_{\lambda}\mathbf{u})|^pdx\bigg)^{1/p}.
\end{equation*}
Therefore,
\begin{equation}
\begin{split}
&|I_1|\leqslant c\|\mathcal{D}(\mathbf{u})+\mathcal{D}(\mathbf{a})\|_{p,\Omega'_{3r}}^{p-1}\left(\frac{\lambda^2}{r^2}\|\nabla\mathbf{u}\|_{p,\Omega'_{3r}}
+\frac{|\lambda|}{r}\left(\int_{\Omega'_{3r}}|\eta\mathcal{D}\left(\Delta_{\lambda}\mathbf{u}\right)|^pdx\right)^{1/p}\right)\\
&\leqslant c(r,|\Omega'_{4r}|,\|\nabla\mathbf{u}\|_{p,\Omega'_{4r}})|\la|^{p'}
+\ve\int_{\Omega'_{3r}}|\eta\mathcal{D}(\Delta_{\lambda}\mathbf{u})|^pdx.
\end{split}
\end{equation}
While
\begin{equation}
\begin{split}
|I_2|&\leqslant \mu_1\||\mcD(\mathbf{u})|^{p-2}\mcD(\mathbf{u})\|_{p',\Omega'_{3r}}\|\Delta_{-\la}(\eta^p\Delta_{\la}\mcD(\mathbf{a}))\|_{p,\Omega'_{3r}}\\
&\leqslant c(r,|\Omega'_{4r}|,\|\nabla\mathbf{u}\|_{p,\Omega'_{4r}})|\la|^2.
\end{split}
\end{equation}
On the other hand,
\begin{align*}
&\Big{|}\int_{\Omega'_{3r}}\pi\nabla\cdot(\Delta_{-\lambda}(\eta^p\Delta_{\lambda}\mathbf{u}))dx\Big{|}
=p\Big{|}\int_{\Omega'_{3r}}\pi\Delta_{-\lambda}(\eta^{p-1}\partial_i\eta\Delta_{\lambda}u_i)dx\Big{|}\\
&\leqslant p\|\pi\|_{p',\Omega'_{3r}}|\lambda|\|\partial_k(\eta^{p-1}\partial_i\eta\Delta_{\lambda}u_i)\|_{p,\Omega'_{3r}}\\
&\leqslant c\|\pi\|_{p',\Omega_{3r}}|\lambda|\bigg(\frac{1}{r^2}\|\Delta_{\lambda}u_i\|_{p,\Omega'_{3r}}
+\frac{1}{r}\|\eta\Delta_{\lambda}\nabla\mathbf{u}\|_{p,\Omega'_{3r}}\bigg)\\
&\leqslant \frac{c}{r}\|\pi\|_{p',\Omega'_{3r}}|\lambda|\bigg(\frac{1}{r}|\lambda|\|\nabla\mathbf{u}\|_{p,\Omega'_{4r}}
+\|\eta\Delta_{\lambda}\mathcal{D}\mathbf{u}\|_{p,\Omega'_{3r}}\bigg)\\
&\leqslant \frac{c}{r^2}\lambda^2\|\pi\|_{p',\Omega'_{3r}}\|\nabla\mathbf{u}\|_{p,\Omega'_{4r}}+\frac{c}{\varepsilon r^2}\lambda^{p'}\|\pi\|^{p'}_{p',\Omega'_{3r}}
+\varepsilon\|\eta\Delta_{\lambda}\mathcal{D}\mathbf{u}\|^p_{p,\Omega'_{3r}}.\\
\end{align*}
Then we arrive at
\begin{equation}
|I_3|\leqslant c(|\Omega'_{4r}|,r,\|\nabla\mathbf{u}\|_{p,\Omega'_{4r}})|\lambda|^{p'}+\ve\|\eta\Delta_{\lambda}\mathcal{D}\mathbf{u}\|^p_p.
\end{equation}
Since
\begin{equation*}
\|\bs\Phi\|_{p',\Omega'_{4r}}\leqslant\|\mathbf{u}\cdot\nabla\mathbf{a}\|_{p',\Omega'_{4r}}+\|\mathbf{a}\cdot\nabla\mathbf{u}\|_{p',\Omega'_{4r}}
+\|\mathbf{a}\cdot\nabla\mathbf{a}\|_{p',\Omega'_{4r}} \le c(|\Omega'_{4r}|)(1+\|\nabla\mathbf{u}\|_{p,\Omega'_{4r}}),
\end{equation*}
one has
\[
\begin{split}
&|I_4|\leqslant \|\bs\Phi\|_{p',\Omega'_{3r}}\|\Delta_{-\la}(\eta^p\Delta_{\lambda}\mathbf{u})\|_{p,\Omega'_{3r}}\\
&\leqslant  \|\bs\Phi\|_{p',\Omega'_{3r}}|\lambda|\bigg(\int_{\Omega'_{3r}}|\partial_k(\eta^p\Delta_{\lambda}\mathbf{u})|^pdx\bigg)^{1/p}\\
&\leqslant   \|\bs\Phi\|_{p',\Omega'_{3r}}|\lambda|\bigg(\bigg(\int_{\Omega'_{3r}}|p\eta^{p-1}\nabla\eta\Delta_{\lambda}\mathbf{u}|^p\bigg)^{1/p}
+\bigg(\int_{\Omega'_{3r}}|\eta^p\partial_k\Delta_{\lambda}\mathbf{u}|^pdx\bigg)^{1/p}\bigg)\\
&\leqslant  c\|\bs\Phi\|_{p',\Omega'_{3r}}|\lambda|\bigg(\frac{1}{r}|\lambda|\|\nabla\mathbf{u}\|_{p,\Omega'_{4r}}
+(\int_{\Omega'_{3r}}|\eta\nabla\Delta_{\lambda}\mathbf{u}|^pdx)^{1/p}\bigg)\\
&\leqslant \frac{c}{r}\lambda^2\|\bs\Phi\|_{p',\Omega'_{3r}}\|\nabla\mathbf{u}\|_{p,\Omega'_{4r}}
+c\|\Phi\|_{p',\Omega'_{3r}}|\lambda|\bigg(\frac{1}{r}|\lambda|\|\nabla\mathbf{u}\|_{p,\Omega'_{4r}}+
(\int_{\Omega'_{3r}}|\eta\Delta_{\lambda}\big{(}\mathcal{D}(\mathbf{u})\big{)}|^pdx)^{1/p}\bigg)\\
&\leqslant \frac{c}{r}\lambda^2\|\bs\Phi\|_{p',\Omega'_{3r}}\|\nabla\mathbf{u}\|_{p,\Omega'_{4r}}+c|\lambda|^{p'}\|\Phi\|^{p'}_{p',\Omega'_{3r}}
+\ve\int_{\Omega'_{3r}}|\eta\Delta_{\lambda}\big{(}\mathcal{D}(\mathbf{u})\big{)}|^pdx.\\
\end{split}
\]
Hence, we have
\begin{equation}
|I_4|\leqslant c(|\Omega'_{4r}|,r,\|\nabla\mathbf{u}\|_{p,\Omega'_{4r}})|\lambda|^{p'}+\ve\|\eta\Delta_{\lambda}\big{(}\mathcal{D}(\mathbf{u})\big{)}\|^p_{p,\Omega'_{3r}}.
\end{equation}
Finally, we start to deal with $I_5$. Note that
\begin{equation}
\begin{split}
&-I_5=\int_{\Omega'_{3r}}u_i\partial_iu_j\Delta_{-\lambda}(\eta^p\Delta_{\la}u_j)dx
=\int_{\Omega'_{3r}}\Delta_{\lambda}(u_i\partial_iu_j)\eta^p\Delta_{\la}u_jdx\\
&=\int_{\Omega'_{3r}}\Delta_{\lambda}u_i(\partial_iu_j)(x+\lambda e_k)\eta^p\Delta_{\la}u_jdx
+\int_{\Omega'_{3r}}u_i\Delta_{\lambda}(\partial_iu_j)\eta^p\Delta_{\la}u_jdx.\\
&\equiv I_{5,1}+I_{5,2}.
\end{split}
\end{equation}
It follows from a direct computation that
\begin{equation*}
|I_{5,1}|\leqslant c\|\Delta_{\lambda}\mathbf{u}\|^2_{2p',\Omega'_{3r}}\|\nabla\mathbf{u}\|_{p,\Omega'_{4r}}.
\end{equation*}
If $p<d$, we set $\theta=\frac{(d+2)p-3d}{2p}$. Then
\begin{equation*}
\|\Delta_{\lambda}\mathbf{u}\|^2_{2p',\Omega'_{2r}}\leqslant\|\Delta_{\lambda}\mathbf{u}\|^{2(1-\theta)}_{p*,\Omega'_{3r}}
\|\Delta_{\lambda}\mathbf{u}\|^{2\theta}_{p,\Omega'_{3r}}
\leqslant|\lambda|^{2\theta}\|\nabla\mathbf{u}\|^2_{p,\Omega'_{4r}},
\end{equation*}
where $p*=\frac{dp}{d-p}$.\\
If $d\leqslant p<3$, by $W^{1,d}\hookrightarrow L^q$ for any $q<\infty$, we then have that for any $\theta<1$,
\begin{equation*}
\|\Delta_{\lambda}\mathbf{u}\|^2_{2p',\Omega'_{3r}}\leqslant\|\Delta_{\lambda}
\mathbf{u}\|^{2(1-\theta)}_{q(\theta),\Omega'_{3r}}\|\Delta_{\lambda}\mathbf{u}\|^{2\theta}_{p,\Omega'_{3r}}
\leqslant c(|\Omega'_{4r}|)|\lambda|^{2\theta}\|\nabla\mathbf{u}\|^2_{p,\Omega'_{4r}}.
\end{equation*}
If $p\geqslant3$, due to $2p'\leqslant p$, one then has
\begin{equation*}
\|\Delta_{\lambda}\mathbf{u}\|^2_{2p',\Omega'_{3r}}\leqslant c(|\Omega'_{4r}|)\|\Delta_{\lambda}\mathbf{u}\|^2_{p,\Omega'_{3r}}
\leqslant c(|\Omega'_{4r}|)|\lambda|^2\|\nabla\mathbf{u}\|^2_{p,\Omega'_{4r}}.
\end{equation*}
Hence, we have
\begin{equation}
|I_{5,1}|\leqslant|\lambda|^{2\theta}\|\nabla\mathbf{u}\|^3_{p,\Omega'_{4r}}.
\end{equation}
In addition,
\begin{equation}
|I_{5,2}|=\f{p}{2}\Big{|}\int_{\Omega'_{3r}}u_i(\Delta_{\la}u_j)^2\eta^{p-1}\partial_i\eta dx\Big{|}
\leqslant\frac{c(|\Omega'_{4r}|)}{r}\|\Delta_{\la}\mathbf{u}\|^2_{2p',\Omega'_{3r}}\|\nabla\mathbf{u}\|_{p,\Omega'_{3r}}
\leqslant\frac{c(|\Omega'_{4r}|)}{r}|\lambda|^{2\theta}\|\nabla\mathbf{u}\|^3_{p,\Omega'_{4r}}.
\end{equation}
Combining (5.16) with (5.17) yields
\begin{equation}
|I_5|\leqslant c(|\Omega'_{4r}|)(1+\frac{1}{r})|\lambda|^{2\theta}\|\nabla\mathbf{u}\|^3_{\Omega'_{4r}}\leqslant c(|\Omega'_{4r}|,r,\|\nabla\mathbf{u}\|_{p,\Omega'_{4r}})|\lambda|^{2\theta}.
\end{equation}
Collecting all above estimates on $I_i$ ($1\le i\le 5$) and setting $\hat{\theta}=min\{\frac{p'}{2},\theta\}$, we
eventually obtain
\begin{equation*}
||\Delta_{\lambda}\mathcal{D}(u)\|_{p,\Omega'}\leqslant c(|\Omega'_{4r}|,r,\|\nabla\mathbf{u}\|_{p,\Omega_{4r}})|\lambda|^{\frac{2\hat{{\theta}}}{p}}.
\end{equation*}
Thus, by the characterization of the fractional order Sobolev space (see \cite{Adams} or \cite{Triebel}), we have
that for any  $\kappa\in [0,\frac{2\hat{\theta}}{p})$,
\begin{equation*}
\nabla \mathbf{u}\in W^{\kappa,p}(\Omega')
\end{equation*}
and
\begin{equation*}
\|\nabla\mathbf{u}\|_{\kappa,p,\Omega'}\leqslant c(\kappa,r,|\Omega'_{4r}|,\|\nabla\mathbf{u}\|_{p,\Omega'_{4r}}).
\end{equation*}

\vskip 0.2 true cm
$\mathbf{Case \ II}$.  $\mu_0>0$, $1<p<2$ \\

\vskip 0.1 true cm

As in Case I, we can find a function $\pi\in L^2(\Omega'_{4r})$ such that  for any
$\boldsymbol{\psi}\in D^{1,2}_0(\Omega'_{4r})$,
\begin{equation}
F(\boldsymbol{\psi})=(\pi,\nabla\cdot\boldsymbol{\psi})
\end{equation}
and
\begin{equation*}
\|\pi\|_{2,\Omega'_{4r}}\leqslant \|F\|_{D^{-1,2}_0({\Omega'_{4r}})}\leqslant c(|\Omega'_{4r}|)(1+\|\nabla\mathbf{u}\|_{2,{\Omega'_{4r}}}
+\|\nabla\mathbf{u}\|^{p-1}_{p,{\Omega'_{4r}}}+\|\nabla\mathbf{u}\|^2_{p,{\Omega'_{4r}}}+\|\nabla\mathbf{u}\|_{p,{\Omega'_{4r}}}).
\end{equation*}
Choosing $\Delta_{-\lambda}(\eta^2\Delta_{\lambda}\mathbf{u})$ as a test function in (5.19)
yields that
\begin{align}
&\mu_0\int_{\Omega'_{3r}}\mcD(\mathbf{u}):\mathcal{D}(\Delta_{-\lambda}(\eta^2\Delta_{\lambda}\mathbf{u}))dx\no\\
&\quad +\mu_1\int_{\Omega'_{3r}}|\mcD(\mathbf{u})+\mcD(\mathbf{a})|^{p-2}(\mcD(\mathbf{u})
+\mcD(\mathbf{a})):\mathcal{D}(\Delta_{-\lambda}(\eta^2\Delta_{\lambda}\mathbf{u}))dx\no\\
&=-\int_{\Omega'_{3r}}\mathbf{u}\cdot\nabla\mathbf{u}\cdot(\Delta_{-\la}(\eta^2\Delta_{\lambda}\mathbf{u}))dx
+\int_{\Omega'_{3r}}\pi\nabla\cdot(\Delta_{-\lambda}(\eta^p\Delta_{\la}\mathbf{u}))dx\no\\
&\quad -\int_{\Omega'_{3r}}\mathbf{u}\cdot\nabla\mathbf{a}\cdot(\Delta_{-\lambda}(\eta^2\Delta_{\lambda}\mathbf{u}))
-\int_{\Omega'_{3r}}\mathbf{a}\cdot\nabla\mathbf{u}\cdot(\Delta_{-\lambda}(\eta^2\Delta_{\lambda}\mathbf{u}))dx\no\\
&\quad -\int_{\Omega'_{3r}}\mathbf{a}\cdot\nabla\mathbf{a}\cdot(\Delta_{-\lambda}(\eta^2\Delta_{\lambda}\mathbf{u}))dx
+\mu_0\int_{\Omega'_{3r}}\Delta\mathbf{a}:\Delta_{-\lambda}(\eta^2\Delta_{\lambda}\mathbf{u})dx.
\end{align}
It follows from a direct computation that
\begin{multline*}
\mu_0\int_{\Omega'_{3r}}\mcD(\mathbf{u}):\mcD(\Delta_{-\lambda}(\eta^2\Delta_{\la}\mathbf{u}))dx\\
=2\mu_0\int_{\Omega'_{3r}}\Delta_{\lambda}(\mcD(\mathbf{u})):Sym(\Delta_{\lambda}\mathbf{u}\otimes\eta\nabla\eta)dx
+\mu_0\int_{\Omega'_{3r}}\eta^2(\Delta_{\la}\mcD(\mathbf{u}))^2dx.
\end{multline*}
Then we have
\begin{equation}
\begin{split}
&\mu_0\int_{\Omega'_{3r}}\eta^2(\Delta_{\la}\mcD(\mathbf{u}))^2dx\\
&\leqslant \frac{2\mu_0}{r}\bigg(\int_{\Omega'_{3r}}\eta^2(\Delta_{\lambda}\mathcal{D}(\mathbf{u}))^2dx\bigg)^{1/2}
\bigg(\int_{\Omega'_{3r}}|\Delta_{\lambda}\mathbf{u}|^2dx\bigg)^{1/2}
+\mu_0\int_{\Omega'_{3r}}\mathcal{D}(\mathbf{u})\cdot\mathcal{D}(\Delta_{-\lambda}(\eta^2\Delta_{\lambda}\mathbf{u}))dx\\
&\leqslant \frac{2\mu_0\varepsilon}{r}\int_{\Omega'_{3r}}\eta^2(\Delta_{\lambda}\mathcal{D}(\mathbf{u}))^2dx+\frac{\lambda^2}{4\ve}\int_{\Omega'_{4r}}|\nabla\mathbf{u}|^2dx
+\mu_0\int_{\Omega'_{3r}}\mathcal{D}(\mathbf{u})\cdot\mathcal{D}(\Delta_{-\lambda}(\eta^2\Delta_{\lambda}\mathbf{u}))dx.
\end{split}
\end{equation}
In addition
\begin{equation*}
\begin{split}
&\mu_1\int_{\Omega'_{3r}}|\mcD(\mathbf{u})+\mcD(\mathbf{a})|^{p-2}(\mcD(\mathbf{u})
+\mcD(\mathbf{a})):\mathcal{D}(\Delta_{-\lambda}(\eta^2\Delta_{\lambda}\mathbf{u}))dx\\
&=2\mu_1\int_{\Omega'_{3r}}|\mcD(\mathbf{u})+\mcD(\mathbf{a})|^{p-2}(\mcD(\mathbf{u})
+\mcD(\mathbf{a}))\Delta_{-\lambda}Sym(\Delta_{\lambda}\mathbf{u}\otimes\eta\nabla\eta)dx\\
&\quad +\mu_1\int_{\Omega'_{3r}}\eta^2\Delta_{\la}(|\mcD(\mathbf{u})
+\mcD(\mathbf{a})|^{p-2}(\mcD(\mathbf{u})+\mcD(\mathbf{a}))\Delta_{\lambda}(\mcD(\mathbf{u})+\mcD(\mathbf{a}))dx\\
&\quad -\mu_1\int_{\Omega'_{3r}}\eta^2|\mcD(\mathbf{u})
+\mcD(\mathbf{a})|^{p-2}(\mcD(\mathbf{u})+\mcD(\mathbf{a}))\Delta_{-\la}(\Delta_{\lambda}(\mcD(\mathbf{a})))dx
\end{split}
\end{equation*}
It follows from (5.8) that
\begin{equation}
\begin{split}
&\mu_1\int_{\Omega'_{3r}}|\mcD(\mathbf{u})+\mcD(\mathbf{a})|^{p-2}(\mcD(\mathbf{u})
+\mcD(\mathbf{a})):\mathcal{D}(\Delta_{-\lambda}(\eta^2\Delta_{\lambda}\mathbf{u}))dx\\
&\geqslant \f{\delta}{2} I(\mathbf{u})-\delta I(\mathbf{a}) -2\mu_1\|\mcD(\mathbf{u})+\mcD(\mathbf{a})\|^{p-1}_{p,\Omega'_{3r}}\|\Delta_{-\lambda}
Sym(\Delta_{\lambda}\mathbf{u}\otimes\eta\nabla\eta)\|_{p,\Omega'_{3r}}\\
&\quad -\mu_1\|\mcD(\mathbf{u})+\mcD(\mathbf{a})\|^{p-1}_{p,\Omega'_{3r}}\|\Delta_{-\la}(\Delta_{\lambda}(\mcD(\mathbf{a})))\|_{p,\Omega'_{3r}},
\end{split}
\end{equation}
where
\begin{equation*}
I(\mathbf{u})=\int_{\Omega'_{3r}}\eta^2(|(1+\mcD(\mathbf{u})+\mcD(\mathbf{a}))(x+\lambda e_k)|+|(\mcD(\mathbf{u})+\mcD(\mathbf{a}))(x)|)^{p-2}|\Delta_{\lambda}\mcD(\mathbf{u})|^2dx,
\end{equation*}
and
\begin{equation*}
I(\mathbf{a})=\int_{\Omega'_{3r}}\eta^2(|1+(\mcD(\mathbf{u})+\mcD(\mathbf{a}))(x+\lambda e_k)|+|(\mcD(\mathbf{u})+\mcD(\mathbf{a}))(x)|)^{p-2}|\Delta_{\lambda}\mcD(\mathbf{a})|^2dx.
\end{equation*}
From (5.20-(5.22), set $\ve=\f{r}{4}$, one has \begin{equation}
\begin{split}
&\frac{\mu_0}{2}\int_{\Omega'_{3r}}\eta^2(\Delta_{\lambda}\mathcal{D}(\mathbf{u}))^2dx+\f{\delta}{2} I(u)\\
&\leqslant \frac{\lambda^2}{r}\int_{\Omega'_{4r}}|\nabla\mathbf{u}|^2dx+\delta I(\mathbf{a})\\ &\quad +2\mu_1\|\mcD(\mathbf{u})+\mcD(\mathbf{a})\|^{p-1}_{p,\Omega'_{3r}}\|\Delta_{-\lambda}
Sym(\Delta_{\lambda}\mathbf{u}\otimes\eta\nabla\eta)\|_{p,\Omega'_{3r}}\\
&\quad +\mu_1\|\mcD(\mathbf{u})+\mcD(\mathbf{a})\|^{p-1}_{p,\Omega'_{3r}}\|\Delta_{-\la}(\Delta_{\lambda}(\mcD(\mathbf{a})))\|_{p,\Omega'_{3r}}\\
&\quad +\int_{\Omega'_{3r}}\pi\nabla\cdot(\Delta_{-\lambda}(\eta^p\Delta_{\lambda}\mathbf{u}))dx\\
&\quad +\int_{\Omega'_{3r}}\bs{\Phi}\cdot(\Delta_{-\lambda}(\eta^p\Delta_{\lambda}\mathbf{u}))dx\\
&\quad -\int_{\Omega'_{3r}}\mathbf{u}\cdot\nabla\mathbf{u}\cdot(\Delta_{-\lambda}(\eta^p\Delta_{\lambda}\mathbf{u}))dx\\
&\equiv \frac{\la^2}{r}\int_{\Omega'_{4r}}|\nabla\mathbf{u}|^2dx+\delta I(\mathbf{a})+I_1+...+I_5,
\end{split}
\end{equation}
where $\bs{\Phi}=\mu_0\Delta\mathbf{a}-\mathbf{u}\cdot\nabla\mathbf{a}
-\mathbf{a}\cdot\nabla\mathbf{u}-\mathbf{a}\cdot\nabla\mathbf{a}$.
Similarly to the treatment in Case I, we have
\begin{equation*}
\|\Delta_{-\lambda}Sym(\Delta_{\lambda}\mathbf{u}\otimes\eta\nabla\eta)\|_{p,\Omega'_{3r}}
\leqslant c\frac{\lambda^2}{r^2}\|\nabla\mathbf{u}\|_{p,\Omega'_{4r}}
+ c\frac{|\lambda|}{r}(\int_{\Omega'_{3r}}|\eta\mathcal{D}(\Delta_{\lambda}\mathbf{u})|^pdx)^{1/p}.
\end{equation*}
Since
\begin{equation}
\begin{split}
&\int_{\Omega'_{3r}}|\eta\mathcal{D}(\Delta_{\lambda}\mathbf{u})|^pdx\\
&\leqslant \int_{\Omega'_{3r}}\bigg(1+|(\mcD(\mathbf{u})+\mcD(\mathbf{a}))(x+\la e_k)|+|(\mcD(\mathbf{u})+\mcD(\mathbf{a}))(x)|\bigg)^{\frac{(p-2)p}{2}}|\eta\Delta_{\la}\mcD(\mathbf{u})|^p\\
&\qquad \times\bigg(1+|(\mcD(\mathbf{u})+\mcD(\mathbf{a}))(x+\lambda e_k)|+|(\mcD(\mathbf{u})+\mcD(\mathbf{a}))(x)|\bigg)^{\frac{(2-p)p}{2}}dx\\
&\leqslant 2I(\mathbf{u})^{\frac{p}{2}}\bigg(\int_{\Omega'_{4r}}(1+|\mcD(\mathbf{u})|+|\mcD(\mathbf{a})|)^pdx\bigg)^{\frac{2-p}{2}},
\end{split}
\end{equation}
we get
\begin{equation}
\begin{split}
&|I_1|\leqslant 2\mu_1\|\mcD(\mathbf{u})+\mcD(\mathbf{a})\|^{p-1}_{p,\Omega'_{3r}}
\bigg(\frac{\lambda^2}{r^2}\|\nabla\mathbf{u}\|_{p,\Omega'_{4r}}
+\frac{|\lambda|}{r}(\int_{\Omega'_{3r}}|\eta\mathcal{D}(\Delta_{\lambda}\mathbf{u})|^pdx)^{1/p}\bigg)\\
&\leqslant c(|\Omega'_{4r}|,r,\|\nabla\mathbf{u}\|_{p,\Omega_{4r}})\la^2
+c\frac{|\lambda|}{r}I(\mathbf{u})^{\frac{1}{2}}\bigg(\int_{\Omega'_{4r}}(1+|\mcD(\mathbf{u})|
+|\mathcal{D}(\mathbf{a})|)^pdx\bigg)^{\frac{1}{2}}\\
&\leqslant  c(|\Omega'_{4r}|,r,\|\nabla\mathbf{u}\|_{p,\Omega_{4r}})\la^2+\f{\delta}{2} I(\mathbf{u}).
\end{split}
\end{equation}
It is easy to get that
\begin{equation*}
|I(\mathbf{a})|\leqslant c(|\Omega'_{4r}|,r,\|\nabla\mathbf{u}\|_{p,\Omega_{4r}})\la^2
\end{equation*}
and
\begin{equation*}
|I_2|\leqslant c(|\Omega'_{4r}|,r,\|\nabla\mathbf{u}\|_{p,\Omega_{4r}})\la^2.
\end{equation*}
In addition,
\begin{align*}
&|I_3|=2\Big{|}\int_{\Omega'_{3r}}\pi\Delta_{-\lambda}(\eta\partial_i\eta\Delta_{\lambda}u_i)dx\Big{|}\\
&\leqslant 2\|\pi\|_{2,\Omega'_{3r}}|\lambda|\|\partial_k(\eta\partial_i\eta\Delta_{\lambda}u_i)\|_{2,\Omega'_{3r}}\\
&\leqslant \|\pi\|_{2,\Omega'_{3r}}|\lambda|\bigg(\frac{1}{r^2}\|\Delta_{\lambda}u_i\|_{2,\Omega'_{3r}}
+\frac{1}{r}\|\eta\Delta_{\lambda}\nabla\mathbf{u}\|_{2,\Omega'_{3r}}\bigg)\\
&\leqslant \|\pi\|_{2,\Omega'_{3r}}|\lambda|\bigg(\frac{1}{r^2}\|\Delta_{\lambda}u_i\|_{2,\Omega'_{3r}}
+\frac{|\la|}{r^2}\|\nabla\mathbf{u}\|_{2,\Omega'_{4r}}+\f{1}{r}\int_{\Omega'_{3r}}|\eta\Delta_{\lambda}\mcD(\mathbf{u})|^2dx\bigg)\\
&\leqslant \frac{c}{r^2}\|\pi\|_{2,\Omega'_{3r}}\|\nabla\mathbf{u}\|_{2,\Omega'_{4r}}\lambda^2+\frac{c}{\varepsilon r^2}\|\pi\|^2_{2,\Omega'_{3r}}\lambda^2
+\ve\int_{\Omega'_{3r}}|\eta\Delta_{\lambda}\mcD(\mathbf{u})|^2dx\\
&\leqslant c(|\Omega'_{4r}|,r,\|\nabla\mathbf{u}\|_{2,\Omega'_{4r}})\la^2+\ve\int_{\Omega'_{3r}}|\eta\Delta_{\lambda}\mcD(\mathbf{u})|^2dx.
\end{align*}
And
\[
\begin{split}
&|I_4|\leqslant \|\bs{\Phi}\|_{2,\Omega'_{3r}}\|\Delta_{-\la}(\eta^2\Delta_{\lambda}\mathbf{u})\|_{2,\Omega'_{3r}}\\
&\leqslant  \|\bs{\Phi}\|_{2,\Omega'_{3r}}|\lambda|\bigg(\int_{\Omega'_{3r}}|\partial_k(\eta^2\Delta_{\lambda}\mathbf{u})|^2dx\bigg)^{1/2}\\
&\leqslant   \|\bs{\Phi}\|_{2,\Omega'_{3r}}|\lambda|\bigg((\int_{\Omega'_{3r}}(\eta|\nabla\eta||\Delta_{\lambda}\mathbf{u}|)^2)^{1/2}
+(\int_{\Omega'_{3r}}|\eta^2\partial_k\Delta_{\lambda}\mathbf{u}|^2dx)^{1/2}\bigg)\\
&\leqslant  \|\bs{\Phi}\|_{2,\Omega'_{3r}}|\lambda|\bigg(\frac{1}{r}|\lambda|\|\nabla\mathbf{u}\|_{2,\Omega'_{3r}}
+(\int_{\Omega'_{3r}}|\eta\nabla\Delta_{\lambda}\mathbf{u}|^2dx)^{1/2}\bigg)\\
&\leqslant \frac{1}{r}\lambda^2\|\bs{\Phi}\|_{2,\Omega'_{3r}}\|\nabla\mathbf{u}\|_{2,\Omega'_{3r}}
+\|\bs{\Phi}\|_{2,\Omega'_{3r}}|\lambda|\bigg(\frac{1}{r}|\lambda|\|\nabla\mathbf{u}\|_{2,\Omega'_{3r}}+
(\int_{\Omega'_{3r}}|\eta\Delta_{\lambda}\mathcal{D}(\mathbf{u})|^2dx)^{1/2}\bigg)\\
&\leqslant c(|\Omega'_{4r}|,r,\|\nabla\mathbf{u}\|_{2,\Omega'_{4r}})\la^2+\ve\int_{\Omega'_{3r}}|\eta\Delta_{\lambda}\mathcal{D}(\mathbf{u})|^2dx.
\end{split}
\]
Finally, we start to deal with the term $I_5$. It follows from a direct computation that
\begin{equation}
\begin{split}
&-I_5=\ds\sum_{i,j=1}^d\int_{\Omega'_{3r}}u_i\partial_iu_j\Delta_{-\lambda}(\eta^2\Delta_{\lambda}u_j)dx
=\ds\sum_{i,j=1}^d\int_{\Omega'_{3r}}\Delta_{\lambda}(u_i\partial_iu_j)\eta^2\Delta_{\lambda}u_jdx\\
&=\ds\sum_{i,j=1}^d\int_{\Omega'_{3r}}\Delta_{\lambda}u_i(\partial_iu_j)(x+\lambda e_k)\eta^2\Delta_{\lambda}u_jdx
+\ds\sum_{i,j=1}^d\int_{\Omega'_{3r}}u_i\Delta_{\lambda}(\partial_iu_j)\eta^2\Delta_{\lambda}u_jdx\\
&\equiv I_{5,1}+I_{5,2}.
\end{split}
\end{equation}
Obviously,
\begin{equation}
|I_{5,1}|\leqslant\|\Delta_{\lambda}\mathbf{u}\|^2_{4,\Omega'_{3r}}\|\nabla\mathbf{u}\|_{2,\Omega'_{4r}}.
\end{equation}
If $d=2$, by $W^{1,2}\hookrightarrow L^{q}$ for any $q<\infty$, we then have that for any $\theta<1$,
\begin{equation*}
\|\Delta_{\la}\mathbf{u}\|^2_{4,\Omega'_{3r}}\leqslant\|\Delta_{\lambda}\mathbf{u}\|^{2(1-\theta)}_{q(\theta),\Omega'_{3r}}
\|\Delta_{\lambda}\mathbf{u}\|^{2\theta}_{2,\Omega'_{3r}}
\leqslant c(|\Omega'_{4r}|)|\la|^{2\theta}\|\nabla\mathbf{u}\|^2_{2,\Omega'_{4r}}.
\end{equation*}
If $d=3$, we set $\theta=\f{1}{4}$. Then
\begin{equation*}
\|\Delta_{\lambda}\mathbf{u}\|^2_{4,\Omega'_{3r}}\leqslant\|\Delta_{\lambda}\mathbf{u}\|^{2(1-\theta)}_{6,\Omega'_{3r}}
\|\Delta_{\lambda}\mathbf{u}\|^{2\theta}_{2,\Omega'_{3r}}
\leqslant c(|\Omega'_{4r}|)|\lambda|^{2\theta}\|\nabla\mathbf{u}\|^2_{2,\Omega'_{4r}}
\end{equation*}
and
\begin{equation}
|I_{5,1}|\leqslant c(|\Omega'_{4r}|)|\la|^{2\theta}\|\nabla\mathbf{u}\|^3_{2,\Omega'_{4r}}.
\end{equation}
On the other hand,
\begin{equation*}
\begin{split}
&|I_{5,2}|=\f{p}{2}|\int_{\Omega'_{3r}}u_i(\Delta_{\la}u_j)^2\eta^{p-1}\partial_i\eta dx|
\leqslant\frac{c}{r}\|\Delta_{\la}\mathbf{u}\|^2_{4,\Omega'_{3r}}\|\mathbf{u}\|_{2,\Omega'_{3r}}\\
&\leqslant c(|\Omega'_{4r}|,r)\|\Delta_{\la}\mathbf{u}\|^2_{4,\Omega'_{3r}}\|\nabla\mathbf{u}\|_{2,\Omega'_{4r}}.
\end{split}
\end{equation*}
Hence, we can obtain that
\begin{equation}
|I_5|\leqslant  c(\ve,|\Omega'_{4r}|,r,\|\nabla\mathbf{u}\|_{p,\Omega'_{4r}}) |\lambda|^{2\theta}.
\end{equation}
Collecting all above estimates on $I_i$ ($1\le i\le 4$), we eventually get
\begin{equation*}
||\Delta_{\la}\mathcal{D}(\mathbf{u})\|_{2,\Omega'}\leqslant c(|\Omega'_{4r}|,r,\|\nabla\mathbf{u}\|_{2,\Omega'_{4r}})|\lambda|^{\theta}
\end{equation*}
and
\begin{equation*}
||\Delta_{\la}\mathcal{D}(\mathbf{u})\|_{p,\Omega'}\leqslant c(|\Omega'_{4r}|,r,\|\nabla\mathbf{u}\|_{2,\Omega'_{4r}})|\lambda|^{\f{2\theta}{p}}.
\end{equation*}
Thus, we have that for any $\varepsilon>0$,
\begin{equation*}
\nabla \mathbf{u}\in W^{\f{2\theta}{p}-\varepsilon,p}(\Omega')\cap W^{\theta-\ve,2}(\Omega')
\end{equation*}
and
\begin{equation*}
\|\nabla\mathbf{u}\|_{\theta-\ve,p,\Omega'}+\|\nabla\mathbf{u}\|_{\f{2\theta}{p}-\ve,p,\Omega'}
\leqslant c(\ve,|\Omega'_{4r}|,r,\|\nabla\mathbf{u}\|_{p,\Omega'_{4r}}).
\end{equation*}
\qquad\qquad \qquad \qquad \qquad \qquad \qquad \qquad \qquad \qquad \qquad
\qquad \qquad \qquad \qquad \qquad \qquad \qquad \qquad \qquad  $\square$

\section{Solvability of Ladyzhenskaya-Solonnikov Problem I (2.4)}

To solve  Ladyzhenskaya-Solonnikov Problem I (2.4) under some suitable conditions, based on Sections 4-5,
we take the following three parts.

\subsection{Part 1. Uniform estimate of $\|\mathbf{u}^T\|_{p,\Omega(t)}$}

In what follows, for convenience and without loss of generality, we
assume that
\begin{equation*}
\Omega=\{x: x_1\in\mathbb{R}, x'\in\Sigma(x_1)\}
\end{equation*}
and
\begin{equation*}
\Omega(t)=\{x\in\Omega: -t<x_1<t\}.
\end{equation*}
Taking the inner product of \eqref{eq:truncated doamin}$_1$
with $\mathbf{u}^T$ and integrating by parts over $\Omega(t)$ yield
\begin{equation}
\begin{split}
&\mu_0\|\mcD(\mathbf{v}^T)\|^2_{2,\Omega(t)}+\mu_1\|\mcD(\mathbf{v}^T)\|^p_{p,\Omega(t)}\\
&=I_1+I_2+I_3+I_4+I_5+I_6
\end{split}
\end{equation}
with
\begin{align*}
&I_1=\mu_0\int_{\Sigma(t)}\mathbf{u}^T\cdot\mcD(\mathbf{v}^T)\cdot\mathbf{e}_1dS
+\mu_1\int_{\Sigma(t)}\mathbf{u}^T\cdot(|\mcD(\mathbf{v}^T)|^{p-2}\mcD(\mathbf{v}^T))\cdot \mathbf{e}_1dS,\\
&I_2=-\mu_0\int_{\Sigma(-t)}\mathbf{u}^T\cdot\mcD(\mathbf{v}^T)\cdot \mathbf{e}_1dS
-\mu_1\int_{\Sigma(-t)}\mathbf{u}^T\cdot(|\mcD(\mathbf{v}^T)|^{p-2}\mcD(\mathbf{v}^T))\cdot \mathbf{e}_1dS,\\
&I_3=-\int_{\Omega(t)}\mathbf{v}^T\cdot\nabla\mathbf{v}^T\cdot\mathbf{u}^Tdx,\\
&I_4=-\int_{\Sigma(t)}\pi^T u_1^T dS+\int_{\Sigma(-t)}\pi^T u_1^T dS,\\
&I_5=\mu_0\int_{\Omega(t)}\mcD(\mathbf{v}^T):\mcD(\mathbf{a})dx,\\
&I_6=\mu_1\int_{\Omega(t)}|\mcD(\mathbf{v}^T)|^{p-2}\mcD(\mathbf{v}^T):\mcD(\mathbf{a})dx.\\
\end{align*}
Next we deal with the terms $I_i$ for $1\le i\le 6$ in (6.1). We still divide the related process into the following two cases:
\vskip 0.2 true cm

$\mathbf{Case\ I}$. {\bf $\mu_0>0$, $p>1$}

\vskip 0.1 true cm
Using the H\"{o}lder and Poincar\'{e} inequality we get that
\begin{equation}
\begin{split}
\Big{|}\int_{\eta-1}^{\eta}I_1dt\Big{|}=&\Big{|}\mu_0\int_{\omega^{+}_{\eta}}\mathbf{u}^T\cdot\mcD(\mathbf{v}^T)\cdot\mathbf{e}_1dx
+\mu_1\int_{\omega^{+}_{\eta}}\mathbf{u}^T\cdot(|\mcD({\mathbf{v}^T})|^{p-2}\mcD(\mathbf{v}^T))\cdot \mathbf{e}_1dx\Big{|}\\
\leqslant &\mu_0\|\mathbf{u}^T\|_{2,\omega^{+}_{\eta}}\|\mcD(\mathbf{v}^T)\|_{2,\omega^{+}_{\eta}}+\mu_1 \|\mathbf{u}^T\|_{p,\omega^{+}_{\eta}}\|\mcD(\mathbf{v}^T)\|^{p-1}_{p,\omega^{+}_{\eta}}\\
\leqslant &\mu_0\|\mathbf{u}^T\|^2_{2,\omega^{+}_{\eta}}+\mu_0\|\mcD(\mathbf{v}^T)\|^2_{2,\omega^{+}_{\eta}}
+\mu_1\|\mathbf{u}^T\|^p_{p,\omega^{+}_{\eta}}+\mu_1\|\mcD(\mathbf{v}^T)\|^{p}_{p,\omega^{+}_{\eta}}\\
\leqslant& c\|\nabla\mathbf{u}^T\|^2_{2,\omega^{+}_{\eta}}+c\|\nabla\mathbf{u}^T\|^p_{p,\omega^{+}_{\eta}}+c\al,
\end{split}
\end{equation}
where $\omega^{+}_{\eta}=\{x\in\Omega: \eta-1<x_1<\eta\}$, and $\al=\max\limits_{1\le i\le N}|\al_i|$.
Similarly, we have
\begin{equation}
\Big{|}\int_{\eta-1}^{\eta}I_2 \ dt\Big{|}\leqslant c\|\nabla\mathbf{u}^T\|^2_{2,\omega^{-}_{\eta}}+c\|\nabla\mathbf{u}^T\|^p_{p,\omega^{-}_{\eta}}+c\al,
\end{equation}
where $\omega^{-}_{\eta}=\{x\in\Omega: -\eta<x_1<-\eta+1\}$.

Next, we estimate $I_3$. Note that
\begin{equation}
\begin{split}
\int_{\Omega(t)}\mathbf{v}^T\cdot\nabla\mathbf{v}^T\cdot\mathbf{u}^Tdx
=&\int_{\Sigma(t)}(\mathbf{a}\cdot\mathbf{u}^T)(\mathbf{a}\cdot\mathbf{e}_1)dS
-\int_{\Sigma(-t)}(\mathbf{a}\cdot\mathbf{u}^T)(\mathbf{a}\cdot\mathbf{e}_1)dS\\
&+\int_{\Sigma(t)}\frac{(\mathbf{u}^T)^2}{2}\mathbf{v}^T\cdot\mathbf{e}_1dS
-\int_{\Sigma(-t)}\frac{(\mathbf{u}^T)^2}{2}\mathbf{v}^T\cdot\mathbf{e}_1dS\\
&-\int_{\Omega(t)}\mathbf{u}^T\cdot\nabla\mathbf{u}^T\cdot\mathbf{a}\ dx
+\int_{\Omega(t)}\mathbf{a}\cdot\nabla\mathbf{a}\cdot\mathbf{u}^Tdx.
\end{split}
\end{equation}
By H\"{o}lder inequality and Lemma 3.1 we arrive at
\begin{equation}
\begin{split}
&\Big{|}\int_{\eta-1}^{\eta}\int_{\Sigma(t)}(\mathbf{a}\cdot\mathbf{u}^T)(\mathbf{a}\cdot\mathbf{e}_1)dydt
+\int_{\eta-1}^{\eta}\int_{\Sigma(t)}\frac{(\mathbf{u}^T)^2}{2}\mathbf{v}^T\cdot\mathbf{e}_1dydt\Big{|}\no\\
&\leqslant \|\mathbf{a}\|^2_{4,\omega^+_{\eta}}\|\mathbf{u}^T\|_{2,\omega^+_{\eta}}+\|\mathbf{u}^T\|^2_{4,\omega^+_{\eta}}\|\mathbf{v}^T\|_{2,\omega^+_{\eta}}\\
&\leqslant \|\mathbf{a}\|^2_{4,\omega^+_{\eta}}\|\mathbf{u}^T\|_{2,\omega^+_{\eta}}+\|\mathbf{u}^T\|^2_{4,\omega^+_{\eta}}\|\mathbf{u}^T\|_{2,\omega^+_{\eta}}
+\|\mathbf{a}\|_{2,\omega^+_{\eta}}\|\mathbf{u}^T\|^2_{4,\omega^+_{\eta}}\\
&\leqslant c(\|\mcD(\mathbf{u}^T)\|^3_{2,\omega^+_{\eta}}+\|\mcD(\mathbf{u}^T)\|^2_{2,\omega^+_{\eta}}
+\|\mcD(\mathbf{u}^T)\|_{2,\omega^+_{\eta}}).
\end{split}
\end{equation}
Similarly,
\begin{align}
\Big{|}\int_{\eta-1}^{\eta}\int_{\Sigma(-t)}&(\mathbf{a}\cdot\mathbf{u}^T)(\mathbf{a}\cdot\mathbf{e}_1)dydt
+\int_{\eta-1}^{\eta}\int_{\Sigma(-t)}\frac{\mathbf{u}^2}{2}\mathbf{v}^T\cdot\mathbf{e}_1dydt\Big{|}\no\\
&\leqslant c(\|\mcD(\mathbf{u}^T)\|^3_{2,\omega^-_{\eta}}+\|\mcD(\mathbf{u}^T)\|^2_{2,\omega^-_{\eta}}
+\|\mcD(\mathbf{u}^T)\|_{2,\omega^-_{\eta}}).
\end{align}
In addition, by Lemma 3.5 (iii) and (iv) we can arrive at
\begin{equation*}
\Big{|}\int_{\Omega(t)}\mathbf{u}^T\cdot\nabla\mathbf{u}^T\cdot\mathbf{a}\ dx\Big{|}\leqslant c(\ve)\int_{\Omega(t)}|\mathbf{a}|^2|\mathbf{u}^T|^2dx+\f{\ve}{2}\int_{\Omega(t)}|\nabla\mathbf{u}^T|^2dx
\leqslant \ve\int_{\Omega(t)}|\nabla\mathbf{u}^T|^2dx.
\end{equation*}
Note that
\begin{equation}
\Big{|}\int_{\Omega(t)}\mathbf{a}\cdot\nabla\mathbf{a}\cdot\mathbf{u}^T\ dx\Big{|}\leqslant\ve\|\mathbf{u}\|^2_{2,\Omega(t)}+c\|\mathbf{a}\cdot\nabla\mathbf{a}\|^2_{2,\Omega(t)}
\leqslant c\ve\|\nabla\mathbf{u}^T\|^2_{2,\Omega(t)}+c\al t.
\end{equation}
Hence,
\begin{equation}
\begin{split}
\Big{|}\int_{\eta-1}^{\eta}I_3dt\Big{|}&\leqslant\ve\int_{\eta-1}^{\eta}(\int_{\Omega(t)}|\nabla\mathbf{u}^T|^2dx)dt
+c(\|\nabla\mathbf{u}^T\|^3_{p,\omega^{+}_{\eta}}
+\|\nabla\mathbf{u}^T\|^2_{p,\omega^{+}_{\eta}}+\|\nabla\mathbf{u}^T\|_{p,\omega^{+}_{\eta}})\\
&\quad +c(\|\nabla\mathbf{u}^T\|^3_{p,\omega^{-}_{\eta}}
+\|\nabla\mathbf{u}^T\|^2_{p,\omega^{-}_{\eta}}+\|\nabla\mathbf{u}^T\|_{p,\omega^{-}_{\eta}})+c\al\eta.
\end{split}
\end{equation}
Next, we treat the term $I_4$. By $\int_{\Sigma(\pm t)}u^T_1dS=0$ for any $t>0$, then from Remark 3.2, we can find a vector
$\mathbf{w}\in W_0^{1,p}(\omega^{\pm}_{\eta})\cap W_0^{1,2}(\omega^{\pm}_{\eta})$ such that $div\mathbf{w}=u^T_1$ in $\omega^{\pm}_{\eta}$.
A direct computation yields
\begin{align*}
&\Big{|}\int_{\omega^{\pm}_\eta}\pi^T u_1^T\ dx|=|\int_{\omega^{\pm}_\eta}\pi^T \ div\mathbf{w} \ dx\Big{|}\\
&=\Big{|}\int_{\omega^{\pm}_{\eta}}(\mu_0\mcD(\mathbf{v}^T)+\mu_1|\mcD({\mathbf{v}^T})|^{p-2}\mcD(\mathbf{v}^T)):\mcD(\mathbf{w})dx
-\int_{\omega^{\pm}_{\eta}}\mathbf{v}^T\cdot\nabla\mathbf{w}\cdot\mathbf{v}^Tdx\Big{|}\\
&\leqslant\mu_0\|\mcD(\mathbf{v}^T)\|_{2,\omega^{\pm}_{\eta}}\|\mcD(\mathbf{w})\|_{2,\omega^{\pm}_{\eta}}
+\mu_1\|\mcD(\mathbf{v}^T)\|^{p-1}_{p,\omega^{\pm}_{\eta}}\|\mcD(\mathbf{w})\|_{p,\omega^{\pm}_{\eta}}+\|\mathbf{v}^T\|^2_{4,\omega^{\pm}_{\eta}}\|\nabla\mathbf{w}\|_{2,\omega^{\pm}_{\eta}}\\
&\leqslant c\|\mcD(\mathbf{v}^T)\|_{2,\omega^{\pm}_{\eta}}\|\nabla\mathbf{u}^T\|_{2,\omega^{\pm}_{\eta}}
+c\|\mcD(\mathbf{v}^T)\|^{p-1}_{p,\omega^{\pm}_{\eta}}\|\nabla\mathbf{u}^T\|_{p,\omega^{\pm}_{\eta}}+c\|\mathbf{v}^T\|^2_{4,\omega^{\pm}_{\eta}}\|\nabla\mathbf{u}^T\|_{2,\omega^{\pm}_{\eta}}\\
&\leqslant c\|\nabla\mathbf{u}^T\|^2_{2,\omega^{\pm}_{\eta}}+c\|\nabla\mathbf{u}^T\|_{2,\omega^{\pm}_{\eta}}
+c\|\nabla\mathbf{u}^T\|^p_{p,\omega^{\pm}_{\eta}}+c\|\nabla\mathbf{u}^T\|_{p,\omega^{\pm}_{\eta}}+c\|\nabla\mathbf{u}^T\|^3_{2,\omega^{\pm}_{\eta}}.
\end{align*}
This means that
\begin{equation}
\begin{split}
\Big{|}\int_{\eta-1}^{\eta}I_4 \ dt\Big{|}&\leqslant c(\|\nabla\mathbf{u}^T\|_{2,\omega^{+}_\eta}+\|\nabla\mathbf{u}^T\|^2_{2,\omega^{+}_\eta}
+\|\nabla\mathbf{u}^T\|_{2,\omega^{+}_\eta}^3+\|\nabla\mathbf{u}^T\|_{p,\omega^{+}_\eta}+\|\nabla\mathbf{u}^T\|^p_{p,\omega^{+}_\eta})\\
&\quad +c(\|\nabla\mathbf{u}^T\|_{2,\omega^{-}_\eta}+\|\nabla\mathbf{u}^T\|^2_{2,\omega^{-}_\eta}
+\|\nabla\mathbf{u}^T\|_{2,\omega^{-}_\eta}^3+\|\nabla\mathbf{u}^T\|_{p,\omega^{-}_\eta}+\|\nabla\mathbf{u}^T\|^p_{p,\omega^{-}_\eta}).
\end{split}
\end{equation}
Since
\[
\begin{split}
&\Big{|}\mu_0\int_{\Omega(t)}\mcD(\mathbf{v}^T):\mcD(\mathbf{a})dx+\mu_1\int_{\Omega(t)}
|\mcD(\mathbf{v}^T)|^{p-2}\mcD(\mathbf{v}^T):\mcD(\mathbf{a})dx\Big{|}\\
&\leqslant\ve\|\mcD(\mathbf{v}^T)\|^2_{2,\Omega(t)}+\ve\|\mcD(\mathbf{v}^T)\|^p_{p,\Omega(t)}
+c\|\mcD(\mathbf{a})\|^2_{2,\Omega(t)}+c\|\mcD(\mathbf{a})\|^p_{p,\Omega(t)}\\
&\leqslant\ve\|\mcD(\mathbf{u}^T)\|^2_{2,\Omega(t)}+\ve\|\mcD(\mathbf{u}^T)\|^p_{p,\Omega(t)}+c\al t,
\end{split}
\]
we arrive at
\begin{equation}
\int_{\eta-1}^{\eta}(|I_5|+|I_6|)dt\leqslant\ve\int_{\eta-1}^{\eta}\|\mcD(\mathbf{u}^T)\|^2_{2,\Omega(t)}dt+\ve\int_{\eta-1}^{\eta}\|\mcD(\mathbf{u}^T)\|^p_{p,\Omega(t)}dt+c\al\eta,
\end{equation}
Collecting all above estimates yields
\begin{equation}
\begin{split}
&\int_{\eta-1}^{\eta}(\|\nabla\mathbf{u}^T\|^2_{2,\Omega(t)}+\|\nabla\mathbf{u}^T\|^p_{p,\Omega(t)})dt\\
&\leqslant c(\|\nabla\mathbf{u}^T\|^3_{2,\omega^{+}_{\eta}}+\|\nabla\mathbf{u}^T\|^p_{p,\omega^{+}_{\eta}})
+c(\|\nabla\mathbf{u}^T\|^3_{2,\omega^{-}_{\eta}}+\|\nabla\mathbf{u}^T\|^p_{p,\omega^{-}_{\eta}})+c\al\eta+c\al.
\end{split}
\end{equation}
Let $y(t)=\int_{\Omega(t)}(|\nabla\mathbf{u}^T|^2+|\nabla\mathbf{u}^T|^p)dx$
and $z(\eta)=\int^{\eta}_{\eta-1}y(t)dt$. Then from (6.10)
\begin{equation}
z(\eta)\leqslant c_3(z'(\eta)+z'(\eta)^{\frac{3}{2}})+c_4\al\eta+c_5\al.
\end{equation}
To apply Lemma 3.3 (i), we set $\Psi(\tau)=c_3(\tau+\tau^{\frac{3}{2}})$,
$\dl=\f12$, $t_0=1$,
and $\varphi(\eta)=2c_4 \al\eta+2c_5\al$, where $c_5>0$ satisfies
\begin{equation}
c_4\al+c_5\al \geqslant \Psi(2c_4\al).
\end{equation}
In addition, it follows from the proof of Theorem 4.1 that
\begin{equation}
z(T)\leqslant \varphi(T).
\end{equation}
Therefore, according (6.11)-(6.13) and Lemma 3.3 (i), we arrive at
\begin{equation}
y(\eta-1)\leqslant z(\eta)\leqslant \vp(\eta), \end{equation}
which means  that for any $t\in[1,T]$,
\begin{equation}
\label{uniform estimate}
y(t)\leqslant 2c_4 \al(t+1)+2c_5\al.
\end{equation}
\vskip 0.2 true cm

$\mathbf{Case\ II}$. {\bf $\mu_0=0$, $p>2$}

\vskip 0.2 true cm
As in Case I, by the same calculation we can get
\begin{equation}
\Big{|}\int_{\eta-1}^{\eta}I_1dx_1\Big{|}+\Big{|}\int_{\eta-1}^{\eta}I_2dx_1\Big{|}\leqslant c\|\nabla\mathbf{u}^T\|^p_{p,\omega^{+}_{\eta}}+c\al.
\end{equation}
Next, we estimate $I_3$ as in Case I. Note that
\begin{equation}
\begin{split}
&\Big{|}\int_{\eta-1}^{\eta}\int_{\Sigma(t)}(\mathbf{a}\cdot\mathbf{u}^T)(\mathbf{a}\cdot\mathbf{e}_1)dydt
+\int_{\eta-1}^{\eta}\int_{\Sigma(t)}\frac{(\mathbf{u}^T)^2}{2}\mathbf{v}^T\cdot\mathbf{e}_1dydt\Big{|}\no\\
&\leqslant \|\mathbf{a}\|^2_{2p',\omega^+_{\eta}}\|\mathbf{u}^T\|_{p,\omega^+_{\eta}}+\|\mathbf{u}^T\|^2_{2p',\omega^+_{\eta}}\|\mathbf{v}^T\|_{p,\omega^+_{\eta}}\\
&\leqslant c(\|\nabla\mathbf{u}^T\|^3_{p,\omega^+_{\eta}}+\|\nabla\mathbf{u}^T\|^2_{p,\omega^+_{\eta}}
+\|\nabla\mathbf{u}^T\|_{p,\omega^+_{\eta}}).
\end{split}
\end{equation}
Similarly,
\begin{equation}
\begin{split}
&\Big{|}\int_{\eta-1}^{\eta}\int_{\Sigma(-t)}(\mathbf{a}\cdot\mathbf{u}^T)(\mathbf{a}\cdot\mathbf{e}_1)dydt
+\int_{\eta-1}^{\eta}\int_{\Sigma(-t)}\frac{(\mathbf{u}^T)^2}{2}\mathbf{v}^T\cdot\mathbf{e}_1dydt\Big{|}\no\\
&\leqslant c(\|\nabla\mathbf{u}^T\|^3_{p,\omega^-_{\eta}}+\|\nabla\mathbf{u}^T\|^2_{p,\omega^-_{\eta}}
+\|\nabla\mathbf{u}^T\|_{p,\omega^-_{\eta}}).
\end{split}
\end{equation}
In addition,
\begin{align*}
\Big{|}\int_{\Omega(t)}\mathbf{u}^T\cdot\nabla\mathbf{u}^T\cdot\mathbf{a}dx\Big{|}
&\leqslant\ve\int_{\Omega(t)}|\mathbf{u}^T|^pdx
+\ve\int_{\Omega(t)}|\nabla\mathbf{u}^T|^pdx
+c\int_{\Omega(t)}|\mathbf{a}|^{\frac{p}{p-2}}dx\\
&\leqslant \ve\int_{\Omega(t)}|\nabla\mathbf{u}^T|^pdx+c\al t
\end{align*}
and
\begin{align*}
\Big{|}\int_{\Omega(t)}\mathbf{a}\cdot\nabla\mathbf{a}\cdot\mathbf{u}^Tdx\Big{|}&\leqslant c\int_{\Omega(t)}|\mathbf{a}\cdot\nabla\mathbf{a}|^{p'}dx
+\ve\int_{\Omega(t)}|\mathbf{u}^T|^pdx\\
&\leqslant \ve\int_{\Omega(t)}|\nabla\mathbf{u}^T|^{p}dx+c\al t.
\end{align*}
Hence it follows from the expression of $I_3$ and the estimates above that
\begin{equation}
\begin{split}
\Big{|}\int_{\eta-1}^{\eta}I_3dt\Big{|}&\leqslant\ve\int_{\eta-1}^{\eta}(\int_{\Omega(t)}|\nabla\mathbf{u}^T|^{p}dx)dt
+c(\|\nabla\mathbf{u}^T\|^3_{p,\omega^{+}_{\eta}}+\|\nabla\mathbf{u}^T\|^2_{p,\omega^{+}_{\eta}}
+\|\nabla\mathbf{u}^T\|_{p,\omega^{+}_{\eta}})\\
&+c(\|\nabla\mathbf{u}^T\|^3_{p,\omega^{-}_{\eta}}+\|\nabla\mathbf{u}^T\|^2_{p,\omega^{-}_{\eta}}
+\|\nabla\mathbf{u}^T\|_{p,\omega^{-}_{\eta}})+c\al\eta.
\end{split}
\end{equation}
In addition, we can find a function $\mathbf{w}\in W_0^{1,p}(\omega^{\pm}_{\eta})$ such that $div\mathbf{w}=u^T_1$ in $\omega^{\pm}_{\eta}$.
Then we have that
\begin{align*}
&\Big{|}\int_{\omega^{\pm}_\eta}\pi^T u_1^Tdx\Big{|}=
\Big{|}\mu_1\int_{\omega^{\pm}_{\eta}}|\mcD({\mathbf{v}^T})|^{p-2}\mcD(\mathbf{v}^T):\mcD(\mathbf{w})dx
-\int_{\omega^{\pm}_{\eta}}\mathbf{v}^T\cdot\nabla\mathbf{w}\cdot\mathbf{v}^Tdx\Big{|}\\
&\leqslant \|\mcD(\mathbf{v}^T)\|^{p-1}_{p,\omega^{\pm}_{\eta}}\|\mcD(\mathbf{w})\|_{p,\omega^{\pm}_{\eta}}+
\|\mathbf{v}^T\|^2_{4,\omega^{\pm}_{\eta}}\|\nabla\mathbf{w}\|_{2,\omega^{\pm}_{\eta}}\\
&\leqslant \|\mcD(\mathbf{v}^T)\|^p_{p,\omega^{\pm}_{\eta}}+\|\mcD(\mathbf{w})\|^p_{p,\omega^{\pm}_{\eta}}+
\|\nabla\mathbf{u}^T\|^3_{2,\omega^{\pm}_{\eta}}+\|\mathbf{a}\|^2_{4,\omega^{\pm}_{\eta}}\|\nabla\mathbf{u}^T\|_{2,\omega^{\pm}_{\eta}}\\
&\leqslant c(\|\nabla\mathbf{u}^T\|^p_{p,\omega^{\pm}_\eta}+\|\nabla\mathbf{u}^T\|_{p,\omega^{\pm}_\eta}^3)+c\al.
\end{align*}
This yields
\begin{equation}
\Big{|}\int_{\eta-1}^{\eta}I_4 \ dt\Big{|}\leqslant c(\|\nabla\mathbf{u}^T\|^p_{p,\omega^{+}_\eta}
+\|\nabla\mathbf{u}^T\|_{p,\omega^{+}_\eta}^3)+ c(\|\nabla\mathbf{u}^T\|^p_{p,\omega^{-}_\eta}
+\|\nabla\mathbf{u}^T\|_{p,\omega^{-}_\eta}^3)+c\al.
\end{equation}
Finally, due to
\[
\begin{split}
&\Big{|}\mu_1\int_{\Omega(t)}
|\mcD(\mathbf{v}^T)|^{p-2}\mcD(\mathbf{v}^T):\mcD(\mathbf{a})dx\Big{|}\\
&\leqslant\ve\|\mcD(\mathbf{v}^T)\|^p_{p,\Omega(t)}+c\|\mcD(\mathbf{a})\|^p_{p,\Omega(t)}\\
&\leqslant\ve\|\mcD(\mathbf{u}^T)\|^p_{p,\Omega(t)}+c\al t,
\end{split}
\]
we arrive at
\begin{equation}
\int_{\eta-1}^{\eta}(|I_5|+|I_6|)dt\leqslant\ve\int_{\eta-1}^{\eta}\|\mcD(\mathbf{u}^T)\|^p_{p,\Omega(t)}dt+c\al\eta,
\end{equation}
Collecting all above estimates yields
\begin{equation}
\int_{\eta-1}^{\eta}\|\nabla\mathbf{u}^T\|^p_{p,\Omega(t)}dt\leqslant c(\|\nabla\mathbf{u}^T\|^3_{2,\omega^{+}_{\eta}}+\|\nabla\mathbf{u}^T\|^p_{p,\omega^{+}_{\eta}})
+c(\|\nabla\mathbf{u}^T\|^3_{2,\omega^{-}_{\eta}}+\|\nabla\mathbf{u}^T\|^p_{p,\omega^{-}_{\eta}})+c\al\eta+c\al.
\end{equation}
Let $y(t)=\int_{\Omega(t)}|\nabla\mathbf{u}^T|^pdx$
and $z(\eta)=\int^{\eta}_{\eta-1}y(t)dt$. Then
\begin{equation}
z(\eta)\leqslant c_6(z'(\eta)+z'(\eta)^{\frac{3}{p}})+c_7\al\eta+c_8\al.
\end{equation}
To apply Lemma 3.3 (i), we set $\Psi(\tau)=c_6(\tau+\tau^{\frac{3}{p}})$,
$\dl=\f12$, $t_0=1$,
and $\varphi(\eta)=2c_7 \al\eta+2c_8\al$, where $c_8$ satisfies
\begin{equation}
c_7\al+c_8\al \geqslant \Psi(2c_7\al).
\end{equation}
In addition, it follows from the proof of Theorem 4.1 that
\begin{equation}
z(T)\leqslant \varphi(T).
\end{equation}
Therefore, according (6.22)-(6.24) and Lemma 3.3 (i), we arrive at
\begin{equation}
y(\eta-1)\leqslant z(\eta)\leqslant \vp(\eta),
\end{equation}
which means  that for any $t\in[1,T]$,
\begin{equation}
y(t)\leqslant 2c_7 \al(t+1)+2c_8\al.
\end{equation}
\subsection{Part 2. Existence of solution $\mathbf{v}$ to problem (2.4)}

\vskip 0.1 true cm

{\bf Theorem 6.1.} {\it Let $\mu_0>0$ and $p>1$ or $\mu_0=0$ and $p>2$. Then problem (2.4) at least has a weak solution. }

\vskip 0.1 true cm

{\bf Remark 6.1.} {\it Here point out that in the case of $\mu_0=0$ and $p>2$, Theorem 6.1 has been proved in \cite{Gil-mar}
by different methods.}

\vskip 0.1 true cm
{\bf Proof.}
We just only treat the case of $\mu_0=0$ and $p>2$, the treatment for $\mu_0>0$ and $p>1$ is similar. Let $T=k$ and $\mathbf{u}^k=0$ in $\Omega\backslash\Omega_k$.
By (6.25) and a diagonalization process, we obtain a subsequence $\{\mathbf{u}^k\}$,
which is still denoted by $\{\mathbf{u}^k\}$, and a vector filed $\mathbf{u}\in W_{loc}^{1,p}(\bar{\Omega})$ such that
for any $t>0$,
\begin{equation}
\begin{split}
&\mathbf{u}^k\rightharpoonup\mathbf{u}\quad \text{in $W^{1,p}(\Omega(t))$},\\
&\mathbf{u}^k\rightarrow\mathbf{u}\quad  \text{in  $L^{p}(\Omega(t))$}.
\end{split}
\end{equation}
Next we show that for all compact subset $K\subset\subset\Omega$,
\begin{equation}
|\mathcal{D}(\mathbf{u}^k)+\mathcal{D}(\mathbf{a})|^{p-2}(\mathcal{D}(\mathbf{u}^k)+\mathcal{D}(\mathbf{a}))
\rightharpoonup|\mathcal{D}(\mathbf{u})+\mathcal{D}(\mathbf{a})|^{p-2}(\mathcal{D}(\mathbf{u})
+\mathcal{D}(\mathbf{a}))\quad \text{in $L^{p'}(K)$}.
\end{equation}
In fact, from Theorem 5.1 and (6.25), we have that there is a constant $c(K)$ which is independent of $k$
\begin{equation}
\|\nabla\mathbf{u}^k\|_{\kappa,p,K}\leqslant c(K).
\end{equation}
Thus, by the Rellich-Kondrachov theorem we have that
there exist a subsequence, which we still denote $\{\mathbf{u}^k\}$, and $\bar{\mathbf{u}}\in W_{loc}^{1,p}(\Omega)$ such that
\begin{equation}
\mathbf{u}^k\rightarrow \bar{\mathbf{u}}\quad\text{in $W^{1,p}(K)$}.
\end{equation}
Hence,
\begin{equation}
\nabla\mathbf{u}^k\rightarrow \nabla\bar{\mathbf{u}}\quad\text{in $L^p(K)$}.
\end{equation}
On the other hand,
\begin{equation}
\nabla\mathbf{u}^k\rightharpoonup \nabla\mathbf{u}\quad\text{ in $L^p(K)$}.
\end{equation}
Hence we have $\nabla\bar{\mathbf{u}}=\nabla\mathbf{u}$ and
\begin{equation}
\mathcal{D}(\mathbf{u}^k)\rightarrow\mathcal{D}(\mathbf{u})\quad\text{a.e. in K}.
\end{equation}
Note that
\begin{equation}
\||\mathcal{D}{\mathbf{u}^k}+\mathcal{D}{\mathbf{a}}|^{p-2}(\mathcal{D}{\mathbf{u}^k}+\mathcal{D}{\mathbf{a}})\|_{p',K}
\leqslant\|\mathcal{D}{\mathbf{u}}^k+\mathcal{D}\mathbf{a}\|^{p-1}_{p,K}\leqslant c(K).
\end{equation}
Therefore, according to Lemma I.1.3 of \cite{Lions}, we have
\begin{equation}
|\mathcal{D}{\mathbf{u}^k}+\mathcal{D}{\mathbf{a}}|^{p-2}(\mathcal{D}{\mathbf{u}^k}+\mathcal{D}{\mathbf{a}})
\rightharpoonup|\mathcal{D}{\mathbf{u}}+\mathcal{D}{\mathbf{a}}|^{p-2}(\mathcal{D}{\mathbf{u}}+\mathcal{D}{\mathbf{a}})\quad\text{ in $L^{p'}(K)$}.
\end{equation}
From (6.26) and (6.34), let $k\rightarrow\infty$, one derives that for all $\boldsymbol{\psi}\in\mcD(\Omega)$
\begin{multline}
\mu_0(\mcD(\mathbf{u})+\mcD(\mathbf{a}),\mcD(\boldsymbol{\psi}))+\mu_1(|\mcD(\mathbf{u})
+\mcD(\mathbf{a})|^{p-2}(\mcD(\mathbf{u})+\mcD(\mathbf{a})),\mcD(\boldsymbol{\psi}))\\
=(\mathbf{u}\cdot\nabla\boldsymbol{\psi},\mathbf{u})+(\mathbf{u}\cdot\nabla\boldsymbol{\psi},\mathbf{a})
+(\mathbf{a}\cdot\nabla\boldsymbol{\psi},\mathbf{u})+(\mathbf{a}\cdot\nabla\boldsymbol{\psi},\mathbf{a}).
\end{multline}
This means that $\mathbf{u}$ is a weak solution of problem (2.4).
\subsection{Part 3. Uniqueness}
At first, as in \cite{La-5}, we show that the dissipation of
the solution $\mathbf{v}$ to Problem (2.4) is distributed uniformly along $\Omega$.

{\bf Lemma 6.1.} {\it
Assume that $\mathbf{v}$ is a solution of problem (2.4). Then there exists a fixed constant
$c_9>0$ such that  if $\mu_0>0$, $p>1$,
\begin{equation}
\int_{\tau}^{\tau+1}dx_1\int_{\Sigma(x_1)}(|\nabla\mathbf{v}^T|^2+|\nabla\mathbf{v}^T|^p)dy\leqslant c_9\al;
\end{equation}
and
if $\mu_0=0$, $p>2$,
\begin{equation}
\int_{\tau}^{\tau+1}dx_1\int_{\Sigma(x_1)}|\nabla\mathbf{v}^T|^pdy\leqslant c_9\al.
\end{equation}}
{\bf Proof.}
Let
\begin{equation*}
\Omega_{\tau}(t)=\{x\in\Omega: \tau-t<x_1<\tau+t\},
\end{equation*}
\begin{align*}
y_{\tau}(t)&=\int_{\Omega_{\tau}(t)}(|\nabla\mathbf{u}^T|^2+|\nabla\mathbf{u}^T|^p)dx,& \text{if $\mu_0>0$, $p>1$,}\\
y_{\tau}(t)&=\int_{\Omega_{\tau}(t)}|\nabla\mathbf{u}^T|^pdx,&\text{if $\mu_0=0$, $p>2$,}
\end{align*}
and
\begin{equation*}
z_{\tau}(\eta)=\int\limits_{\eta-1}^{\eta}y_{\tau}(t)dt \qquad \text{for $\eta\geqslant1$.}
\end{equation*}
Since
\begin{equation*}
z_{\tau}(\tau)\leqslant y_{\tau}(\tau)\leqslant\varphi(2\tau+1),
\end{equation*}
similarly to the proof of (6.15) or (6.25), we have that for any $\eta\in[1,\tau]$,
\begin{equation*}
z_{\tau}(\eta)\leqslant \varphi(2\eta+1),
\end{equation*}
where the definition of $\vp(\eta)$ is given in (6.11) or (6.22). Therefore,
\begin{equation*}
y_{\tau}(\frac{1}{2})\leqslant\varphi(2)\leqslant c\al,
\end{equation*}
while, by Lemma 3.5, the same inequalities also hold for the vector $\mathbf{a}$, hence Lemma 6.1 is proved. $\square$

The following result comes from Lemma 4.2 of \cite{Mar-pa} (one can see also Proposition 4.3 of \cite{Gil-mar}).

{\bf Lemma 6.2.}
{\it Assume that $\mathbf{v}$ is a divergence free vector field in $W^{1,p}_{loc}(\bar{\Omega})$
vanishing on $\partial\Omega$ and satisfying (6.34) for $\mu_0=0$ and $p>2$.
If
\begin{equation}
\label{uniqueness condition}
|\mcD(\mathbf{v})(x_1,y)|\geqslant c|y|^{\frac{1}{p-1}}
\end{equation}
holds for some positive number $c>0$, then there is a constant $C>0$ such that
for all $\mathbf{w}\in\mathcal{D}^{1,p}_{loc}(\Omega)$ and $t>0$,
\begin{equation}
|(\mathbf{w}\cdot\nabla\mathbf{v},\mathbf{w})_{\Omega(t)}|\leqslant
C\alpha^{\f{1}{p}} \||\mathcal{D}(\mathbf{v})|^{\frac{1}{p-1}}\mcD(\mathbf{w})\|^2_{2,\Omega(t)},
\end{equation}
where the number $\al$
has been defined in Lemma 3.5.
}

{\bf Theorem 6.2.} {\it Assume that the flux $|\alpha_i|$ $(1\le i\le N)$ is sufficient small  and
there exists a solution $\mathbf{v}$ satisfying \eqref{uniqueness condition} to problem (2.4).
Then the solution $\mathbf{v}$ of problem (2.4) is unique for $p>2$ and $\mu_0=0$.
When $\mu_0>0$, even if the assumption \eqref{uniqueness condition} on $\mathbf{v}$ is removed,
then the solution $\mathbf{v}$ of problem (2.4) exists uniquely for $p>1$.
}

\vskip 0.1 true cm

{\bf Remark 6.2.} {\it Here point out that in the case of $\mu_0=0$ and $p>2$, Theorem 6.2 has been proved in \cite{Gil-mar}.}

\vskip 0.1 true cm

{\bf Proof.}  Assume that $\mathbf{v}_1$ and $\mathbf{v}_2$ are the solutions of problem (2.4).
Let $\mathbf{w}=\mathbf{v}_1-\mathbf{v}_2$. Then
\begin{multline}
-\mu_0div(\mcD(\mathbf{w}))-\mu_1div\Big{(}|\mcD(\mathbf{w})+\mcD(\mathbf{v}_2)|^{p-2}\big{(}\mcD(\mathbf{w})+\mcD(\mathbf{v}_2)\big{)}
-|\mcD(\mathbf{v}_2)|^{p-2}\mcD(\mathbf{v}_2)\Big{)}\\
+\mathbf{w}\cdot\nabla\mathbf{w}+\mathbf{w}\cdot\nabla\mathbf{v}_2+\mathbf{v}_2\cdot\nabla\mathbf{w}+\nabla(\pi_1-\pi_2)=0,
\end{multline}
From (6.40) and the proof of (6.15) and (6.25), we just need to estimate $(\mathbf{w}\cdot\nabla\mathbf{v}_2,\mathbf{w})_{\Omega(t)}$.

If $\mu_0=0$, $p>2$, since $\mathbf{v}_2$ satisfies \eqref{uniqueness condition}, then according to Lemma 6.2, we have
\begin{equation}
\label{uniqueness 1}
(\mathbf{w}\cdot\nabla\mathbf{v}_2,\mathbf{w})_{\Omega(t)}\leqslant C\alpha^{\f{1}{p}} \||\mathcal{D}(\mathbf{v}_2)|^{\frac{p-2}{2}}
\mcD(\mathbf{w})\|^2_{2,\Omega(t)}.
\end{equation}

If $\mu_0>0$, without loss of generality, we set $t=k$, by (6.36) we can get that
\begin{multline}
\label{uniqueness 2}
|(\mathbf{w}\cdot\nabla\mathbf{v}_2,\mathbf{w})_{\Omega(k)}|
=|\sum\limits_{j=-k}^{k-1}\int_{j}^{j+1}\int_{\Sigma(x_1)}
\mathbf{w}\cdot\nabla\mathbf{v}_2\cdot\mathbf{w}dx|\\
\leqslant\sum\limits_{j=-k}^{k-1}\big{(}\int_{j}^{j+1}\int_{\Sigma(x_1)}|\nabla\mathbf{v}_2|^2dx\big{)}^{1/2}\big{(}\int_{j}^{j+1}
\int_{\Sigma(x_1)}|\mathbf{w}|^4dx\big{)}^{1/2}
\leqslant C\alpha^{\f{1}{2}} \|\nabla\mathbf{w}\|^2_{2,\Omega(t)}.
\end{multline}
By \eqref{uniqueness 1} and \eqref{uniqueness 2}, similarly to the proof of (6.11) or (6.21), we have
\begin{equation}
z(\eta)\leqslant\Psi(z'(\eta)), \end{equation}
where
\begin{equation*}
\text{$z(\eta)=\int_{\eta-1}^{\eta}y(t)dt$\quad with\quad
$y(t)=\int_{\Omega(t)}(|\nabla\mathbf{w}^T|^2+|\nabla\mathbf{w}^T|^p)dx$, \ if $\mu_0>0$, $p>1$},
\end{equation*}
and
\begin{equation*}
\text{$z(\eta)=\int_{\eta-1}^{\eta}y(t)dt$\quad with \quad$y(t)=\int_{\Omega(t)}|\nabla\mathbf{w}^T|^pdx$, \quad if $\mu_0=0$, $p>2$},
\end{equation*}
and the definition of function $\Psi(\tau)$ is given in (6.12) or (6.22).
If $z(t)$ is not identically zero, it then follows from Lemma 3.3 (iii) that when $\mu_0=0$,
\begin{align*}
\liminf\limits_{t\rightarrow\infty}t^{\frac{-3}{3-p}}z(t)>0\qquad \text{if $p<3$},\\
\liminf\limits_{t\rightarrow\infty}e^{-t}z(t)>0\qquad \text{if $p\geqslant3$};
\end{align*}
when $\mu_0>0$,
\begin{equation*}
\liminf\limits_{t\rightarrow\infty}t^{-3}z(t)>0.
\end{equation*}
These contradict with \eqref{problem 1}$_5$. Hence, $z\equiv0$, and further $\mathbf{v}_1=\mathbf{v}_2$.
This completes the proof of Theorem 6.1. \qquad\qquad \qquad \qquad \qquad \qquad \qquad \qquad \qquad
\qquad \qquad \qquad \qquad \qquad \qquad \qquad \qquad \qquad \qquad $\square$

{\bf Remark 6.3.}  {\it  When $\Omega$ admits straight outlets, then the corresponding solution $\mathbf{v}$ of Leray Problem (1.8)
exactly satisfies (6.38) (see \cite{Mar-pa}). Therefore, in this case,
if the flux is sufficiently small, then the solution of problem (2.4) coincides with the solution of Leray Problem (1.8). }

\section{Solvability of Ladyzhenskaya-Solonnikov Problem II (2.5)}

To solve  Ladyzhenskaya-Solonnikov Problem II (2.5) under some suitable conditions, based on Sections 4-5,
we will take the following two parts.

\subsection{Part 1. Existence of solution $\mathbf{v}$ to problem (2.5)}
In the following, we will assume that
\begin{equation}
\Omega_i=\{x^i=(x_1^i, y^i): 0<x_1^i<\infty,|y^i|<g_i(x_1^i)\}.
\end{equation}
Hence it follows from the definition of $I_i(t)$ in (2.5) that
\begin{equation}
\begin{split}
I_i(t)& =\int_0^tg_i^{-(1+d)}(s)ds+\int_0^tg_i^{(1-p)d-1}(s)ds, \quad \text{if $\mu_0>0$;}\\
I_i(t)& =\int_0^tg_i^{(1-p)d-1}(s)ds, \qquad \qquad \qquad \qquad \quad\ \text{if $\mu_0=0$.}
\end{split}
\end{equation}
As in \cite{La-5}, we suppose that some outlets are ``narrow", namely,
\begin{equation}
I_i(\infty)=\infty, \  i=1,...,m,
\end{equation}
meanwhile other outlets are ``wide", that is,
\begin{equation}
I_i(\infty)<\infty, \  i=m+1,...,N.
\end{equation}
In addition, we assume that
\begin{equation}
\text{$|g'_i(t)|\leqslant(2\beta)^{-1}$ and $g_i(t)\geqslant g_0>0$}.
\end{equation}
In this case, $g_i(t)$ has the properties as follows
\begin{align*}
&g_i(t)\leqslant  g_i(0)+(2\beta)^{-1}t,\\
&t-\beta g_i(t)\geqslant\frac{1}{2}t-\beta g_i(0),\\
t&-\beta g_i(t)\geqslant 0\quad\text{for $t\geqslant t^{*}=2\beta \max\limits_{i=1,...,m}g_i(0)$},\\
\frac{1}{2}g_i(t) &\leqslant g(s)\leqslant\frac{3}{2}g_i(t)\quad\text{for $s\in[t-\beta g_i(t),t]$ and $t\geqslant t^{*}$}.
\end{align*}
Let
\begin{equation*}
\Omega(t)=\Omega_0\cup\big{(}\ds\cup_{i=1}^m\Omega_i(h_i(t))\big{)}\cup(\ds\cup_{i=m+1}^{N}\Omega_i),
\end{equation*}
where  $\Omega_i(s)=\{x\in\Omega_i:0<x_1^i<s\}$ for $1\le i\le N$, and the function $h_i(t)$
$(1\le i\le m)$ is determined by the equations
\begin{equation}
\begin{split}
\text{$h'_i(t)=g_i^{\f{7-d}{3}}(h_i(t))$, \qquad if $\mu_0>0$;}\\
\text{$h'_i(t)=g_i^{\f{2p+3}{3}-(1-\f{p}{3})d}(h_i(t))$, \qquad if $\mu_0=0$}.
\end{split}
\end{equation}
Actually, $h_i(t)$ ($1\le i\le m$) is the inverse function  of
$k_i(t)=\int_0^t g_i^{\f{d-7}{3}}(s)ds$ (for the case of $\mu_0>0$) or $k_i(t)=\int_0^t g_i^{-\f{2p+3}{3}+(1-\f{p}{3})d}(s)ds$
(for the case of $\mu_0=0$)
for $t>0$. It is worth noting that we need $h_i(t)\rightarrow\infty$ as $t\rightarrow\infty$. If $p\geqslant\f{4+2d}{3d}$, since $\f{7-d}{3}\leqslant(p-1)d+1$ and $I_i(\infty)=\infty$,
automatically, we have that $h_i(t)\rightarrow\infty$ as $t\rightarrow\infty$. While, if $1<p<\f{4+2d}{3d}$, we will assume that
\begin{equation}
\begin{split}
\text{$h_i(t)\rightarrow\infty$ as $t\rightarrow\infty$},\\
\text{or $k_i(t)=\int_0^t g_i^{\f{d-7}{3}}(s)ds\rightarrow\infty$ as $t\rightarrow\infty$}.
\end{split}
\end{equation}
Since $h_i(t)\rightarrow\infty$ as $t\rightarrow\infty$, we can introduce the truncated domains
\begin{equation*}
\Omega(t,R)=\Omega_0\cup(\cup_{i=1}^m\Omega_i(h_i(t)))\cup(\cup_{i=m+1}^N\Omega_i(R))
\end{equation*}
and the ``truncating" functions $\zeta(x,t)$
\begin{equation}
\zeta(x,t)=
\begin{cases}
\beta &\text{for $x\in\Omega_0\cup\Big{(}\cup_{i=1}^m\Omega_i\big{(}h_i(t)-\beta g_i(h(t))\big{)}\Big{)}\cup(\cup_{i=m+1}^N\Omega_i)$},\\
 0 &\text{for $x\in\cup_{i=1}^m(\Omega_i\backslash\Omega_i(h_i(t)))$},\\
\frac{h_i(t)-x_1^i}{g_i(h_i(t))}&\text{ for $x\in\omega_i(t)$},
\end{cases}
\end{equation}
where  $\omega_i(t)=\{x\in\Omega:h_i(t)-\beta g_i(h_i(t))<x_1^i<h_i(t)\}$ for $1\le i\le m$.
It is easy to check that \begin{equation}
|\nabla_x\zeta|\leqslant[g_i(h_i(t))]^{-1}
\end{equation}
and
\begin{equation}
\frac{\partial\zeta(x,t)}{\partial t}=\frac{h'_i(t)}{g_i(h_i(t))}[1-\frac{h_i(t)-x_1^i}{g_i(h_i(t))}g'_i(h_i(t))]
\geqslant\frac{1}{2}\frac{h'(t)}{g_i(h_i(t))}.
\end{equation}
This yields
\begin{equation}
\begin{split}
&\text{$\frac{\partial\zeta(x,t)}{\partial t}\geqslant\frac{1}{2}g_i^{\f{4-d}{3}}(h_i(t))$,\qquad\qquad\qquad\qquad\ \ \ if $\mu_0>0$;}\\
&\text{$\frac{\partial\zeta(x,t)}{\partial t}\geqslant\frac{1}{2}g^{\f{2p}{3}-(1-\f{p}{3})d}_i(h_i(t))$,\qquad\qquad\qquad if $\mu_0=0$}.
\end{split}
\end{equation}
We now set
\begin{equation}
\begin{split}
&\text{$Q(t)=\int_{\Omega(t)}(|\nabla\mathbf{v}|^2+|\nabla\mathbf{v}|^p)\zeta(x,t)dx$, \qquad\quad \ \ if $\mu_0>0$;}\\
&\text{$Q(t)=\int_{\Omega(t)}|\nabla\mathbf{v}|^p\zeta(x,t)dx$,\qquad\qquad\quad\qquad\quad  if $\mu_0=0$,}
\end{split}
\end{equation}
and
\begin{equation}
\begin{split}
&\text{$\phi(t)=\sum\limits_{i=1}^m \int\limits_0^{h_i(t)}(g^{-(d+1)}(s)+g^{(1-p)d-1}(s))ds$, \qquad if $\mu_0>0$;}\\
&\text{$\phi(t)=\sum\limits_{i=1}^m \int\limits_0^{h_i(t)}g^{(1-p)d-1}(s)ds$, \qquad\qquad\qquad\qquad\ if $\mu_0=0$.}
\end{split}
\end{equation}
Then we have

{\bf Theorem 7.1.} {\it  if $\mu_0>0$, $p\geqslant\frac{4+2d}{3d}$ or $\mu_0=0$, $2< p\leqslant 3-\f{2}{d}$,
under the assumptions (7.3)-(7.5) and (7.7), there exists at least one solution $\mathbf{u}$  satisfying $(2.5)_1-(2.5)_4$, moreover,
\begin{equation}
Q(t)\leqslant \bar{c}_{10}\al\phi(t)+\bar{c}_{11}\al,\quad\text{for any $t\geqslant \bar{t}=\max\limits_{i=1,...,m}k_i(t^{*})$,}
\end{equation}
where $\bar{c}_{10}>0$ and $\bar{c}_{11}>0$ are some constants, and the definitions of $Q(t)$ and $\phi(t)$ are given in
(7.12) and (7.13) respectively. When $1<p<\frac{4+2d}{3d}$, if  we assume that $g_i(t)$ satisfies the condition
\begin{equation}
|g'_i(t)g^{\f{4-pd}{2}}_i(t)|<\gamma,
\end{equation}
where $\gamma>0$ is a sufficient small constant,
then the corresponding conclusion (7.14) still holds.\\
}

{\bf Proof.}   Let $\mathbf{v}^{T,R}=\mathbf{u}^{T,R}+\mathbf{a}$ be the solution of problem (2.5) in the bounded domains $\Omega(T,R)$.
Note that the existence of $\mathbf{v}^{T,R}$ is ensured by Theorem 4.1. Moreover, since $|\mathbf{a}(x_1,x')|\leqslant cg_i^{-(d-1)}(x_1)$ and
$|\mathbf{a}(x_1,x')|\leqslant cg_i^{-d}(x_1)$ for $x\in\Omega_i$, if $\mu_0>0$, by (4.10) we have
\begin{equation*}
\|\nabla\mathbf{u}^{T,R}\|^2_{2,\Omega(T,R)}+\|\nabla\mathbf{u}^{T,R}\|^p_{p,\Omega(T,R)}\leqslant c\al\phi(T).
\end{equation*}
Let $R\rightarrow\infty$, we can find a vector function $\mathbf{u}^T$ such that $\mathbf{v}^T=\mathbf{u}^T+\mathbf{a}$ is a
solution of Problem (2.5) in $\Omega(T)$ and
\begin{equation*}
\|\nabla\mathbf{u}^T\|^2_{2,\Omega(T)}+\|\nabla\mathbf{u}^T\|^p_{p,\Omega(T)}\leqslant c\al\phi(T).
\end{equation*}
If $\mu_0=0$, by (4.29) we obtain that
\begin{equation*}
\|\nabla\mathbf{u}^T\|^p_{p,\Omega(T)}\leqslant c\big{(}\|\nabla\mathbf{a}\|^p_{p,\Omega(T)}
+\|\mathbf{a}\|^{\frac{p(d-2)}{(p-2)(d-1)}}_{\frac{p(d-2)}{(p-2)(d-1)},\Omega(T)}+\|\mathbf{a}\|^{2p'}_{2p',\Omega(T)}\big{)}.
\end{equation*}
Since $2< p \leqslant 3-\f{2}{d}$, a direct calculation derives
\begin{equation*}
\|\nabla\mathbf{a}\|^p_{p,\Omega(T)}+\|\mathbf{a}\|^{\frac{p(d-2)}{(p-2)(d-1)}}_{\frac{p(d-2)}{(p-2)(d-1)},\Omega(T)}
+\|\mathbf{a}\|^{2p'}_{2p',\Omega(T)}\leqslant c\al\phi(T).
\end{equation*}
Therefore,
\begin{equation*}
\|\nabla\mathbf{u}^T\|^p_{p,\Omega(T)}\leqslant c\al\phi(T).
\end{equation*}
Next we derive (7.14).
Multiplying the equation in \eqref{problem 2} by $\mathbf{u}(x)\zeta(x,t)$
and integrating over $\Omega(t)$ yield
\begin{equation}
\int_{\Omega(t)}(\mu_0\mcD(\mathbf{v})+\mu_1|\mcD(\mathbf{v})|^{p-2}\mcD(\mathbf{v})):\mcD(\zeta\mathbf{v})dx
+\int_{\Omega(t)}\mathbf{v}\cdot\nabla\mathbf{v}\cdot\zeta\mathbf{u}dx
+\int_{\Omega(t)}\nabla \pi\cdot\zeta\mathbf{u}dx=0, \end{equation}
where and below the superscript $T$ of $\mathbf{u}$ is omitted for notational convenience. It is easy to check that
\begin{equation*}
\int_{\Omega(t)}\mcD(\mathbf{v}):\mcD(\zeta\mathbf{u})dx=\int_{\Omega(t)}\zeta|\mcD(\mathbf{u})|^2dx
+\int_{\Omega(t)}\zeta\mcD(\mathbf{a}):\mcD(\mathbf{u})dx+\int_{\Omega(t)}\mcD(\mathbf{v}):Sym(\mathbf{u}\otimes\nabla\zeta)dx
\end{equation*}
and
\begin{multline*}
\int_{\Omega(t)}|\mcD(\mathbf{v})|^{p-2}\mcD(\mathbf{v}):\mcD(\zeta\mathbf{u})dx
=\int_{\Omega(t)}\zeta|\mcD(\mathbf{v})|^{p-2}\mcD(\mathbf{v}):\mcD(\mathbf{u})dx\\
+\int_{\Omega(t)}|\mcD(\mathbf{v})|^{p-2}\mcD(\mathbf{v}):Sym(\mathbf{u}\otimes\nabla\zeta)dx.
\end{multline*}
In addition, by integrating by parts we deduce
\begin{multline*}
\int_{\Omega(t)}\mathbf{v}\cdot\nabla\mathbf{v}\cdot\zeta\mathbf{u}dx
=\int_{\Omega(t)}\mathbf{v}\cdot\nabla\mathbf{u}\cdot\zeta\mathbf{u}dx
+\int_{\Omega(t)}\mathbf{v}\cdot\nabla\mathbf{a}\cdot\zeta\mathbf{u}dx\\
=-\frac{1}{2}\int_{\Omega(t)}\mathbf{u}^2\mathbf{v}\cdot\nabla\zeta dx
+\int_{\Omega(t)}\mathbf{v}\cdot\nabla\mathbf{a}\cdot\zeta\mathbf{u}dx
\end{multline*}
and
\begin{align*}
\int_{\Omega(t)}&\mathbf{v}\cdot\nabla\mathbf{a}\cdot\zeta\mathbf{u}dx
=\int_{\Omega(t)}\mathbf{u}\cdot\nabla\mathbf{a}\cdot\zeta\mathbf{u}dx
+\int_{\Omega(t)}\mathbf{a}\cdot\nabla\mathbf{a}\cdot\zeta\mathbf{u}dx\\
&=-\int_{\Omega(t)}(\mathbf{u}\cdot\nabla\zeta)(\mathbf{u}\cdot\mathbf{a})dx
-\int_{\Omega(t)}\mathbf{u}\cdot\nabla\mathbf{u}\cdot\mathbf{a}\zeta dx
-\int_{\Omega(t)}(\mathbf{a}\cdot\nabla\zeta)(\mathbf{u}\cdot\mathbf{a})dx
-\int_{\Omega(t)}\mathbf{a}\cdot\nabla\mathbf{u}\cdot\mathbf{a}\zeta dx.
\end{align*}
Meanwhile,
\[
\begin{split}
\int_{\Omega(t)}\zeta|\mcD(\mathbf{v})|^{p-2}\mcD(\mathbf{v}):\mcD(\mathbf{u}) dx-\int_{\Omega(t)}\zeta|\mcD(\mathbf{a})|^{p-2}\mcD(\mathbf{a}):\mcD(\mathbf{u}) dx
\geqslant \delta\int_{\Omega(t)}\zeta|\mcD(\mathbf{u})|^p dx,
\end{split}
\]
and
\[
\begin{split}
\Big{|}\int_{\Omega(t)}\zeta|\mcD(\mathbf{a})|^{p-2}\mcD(\mathbf{a}):\mcD(\mathbf{u}) dx\Big{|}
\leqslant c\int_{\Omega(t)}\zeta |\mcD(\mathbf{a})|^p dx
+\ve\int_{\Omega(t)}\zeta |\mcD(\mathbf{u})|^pdx,
\end{split}
\]
which yields
\begin{equation*}
\int_{\Omega(t)}\zeta|\mcD(\mathbf{v})|^{p-2}\mcD(\mathbf{v}):\mcD(\mathbf{u}) dx\geqslant c\int_{\Omega(t)}\zeta|\mcD(\mathbf{u})|^p dx
-\f{\delta}{2}\int_{\Omega(t)}\zeta |\mcD(\mathbf{a})|^p dx.
\end{equation*}
Hence,
\begin{equation}
\mu_0\int_{\Omega(t)}\zeta|\mcD(\mathbf{u})|^2 dx
+\f{\delta}{2}\mu_1\int_{\Omega(t)}\zeta|\mcD(\mathbf{u})|^p dx\leqslant c\mu_1\int_{\Omega(t)}\zeta |\mcD(\mathbf{a})|^pdx+I_1+...+I_9,
\end{equation}
where
\[
\begin{split}
&I_1=-\mu_0\int_{\Omega(t)}\zeta\mcD(\mathbf{u}):\mcD(\mathbf{a})dx,\\
&I_2=-\mu_0\int_{\Omega(t)}\mcD(\mathbf{v}):Sym(\mathbf{u}\otimes\nabla\zeta)dx,\\
&I_3=-\mu_1\int_{\Omega(t)}|\mcD(\mathbf{v})|^{p-2}\mcD(\mathbf{v}):Sym(\mathbf{u}\otimes\nabla\zeta)dx,\\
&I_4=\frac{1}{2}\int_{\Omega(t)}\mathbf{u}^2\mathbf{v}\cdot\nabla\zeta  dx,\\
&I_5=\int_{\Omega(t)}(\mathbf{u}\cdot\nabla\zeta)(\mathbf{u}\cdot\mathbf{a}) dx,\\
&I_6=\int_{\Omega(t)}\mathbf{u}\cdot\nabla\mathbf{u}\cdot\mathbf{a}\zeta dx,\\
&I_7=\int_{\Omega(t)}(\mathbf{a}\cdot\nabla\zeta)(\mathbf{u}\cdot\mathbf{a}) dx,\\
&I_8=\int_{\Omega(t)}\mathbf{a}\cdot\nabla\mathbf{u}\cdot\mathbf{a}\zeta dx,\\
&I_9=\int_{\Omega(t)}\pi\mathbf{u}\cdot\nabla\zeta dx.\\
\end{split}
\]
Next, we treat the terms $I_i$ $(1\le i\le 9)$ in two cases as follows:\\
\vskip 0.1 true cm

$\mathbf{Case\ I}$. {\bf $\mu_0>0$, $p>1$}\\

\vskip 0.1 true cm
It is easy to check that
\begin{equation*}
|I_1|\leqslant\ve\int_{\Omega(t)}\zeta|\mcD(\mathbf{u})|^2dx+c\int_{\Omega(t)}\zeta|\mcD(\mathbf{a})|^2dx.
\end{equation*}
By H\"{o}lder and Poincar\'{e} inequality, we have
\begin{align*}
|I_2|&\leqslant\ve\int_{\Omega(t)}|\mcD(\mathbf{v})|^2dx+\int_{\Omega(t)}||\mathbf{u}|\nabla\zeta|^2dx\\
&\leqslant2\ve\int_{\Omega(t)}|\mcD(\mathbf{u})|^2dx+ \int_{\omega_i(t)}|\mcD(\mathbf{u})|^2dx+2\ve\int_{\Omega(t)}|\mcD(\mathbf{a})|^2dx.
\end{align*}
Similarly,
\begin{align*}
|I_3|&\leqslant\ve\int_{\Omega(t)}|\mcD(\mathbf{v})|^pdx+\int_{\Omega(t)}||\mathbf{u}|\nabla\zeta|^pdx\\
&\leqslant2^{p-1}\ve\int_{\Omega(t)}|\mcD(\mathbf{u})|^pdx+\int_{\omega_i(t)}|\mcD(\mathbf{u})|^pdx+2^{p-1}\ve\int_{\Omega(t)}|\mcD(\mathbf{a})|^pdx.
\end{align*}
Using H\"{o}lder inequality and Lemma 3.1, we deduce that
\[
\begin{split}
&|I_4|=|\f{1}{2}\int_{\Omega(t)}\mathbf{u}^2\mathbf{u}\cdot\nabla\zeta dx
+\f{1}{2}\int_{\Omega(t)}\mathbf{u}^2\mathbf{a}\cdot\nabla\zeta dx|\\
&\leqslant c\sum\limits_{i=1}^{m}g^{-1}_i(h_i(t))\|\mathbf{u}\|^3_{3,\omega_i(t)}+c\sum\limits_{i=1}^{m}g^{-d}_i(h_i(t))\|\mathbf{u}\|^2_{2,\omega_i(t)}\\
&\leqslant c\sum\limits_{i=1}^{m}g_i^{\frac{4-d}{2}}(h_i(t))\|\nabla\mathbf{u}\|^3_{2,\omega_i(t)}
+c\sum\limits_{i=1}^{m}g_i^{2-d}(h_i(t))\|\nabla\mathbf{u}\|^2_{2,\omega_i(t)}.
\end{split}
\]
By Poincar\'{e} inequality, we have
\[
\begin{split}
|I_5|\leqslant \sum\limits_{i=1}^{m}g^{-d}_i(h_i(t))\|\mathbf{u}\|^2_{2,\omega_i(t)}
\leqslant c\sum\limits_{i=1}^{m} g_i^{2-d}(h_i(t))\|\nabla\mathbf{u}\|^2_{2,\omega_i(t)}.
\end{split}
\]
It follows from Lemma 3.1 and Lemma 3.5 (iii) and (iv) that
\[
\begin{split}
&|I_6|
\leqslant\int_{\Omega(t)}\mathbf{a}^2\mathbf{u}^2\zeta dx+\ve\int_{\Omega(t)}|\nabla\mathbf{u}|^2\zeta dx\\
&\leqslant \beta\int_{\Omega_0\cup\big{(}\cup_{i=1}^m\Omega_i(h_i(t)-\beta g_i(h(t)))\big{)}\cup(\cup_{i=m+1}^N\Omega_i)}\mathbf{a}^2\mathbf{u}^2dx
+\beta\sum\limits_{i=1}^m\int_{\omega_i(t)}\mathbf{a}^2\mathbf{u}^2dx+\ve\int_{\Omega(t)}|\nabla\mathbf{u}|^2\zeta dx\\
&\leqslant 2\ve \int_{\Omega(t)}|\nabla\mathbf{u}|^2\zeta dx+\ve\sum\limits_{i=1}^m\int_{\omega_i(t)}|\nabla\mathbf{u}|^2 dx.
\end{split}
\]
By Schwartz inequality and Lemma 3.1, we arrive at
\begin{align*}
|I_7|&\leqslant\sum\limits_{i=1}^m\||\mathbf{u}|\nabla\zeta\|^2_{2,\omega_i(t)}+\sum\limits_{i=1}^m\|\mathbf{a}\|^2_{4,\omega_i(t)}\\
&\leqslant c\sum\limits_{i=1}^m \|\nabla\mathbf{u}\|^2_{2,\omega_i(t)}
+\sum\limits_{i=1}^m\|\mathbf{a}\|^4_{4,\omega_i(t)}.
\end{align*}
On the other hand, it is easy to get
\begin{equation*}
|I_8|\leqslant \ve\|\zeta\nabla\mathbf{u}\|^2_{2,\Omega(t)}+\|\mathbf{a}\|^4_{4,\Omega(t)}.
\end{equation*}
Finally, we estimate $I_9$.
By Lemma 3.2, we can find $\mathbf{w}\in W_0^{1,2}(\omega_i(t))\cap W_0^{1,p}(\omega_i(t))$
such that $div\mathbf{w}=u_1$ (see Remark 3.2),
and there is a constant $M(d,p)>0$ (see Remark 3.1) such that
\begin{equation}
\begin{split}
&\|\nabla\mathbf{w}\|_{2,\omega_i(t)}\leqslant M(d,p)\|u_1\|_{2,\omega_i(t)},\\
&\|\nabla\mathbf{w}\|_{p,\omega_i(t)}\leqslant M(d,p)\|u_1\|_{p,\omega_i(t)}.
\end{split}
\end{equation}
In this case,
\begin{equation}
\begin{split}
&I_9=g_i^{-1}(h_i(t))\int_{\omega_i(t)}\pi u_1 dx\\
&=-g_i^{-1}(h_i(t))\int_{\omega_i(t)}\mathbf{w}\cdot\nabla \pi dx\\
&=g^{-1}(h_i(t))\int_{\omega_i(t)}[-div(\mu_0\mcD(\mathbf{v})
+\mu_1|\mcD(\mathbf{v})|^{p-2}\mcD(\mathbf{v}))\cdot\mathbf{w}
+\mathbf{v}\cdot\nabla\mathbf{v}\cdot\mathbf{w}] dx\\
&=g^{-1}(h_i(t))\int_{\omega_i(t)}[(\mu_0\mcD(\mathbf{v})
+\mu_1|\mcD(\mathbf{v})|^{p-2}\mcD(\mathbf{v}))\cdot\nabla\mathbf{w}
+\mathbf{v}\cdot\nabla\mathbf{v}\cdot\mathbf{w}] dx.
\end{split}
\end{equation}
Note that by (7.18), one has
\begin{equation*}
\Big{|}\int_{\Omega_i(t)}\mcD(\mathbf{v})\cdot\nabla\mathbf{w} dx\Big{|}\leqslant\|\mcD(\mathbf{v})\|_{2,\omega_i(t)}\|u_1\|_{2,\omega_i(t)}
\end{equation*}
and
\begin{equation*}
\Big{|}\int_{\omega_i(t)}|\mathcal{D}(\mathbf{v})|^{p-2}\mathcal{D}(\mathbf{v})\cdot\nabla\mathbf{w} dx\Big{|}
\leqslant\|\mcD(\mathbf{v})\|^{p-1}_{p,\omega_i(t)}\|\mcD(\mathbf{w})\|_{p,\omega_i(t)}
\leqslant\|\mcD(\mathbf{v})\|^{p-1}_{p,\omega_i(t)}\|u_1\|_{p,\omega_i(t)}.
\end{equation*}
In addition, by Lemma 3.1,
\begin{equation}
\text{$\|u_1\|_{q,\omega_i(t)}\leqslant cg_i(h_i(t))\|\nabla\mathbf{u}\|_{q,\omega_i(t)}$\quad  for any $q>1$.}
\end{equation}
Consequently, we have
\begin{multline}
\Big{|}g^{-1}(h_i(t))\int_{\omega_i(t)}(\mu_0\mcD(\mathbf{v})+\mu_1|\mcD(\mathbf{v})|^{p-2}\mcD(\mathbf{v}))\cdot\nabla\mathbf{w} dx\Big{|}\\
\leqslant c\|\nabla\mathbf{u}\|^2_{2,\omega_i(t)}+c\|\nabla\mathbf{u}\|^p_{p,\omega_i(t)}+c\|\nabla\mathbf{a}\|^2_{2,\omega_i(t)}
+c\|\nabla\mathbf{a}\|^p_{p,\omega_i(t)}.
\end{multline}
On the other hand, by Lemma 3.1 and (7.18),
\begin{equation}
\begin{split}
&\Big{|}\int_{\omega_i(t)}\mathbf{v}\cdot\nabla\mathbf{v}\cdot\mathbf{w}dx\Big{|}
=\Big{|}\int_{\omega_i(t)}\mathbf{v}\cdot\nabla\mathbf{w}\cdot\mathbf{v}dx\Big{|}\\
&\leqslant2\|\mathbf{u}\|^2_{4,\omega_i(t)}\|\nabla\mathbf{w}\|_{2,\omega_i(t)}+2\|\mathbf{a}\|^2_{4,\omega_i(t)}\|\nabla\mathbf{w}\|_{2,\omega_i(t)}\\
&\leqslant cg_i^{\frac{4-d}{2}}(h_i(t))\|\nabla\mathbf{u}\|^2_{2,\omega_i(t)}\|u_1\|_{2,\omega_i(t)}+c\|\mathbf{a}\|^2_{4,\omega_i(t)}\|u_1\|_{2,\omega_i(t)}.
\end{split}
\end{equation}
Therefore,
\begin{equation}
\begin{split}
|I_9|&=\Big{|}\sum\limits_{i=1}^m\int_{\Omega_i(t)}\pi\mathbf{u}\cdot\nabla\zeta dx\Big{|}\leqslant c\sum\limits_{i=1}^m\Big{(}\|\nabla\mathbf{u}\|^2_{2,\omega_i(t)}+\|\nabla\mathbf{u}\|^p_{p,\omega_i(t)}\\
&\quad +g_i^{\f{4-d}{2}}(h_i(t))\|\nabla\mathbf{u}\|^3_{2,\omega_i(t)}+\|\nabla\mathbf{a}\|^2_{2,\omega_i(t)}+\|\nabla\mathbf{a}\|^4_{4,\omega_i(t)}
+\|\nabla\mathbf{a}\|^p_{p,\omega_i(t)}\Big{)}.
\end{split}
\end{equation}
Collecting the estimates on $I_i$ $(1\le i\le 9)$, together with Korn inequality and (7.17), we have
\begin{multline}
\|\nabla\mathbf{u}\|^2_{2,\Omega(t)}+\|\nabla\mathbf{u}\|^p_{p,\Omega(t)}
\leqslant c\sum\limits_{i=1}^m\big{(}\|\nabla\mathbf{u}\|^2_{2,\omega_i(t)}+
g_i^{\frac{4-d}{2}}(h_i(t))\|\nabla\mathbf{u}\|^3_{2,\omega_i(t)}+\|\nabla\mathbf{u}\|^p_{p,\omega_i(t)}\big{)}\\
+c\|\nabla\mathbf{a}\|^2_{2,\Omega(t)}+c\|\nabla\mathbf{a}\|^p_{p,\Omega(t)}+c\|\mathbf{a}\|^4_{4,\Omega(t)}.
\end{multline}
In addition, it is easy to get
\begin{equation}
\frac{dy(t)}{dt}\geqslant\sum\limits_{i=1}^m \f{1}{2}g_i^{\f{4-d}{3}}(h(t))\int_{\omega_i(t)}(|\nabla\mathbf{u}|^2+|\nabla\mathbf{u}|^p)dx,
\end{equation}
where $y(t)=\int_{\Omega(t)}(|\nabla\mathbf{u}|^2+|\nabla\mathbf{u}|^p)\zeta(x,t)dx$.
Hence
\begin{equation}
\sum\limits_{i=1}^m g_i^{\f{4-d}{2}}(h_i(t))\|\nabla\mathbf{u}\|_{\omega_i(t)}^3\leqslant c(\frac{dy(t)}{dt})^{3/2}
\end{equation}
and
\begin{equation}
\sum\limits_{i=1}^m(\|\nabla\mathbf{u}\|^2_{2,\omega_i(t)}+\|\nabla\mathbf{u}\|^p_{p,\omega_i(t)})\leqslant g^{\f{d-4}{3}}_0\frac{dy(t)}{dt}.
\end{equation}
Since
\begin{equation}
\|\nabla\mathbf{a}\|^2_{2,\Omega(t)}+\|\mathbf{a}\|^4_{4,\Omega(t)}+\|\nabla\mathbf{a}\|^p_{p,\Omega(t)}\leqslant c\al
\sum\limits_{i=1}^m\int\limits_0^{h_i(t)}(g_i^{-(d+1)}(s)+g_i^{(1-p)d-1}(s))ds,
\end{equation}
we have
\begin{equation}
y(t)\leqslant c[y'(t)+y'^{3/2}(t)]+c\al\phi(t).
\end{equation}
Note that
\begin{equation}
\frac{d\phi(t)}{dt}=\sum\limits_{i=1}^{m}\Big{(}\frac{h'_i(t)}{g^{d+1}_i(h_i(t))}+\frac{h'_i(t)}{g_i^{(p-1)d+1}(h_i(t))}\Big{)}
=\sum\limits_{i=1}^{m}\big{(}g_i^{\frac{4(1-d)}{3}}(h_i(t))+g_i^{\f{4}{3}+(\f{2}{3}-p)d}(h_i(t))\big{)}.
\end{equation}
If $p\geqslant\f{4+2d}{3d}$, we then have $\frac{d\phi(t)}{dt}\leqslant m(g_0^{\frac{4(1-d)}{3}}
+g_0^{\f{4}{3}+(\f{2}{3}-p)d})\equiv c_{12}$.
Thus, if we set $\Psi(t)=c_{13}(t+t^{\f{3}{2}})$, $\varphi(t)=2c_{10}\al\phi(t)+2c_{11}\al$ and $\delta=\frac{1}{2}$, $t_0=\bar{t}$, where $c_{11}$ satisfies
\begin{equation*}
c_{10}\phi(\bar{t})\al+c_{11}\al\geqslant \Psi(2c_{10}c_{12}\al),
\end{equation*}
 then it is easy to check that
all the conditions of Lemma 3.3 are satisfied. Hence, by Lemma 3.3 (i) we have
\begin{equation}
y(t)\leqslant 2c_{10}\al\phi(t)+2c_{11}\al.
\end{equation}
If $1<p<\f{4+2d}{3d}$, we then have $\frac{d\phi(t)}{dt}\leqslant mg_0^{\frac{4(1-d)}{3}}+\sum\limits_{i=1}^mg_i^{\f{4}{3}+(\f{2}{3}-p)d}(h_i(t))$. Therefore,
\begin{equation*}
(\phi'(t))^{\f{3}{2}}\leqslant c g_0^{2(1-d)}+c\sum\limits_{i=1}^mg_i^{2+(1-\f{3}{2}p)d}(h_i(t)).
\end{equation*}
By (7.15), we have
\[
\begin{split}
&g_i^{2+(1-\f{3}{2}p)d}(h_i(t))\\
&=g_0^{2+(1-\f{3}{2}p)d}+\int_0^{h_i(t)}\f{d}{dt}g_i^{2+(1-\f{3}{2}p)d}(s)ds\\
&=g_0^{2+(1-\f{3}{2}p)d}+(2+(1-\f{3}{2}p)d)\int_0^{h_i(t)}g'_i(s)g_i^{1+(1-\f{3}{2}p)d}(s)ds\\
&\leqslant g_0^{2+(1-\f{3}{2}p)d}+c\gamma\int_0^{h_i(t)}g_i^{(1-p)d-1}(s)ds\\
&\leqslant c\gamma\phi(t)+c.
\end{split}
\]
Thus, if $\gamma$ is sufficient, then all the conditions of Lemma 3.3 (i) are satisfied. Hence, by Lemma 3.3 (i) we
arrive at
\begin{equation}
y(t)\leqslant 2c_{10}\al\phi(t)+2c_{11}\al.
\end{equation}
Finally, since $\int_{\Omega(t)}(|\nabla\mathbf{a}|^2+|\nabla\mathbf{a}|^p)\zeta(x,t)dx\leqslant c\al\phi(t)$, we get
\begin{equation*}
Q(t)\equiv \int_{\Omega(t)}(|\nabla\mathbf{v}|^2+|\nabla\mathbf{v}|^p)\zeta(x,t)dx\leqslant \bar{c}_{10}\al\phi(t)+\bar{c}_{11}\al.
\end{equation*}
\vskip 0.1 true cm

$\mathbf{Case\ II}$. {\bf $\mu_0=0$, $2<p\leqslant 3-\f{2}{d}$}\\

\vskip 0.1 true cm
In this case $I_1=I_2=0$. By completely analogous treatments in Case I, we can obtain
\begin{align*}
&|I_3|\leqslant\ve\int_{\Omega(t)}|\mcD(\mathbf{u})|^pdx+\int_{\omega_i(t)}|\mcD(\mathbf{u})|^pdx+\int_{\Omega(t)}|\mcD(\mathbf{a})|^pdx.\\
&|I_4|\leqslant c\sum\limits_{i=1}^{m}(g_i^{2-(\f{3}{p}-1)d}(h_i(t))\|\nabla\mathbf{u}\|^3_{p,\omega_i(t)}
+g_i^{2(1-\f{d}{p})}(h_i(t))\|\nabla\mathbf{u}\|^2_{p,\omega_i(t)}),\\
&|I_5|\leqslant c\sum\limits_{i=1}^{m}g_i^{2(1-\f{d}{p})}(h_i(t))\|\nabla\mathbf{u}\|^2_{2,\omega_i(t)},\\
&|I_7|\leqslant c\sum\limits_{i=1}^m \|\nabla\mathbf{u}\|^p_{p,\omega_i(t)}
+\sum\limits_{i=1}^m\|\mathbf{a}\|^{2p'}_{2p',\omega_i(t)},\\
&|I_8|\leqslant \ve\|\zeta\nabla\mathbf{u}\|^p_{p,\Omega(t)}+\|\mathbf{a}\|^{2p'}_{2p',\Omega(t)},\\
&|I_9|\leqslant c\sum\limits_{i=1}^m\Big{(}g_i^{2-(1+\f{1}{p})d}(h_i(t))\|\nabla\mathbf{u}\|_{p,\omega_i(t)}+\|\nabla\mathbf{u}\|^p_{p,\omega_i(t)}
+g_i^{2-(\f{3}{p}-1)d}(h_i(t))\|\nabla\mathbf{u}\|^3_{p,\omega_i(t)}\\
&\qquad +\|\mathbf{a}\|^{2p'}_{2p',\Omega_i(t)}+\|\nabla\mathbf{a}\|^p_{p,\Omega_i(t)}\Big{)}.
\end{align*}
For $I_6$, by Young inequality and Lemma 3.1, we get that
\[
\begin{split}
&|I_6|=\Big{|}\int_{\Omega(t)}\mathbf{u}\cdot\nabla\mathbf{u}\cdot\mathbf{a}\zeta dx\Big{|}\\
&\leqslant \ve\int_{\Omega(t)}\zeta|\mathbf{a}|^{\f{p}{d-1}}|\mathbf{u}|^pdx+\ve\int_{\Omega(t)}\zeta|\nabla\mathbf{u}|^pdx+c\int_{\Omega(t)}|\mathbf{a}|^{\f{(d-2)p}{(d-1)(p-2)}}dx.\\
&\leqslant c\ve\int_{\Omega(t)}\zeta|\nabla\mathbf{u}|^pdx+c\sum\limits_{i=1}^{m}\int_{\omega_i(t)}|\nabla\mathbf{u}|^pdx+c\int_{\Omega(t)}|\mathbf{a}|^{\f{(d-2)p}{(d-1)(p-2)}}dx.
\end{split}
\]
As in Case I, based on the estimates on $I_i$ $(1\le i\le 9)$, we have
\begin{equation}
\begin{split}
\|\nabla\mathbf{u}\|^p_{p,\Omega(t)}
\leqslant &c\sum\limits_{i=1}^m \Big{(}g_i^{2-(1+\f{1}{p})d}(h_i(t))\|\nabla\mathbf{u}\|_{p,\omega_i(t)}
+g_i^{2(1-\f{d}{p})}(h_i(t))\|\nabla\mathbf{u}\|^2_{p,\omega_i(t)}\\
& +g_i^{2-(\f{3}{p}-1)d}(h_i(t))\|\nabla\mathbf{u}\|^3_{p,\omega_i(t)}+\|\nabla\mathbf{u}\|^p_{p,\omega_i(t)}\Big{)}+c\|\nabla\mathbf{a}\|^p_{p,\Omega(t)}\\
&+c\|\mathbf{a}\|^{2p'}_{2p',\Omega(t)}+c\|\mathbf{a}\|^{\f{(d-2)p}{(d-1)(p-2)}}_{\f{(d-2)p}{(d-1)(p-2)},\Omega(t)}.
\end{split}
\end{equation}
In addition, it is easy to get
\begin{equation}
\frac{dy(t)}{dt}\geqslant\f{1}{2}\sum\limits_{i=1}^m g_i^{\f{2p}{3}-(1-\f{p}{3})d}(h_i(t))\int_{\omega_i(t)}|\nabla\mathbf{u}|^pdx,
\end{equation}
where $y(t)=\int_{\Omega(t)}|\nabla\mathbf{u}|^p\zeta(x,t)dx$. Hence
\begin{equation}
\sum\limits_{i=1}^m g_i^{2-(\f{3}{p}-1)d}(h_i(t))\|\nabla\mathbf{u}\|_{p,\omega_i(t)}^3\leqslant c(\frac{dy(t)}{dt})^{3/p},
\end{equation}
\begin{equation}
\sum\limits_{i=1}^m g_i^{2(1-\f{d}{p})}\|\nabla\mathbf{u}\|_{p,\omega_i(t)}^2\leqslant  g_0^{\f{2}{3}(1-d)}\sum\limits_{i=1}^mg_i^{\f{4}{3}-(\f{2}{p}-\f{2}{3})d}\|\nabla\mathbf{u}\|_{p,\omega_i(t)}^2
\leqslant c(\frac{dy(t)}{dt})^{2/p},
\end{equation}
and
\begin{equation}
\sum\limits_{i=1}^m\|\nabla\mathbf{u}\|^p_{p,\omega_i(t)}\leqslant g_0^{-\f{4}{3}+(\f{2}{p}-\f{2}{3})d}(\f{dy(t)}{dt}).
\end{equation}
By virtue of
\begin{equation}
\|\nabla\mathbf{a}\|^{2p'}_{2p',\Omega(t)}+\|\mathbf{a}\|^{\f{(d-2)p}{(d-1)(p-2)}}_{\f{(d-2)p}{(d-1)(p-2)},\Omega(t)}
+\|\nabla\mathbf{a}\|^p_{p,\Omega(t)}\leqslant c\al
\sum\limits_{i=1}^m\int\limits_0^{h_i(t)}g_i^{(1-p)d-1}(s)ds,
\end{equation}
we have
\begin{equation}
y(t)\leqslant c[y'(t)+y'^{2/p}+y'^{3/p}(t)]+c\al\phi(t).
\end{equation}
Note that
\begin{equation}
\frac{d\phi(t)}{dt}=\sum\limits_{i=1}^m\frac{h'_i(t)}{g_i^{(p-1)d+1}(h_i(t))}=\sum\limits_{i=1}^m
g_i^{\f{2p}{3}(1-d)}(h_i(t)),
\end{equation}
we then  have $\frac{d\phi(t)}{dt}\leqslant mg_0^{\f{2p}{3}(1-d)}\equiv c_{14}$.
Thus, if we set $\Psi(t)=c_{15}(t+t^{\f{2}{p}}+t^{\f{3}{p}})$, $\varphi(t)=2c_{10}\al\phi(t)+2c_{11}\al$
and $\delta=\frac{1}{2}$, $t_0=\bar{t}$, where $c_{11}$ satisfies
\begin{equation*}
c_{10}\al\phi(\bar{t})+c_{11}\al\geqslant \Psi(2c_{10}c_{14}\al),
\end{equation*}
then it is easy to check that
all the conditions of Lemma 3.3 are fulfilled. Hence, by Lemma 3.3 (i) we have
\begin{equation}
y(t)\leqslant 2c_{10}\al\phi(t)+2c_{11}\al.
\end{equation}
Finally, since $\int_{\Omega(t)}|\nabla\mathbf{a}|^p\zeta(x,t)dx\leqslant c\al\phi(t)$, we get
\begin{equation*}
Q(t)\equiv \int_{\Omega(t)}(|\nabla\mathbf{v}|^2+|\nabla\mathbf{v}|^p)\zeta(x,t)dx\leqslant \bar{c}_{10}\al\phi(t)+\bar{c}_{11}\al.
\end{equation*}
From (7.32) and (7.41), completely similar to the proof of Part 2 in $\S 6$, we can establish
the existence of the solution $\mathbf{v}$ to problem $(2.5)_1-(2.5)_4$, here the details are omitted.
Thus, Theorem 7.1 is shown. \qed

{\bf Remark 7.1.} {\it
If there is a constant $c>0$ such that $|g_i(t)|<c$ for $1\le i\le N$, then Theorem 7.1 is
also true for $\mu_0=0$ and $p>2$. Actually, from the proof of Theorem 7.1, in the case
of $\mu_0=0$ and $p>2$, we just
only need to treat the term $\int_{\Omega(t)}\mathbf{u}\cdot\nabla\mathbf{u}\cdot\mathbf{a}\zeta dx$
since the other terms can be estimated analogously.
For  $|g_i(t)|<c$, by $h_i(t)\thicksim t$ and $\phi(t)\thicksim t$ one has
\begin{align*}
&\Big{|}\int_{\Omega(t)}\mathbf{u}\cdot\nabla\mathbf{u}\cdot\mathbf{a}\zeta dx\Big{|}
\leqslant c\al\int_{\Omega(t)}|\nabla\mathbf{u}|^2\zeta dx+c\al\int_{\Omega(t)}|\mathbf{u}|^2\zeta dx\\
&\leqslant c\al\int_{\Omega(t)}|\nabla\mathbf{u}|^2\zeta dx+c\al\sum\limits_{i=1}^m\int_{\omega_i(t)}|\nabla\mathbf{u}|^2dx\\
&\leqslant \ve \int_{\Omega(t)}|\nabla\mathbf{u}|^p\zeta dx+c\sum\limits_{i=1}^m\int_{\omega_i(t)}|\nabla\mathbf{u}|^p
+c\al t+c\al\\
&\leqslant  \ve \int_{\Omega(t)}|\nabla\mathbf{u}|^p\zeta dx+c\sum\limits_{i=1}^m\int_{\omega_i(t)}|\nabla\mathbf{u}|^p+c\al\phi(t)+c\al.
\end{align*}
Hence the crucial estimate (7.14) still holds.}

{\bf Remark 7.2.} {\it  If $1<p<\frac{4+2d}{3d}$, in order to get the existence of solution to problem (2.5),
we need both the conditions (7.7) and (7.15).
It is worth noting that there are some $g_i's$ such that
conditions (7.7) and (7.15) are satisfied.
For examples, choosing $g_i(t)=\gamma(t+1)^{\al}$, then it is easy to check that (7.7) and (7.15)
are satisfied when $0\leqslant\al\leqslant\f{2}{6-pd}$ and $\gamma$ is sufficient small.
}

{\bf Theorem 7.2.}
{\it Let $\mathbf{v}$ be a weak solution of the system $(2.5)_1-(2.5)_4$. In addition, we assume that
\begin{equation}
\begin{split}
&\text{$\liminf\limits_{t\rightarrow\infty}t^{-3}\int_{\Omega_i(t)}(|\nabla\mathbf{v}|^2+|\nabla\mathbf{v}|^p)dx=0$,\qquad\quad if $\mu_0>0$, $p>1$;}\\
&\text{$\liminf\limits_{t\rightarrow\infty}t^{-\f{3}{3-p}}\int_{\Omega_i(t)}|\nabla\mathbf{v}|^pdx=0$, \qquad\qquad\qquad if $\mu_0=0$,
$2< p\leqslant 3-\f{2}{d}$}.
\end{split}
\end{equation}
Then
\begin{equation}
\label{problem 2: unique modified estimate}
\begin{split}
&\text{$\int_{\Omega_i(t)}(|\nabla\mathbf{v}|^2+|\nabla\mathbf{v}|^p)dx\leqslant c_{15}\al\phi_i(t)+c_{16}\al$,\qquad if $\mu_0>0$, $p>1$;}\\
&\text{$\int_{\Omega_i(t)}|\nabla\mathbf{v}|^pdx\leqslant c_{15}\al\phi_i(t)+c_{16}\al$, \qquad\qquad\qquad if $\mu_0=0$, $2< p\leqslant 3-\f{2}{d}$},
\end{split}
\end{equation}
where
\[
\begin{split}
&\text{$\phi_i(t)=\int_0^t (g_i^{-(d+1)}(s)+g_i^{(1-p)d-1}(s))ds$, \qquad if $\mu_0>0$, $p>1$;}\\
&\text{$\phi_i(t)=\int_0^t g_i^{(1-p)d-1}(s)ds$, \qquad\qquad\qquad\qquad if $\mu_0=0$, $2< p\leqslant 3-\f{2}{d}$.}
\end{split}
\]
}

{\bf Proof.} Let $t\geqslant \bar{t}_i$, where $\bar{t}_i$ is a constant such that $h_i(\bar{t}_i)\geqslant4\beta g_i(0)$. Then we can define that
\begin{equation}
\hat{\zeta}(x^i,t)=
\begin{cases}
\frac{h_i(t)-x^i_1}{g_i(h_i(t))} &\text{for $x_1^i\in[h_i(t)-\beta g_i(h_i(t)),h_i(t)]$},\\
g_i^{-1}(0)x_1^i &\text{for $x_1^i\in[0,\beta g_i(0)]$},\\
\beta &\text{for $x^i_1\in[\beta g_i(0), h_i(t)-\beta g_i(h_i(t))]$}
\end{cases}
\end{equation}
and
\begin{equation*}
\begin{split}
&\text{$\hat{y}_i(t)=\int_{\hat{\Omega}_i(t)}(|\nabla\mathbf{u}|^2+|\nabla\mathbf{u}|^p)\hat{\zeta}dx$, \qquad if $\mu_0>0$, $p>1$,}\\
&\text{$\hat{y}_i(t)=\int_{\hat{\Omega}_i(t)}|\nabla\mathbf{u}|^p\hat{\zeta}dx$, \qquad\qquad\qquad\ if $\mu_0=0$, $2< p\leqslant 3-\f{2}{d}$,}
\end{split}
\end{equation*}
where $\hat{\Omega}_i(t)=\Omega_i(h_i(t))$.
Similarly to the proof of Theorem 7.1,
we can obtain that\\
if $\mu_0>0$, $p>1$,
\begin{multline}
\hat{y}_i(t) \leqslant c\|\nabla\mathbf{u}\|^2_{2,\omega_i(t)}
+c\|\nabla\mathbf{u}\|^p_{p,\omega_i(t)}+cg^{\f{4-d}{2}}_i(h_i(t))\|\nabla\mathbf{u}\|^3_{2,\omega_i(t)}\\
+c\|\nabla\mathbf{a}\|^2_{2,\hat{\Omega}_i(t)}+c\|\mathbf{a}\|^4_{4,\hat{\Omega}_i(t)}+c\|\nabla\mathbf{a}\|^p_{p,\hat{\Omega}_i(t)};
\end{multline}
if $\mu_0=0$, $2<p\leqslant 3-\f{2}{d}$,
\begin{equation}
\begin{split}
\hat{y}_i(t) \leqslant&c\big{(}g_i^{2-(1+\f{1}{p})d}(h_i(t))\|\nabla\mathbf{u}\|_{p,\omega_i(t)}+g_i^{2(1-\f{d}{p})}(h_i(t))\|\nabla\mathbf{u}\|^2_{p,\omega_i(t)}\\
&+g_i^{2-(\f{3}{p}-1)d}(h_i(t))\|\nabla\mathbf{u}\|^3_{p,\omega_i(t)}
+\|\nabla\mathbf{u}\|^p_{p,\omega_i(t)}\big{)}\\
&+c\big{(}\|\nabla\mathbf{a}\|^p_{p,\hat{\Omega}_i(t)}+\|\mathbf{a}\|^{2p'}_{2p',\hat{\Omega}_i(t)}
+\|\mathbf{a}\|^{\f{(d-2)p}{(d-1)(p-2)}}_{\f{(d-2)p}{(d-1)(p-2)},\hat{\Omega}_i(t)}\big{)}.
\end{split}
\end{equation}
Hence
\begin{equation*}
\begin{split}
&\text{$\hat{y}_i(t)\leqslant c(\hat{y}'_i(t)+(\hat{y}'_i(t))^{3/2})+c\al\hat{\phi}_i(t)$, \qquad\qquad\qquad\quad if $\mu_0>0$, $p>1$};\\
&\text{$\hat{y}_i(t)\leqslant c(\hat{y}'_i(t)+(\hat{y}'_i(t))^{2/p}+(\hat{y}'_i(t))^{3/p})+c\al\hat{\phi}_i(t)$,\qquad if $\mu_0=0$, $2< p\leqslant 3-\f{2}{d}$,}
\end{split}
\end{equation*}
where $\hat{\phi}_i(t)=\phi_i(h_i(t))$. Similarly to the proof of Theorem 7.1, setting $\Psi(\tau)=c(\tau+\tau^{\f{3}{2}})$, if $\mu_0>0$; $\Psi(\tau)=c(\tau+\tau^{\f{2}{p}}+\tau^{\f{3}{p}})$, if $\mu_0=0$,
and taking $\varphi_i(t)=2c\al\hat{\phi}_i(t)+2c\al$. By virtue of (7.42), from Lemma 3.3 (ii),
we can arrive at
\begin{equation}
\text{$\hat{y}_i(t)\leqslant 2c\al\hat{\phi}_i(t)+2c\al$,\quad for $t\geqslant\bar{t}_i$.}
\end{equation}
By (7.47), we have that:\\
if $\mu_0>0$, $p>1$,
\begin{equation}
\beta\int_{\beta g_i(0)}^{x_1-\beta g(x_1)}ds\int_{\Sigma_i(s)}(|\nabla\mathbf{u}|^2
+|\nabla\mathbf{u}|^p)dy\leqslant c\al\int_0^{x_1} (g_i^{-(d+1)}(s)+g_i^{(1-p)d-1}(s))ds+c\al;
\end{equation}
if $\mu_0=0$, $2<p\leqslant 3-\f{2}{d}$,
\begin{equation}
\beta\int_{\beta g_i(0)}^{x_1-\beta g(x_1)}ds\int_{\Sigma_i(s)}
|\nabla\mathbf{u}|^pdy\leqslant c \al\int_0^{x_1} g_i^{(1-p)d-1}(s)ds+c\al.
\end{equation}
Meanwhile, if $\mu_0>0$, $p>1$,
\begin{equation}
\begin{split}
&\int_{x_1-\beta g(x_1)}^{x_1}(g^{-(d+1)}_i(s)+g^{(1-p)d-1}_i(s))ds\no\\
&\leqslant \max\limits_{[x_1-\beta g(x_1),x_1]} (g^{-(d+1)}_i(s)+g^{(1-p)d-1}_i(s))\cdot \beta g_i(x_1)\\
&\leqslant c(g^{-d}_i(x_1)+g^{(1-p)d}_i(x_1))\\
&\leqslant  c(g^{-d}_0+g^{(1-p)d}_0);
\end{split}
\end{equation}
if $\mu_0=0$, $2< p\leqslant 3-\f{2}{d}$,
\begin{equation}
\begin{split}
&\int_{x_1-\beta g(x_1)}^{x_1}g^{(1-p)d-1}_i(s)ds\no\\
&\leqslant \max\limits_{[x_1-\beta g(x_1),x_1]} g^{(1-p)d-1}_i(s)\cdot \beta g_i(x_1)\\
&\leqslant cg^{(1-p)d}_i(x_1)\\
&\leqslant c  g^{(1-p)d}_0.
\end{split}
\end{equation}
Therefore, if $\mu_0>0$, $p>1$,
\begin{equation*}
\int_0^t ds\int_{\Sigma(s)}(|\nabla\mathbf{u}|^2+|\nabla\mathbf{u}|^p)dy\leqslant c_{15}\al\phi_i(t)+c_{16}\al;
\end{equation*}
if $\mu_0=0$, $2\leqslant p\leqslant 3-\f{2}{d}$,
\begin{equation*}
\int_0^t ds\int_{\Sigma(s)}|\nabla\mathbf{u}|^pdy\leqslant c_{15}\al\phi_i(t)+c_{16}\al.
\end{equation*}
Since $\int_{\Omega_i(t)}|\nabla\mathbf{a}|^p\zeta(x,t)dx\leqslant c\al\phi_i(t)$, we have that (7.43) is satisfied.  \qquad  \qquad \qquad
\qquad  \qquad  \qquad $\square$

 From (7.14), it is easy to check that the solutions of Theorem 7.1 satisfy (7.42). Hence, we can get the following result:\\
{\bf Theorem 7.3.} {\it we assume that all the conditions of Theorem 7.1 are satisfied, then problem (2.5) has at least one weak solution.}
\subsection{Part 2. Uniqueness}
At first, we establish the following result:\\

{\bf Lemma 7.1.} {\it Assume that $g_i(t)$ satisfies the conditions
\begin{equation}
\label{II uniqueness: small assumption}
\begin{split}
&\text{$|g'_i(t)g^{\f{4}{3}-\f{d}{3}}_i(t)|<\tilde{\gamma}$, \qquad\qquad\ \ if $\mu_0>0$, $p\geqslant\f{4+2d}{3p}$;}\\
&\text{$|g'_i(t)g^{\f{2p}{3}+(\f{p}{3}-1)d}_i(t)|<\tilde{\gamma}$, \qquad if $\mu_0=0$, $2< p\leqslant 3-\f{2}{d}$,
}
\end{split}
\end{equation}
where $\tilde{\gamma}>0$ is a sufficient small constant. If $\mathbf{v}$ satisfies (7.43),
then we have
\begin{equation}
\begin{split}
&\text{$\int_{\hat{\omega}_i(t)}(|\nabla\mathbf{v}|^2+|\nabla\mathbf{v}|^p)dx\leqslant c_{17}\al (g^{\f{4}{3}(1-d)}_i(t)+g^{\f{4}{3}-(p-\f{2}{3})d)}_i(t))$, \quad if  $\mu_0>0$, $p\geqslant\f{4+2d}{3p}$;}\\
&\text{$\int_{\hat{\omega}_i(t)}|\nabla\mathbf{v}|^pdx\leqslant c_{17}\al g^{\f{2}{3}(1-d)}_i(t)$,\quad\qquad\qquad\qquad\qquad\qquad\qquad if $\mu_0=0$, $2< p\leqslant 3-\f{2}{d}$},
\end{split}
\end{equation}
where $\hat{\omega}_i(t)=\{x\in\Omega_i:t-\beta g_i(t)<x_1^i<t\}$ and $t\geqslant4\bar{t}+6\beta g_i(0)$.}\\

{\bf Proof.}   For fixed $t$, we set
\begin{equation*}
\hat{t}=t-\beta g_i(t),\quad t_1=\hat{t}-\beta g_i(\hat{t}),\quad t_2-\beta g_i(t_2)=t,
\end{equation*}
then $\bar{t}\leqslant t_1<\hat{t}<t<t_2$.\\
\vskip 0.1 true cm

$\mathbf{Case \ I}$.  {\bf $\mu_0>0$, $p\geqslant\f{4+2d}{3p}$}\\

\vskip 0.1 true cm
Introduce the functions $m_1(\tau)$ and $m_2(\tau)$ such that
\begin{equation*}
\frac{dm_1(\tau)}{d\tau}=g^{\f{7-d}{3}}_i(t_1-m_1(\tau)), \qquad  \frac{dm_2(\tau)}{d\tau}=g^{\f{7-d}{3}}_i(t_2+m_2(\tau))
\end{equation*}
and  $m_1(0)=m_2(0)=0$.
In addition, constructing the following truncating function
\begin{equation*}
\chi(x^i,\tau)=
\begin{cases}
\ds\frac{x^i_1-t_1+m_1(\tau)}{g_i(t_1-m_1(\tau))},& \text{for $x^i_1\in[t_1-m_1(\tau), t_1-m_1(\tau)+\beta g_i(t_1-m_1(\tau))]$},\\
\ds \qquad \beta, &\text{for $ x^i_1\in[t_1-m_1(\tau)+\beta g_i(t_1-m_1(\tau)),t_2+m_2(\tau)-\beta g_i(t_2+m_2(\tau))]$},\\
\ds\frac{t_2+m_2(\tau)-x_1^i}{g_i(t_2+m_2(\tau))},&\text{for $x^i_1\in[t_2+m_2(\tau)-\beta g_i(t_2+m_2(\tau)), t_2+m_2(\tau)]$}.
\end{cases}
\end{equation*}
Set
\begin{equation*}
\Omega_i(\tau;t)=\{x\in\Omega_i: t_1-m_1(\tau)<x_1^i<t_2+m_2(\tau)\}.
\end{equation*}
As the proof of Theorem 7.1, we multiply the equation in \eqref{problem 2} by $\mathbf{u}\chi$ and integrate over $\Omega_i(\tau;t)$
to get
\begin{equation}
\text{$z(\tau)\leqslant c_{18}[z'(\tau)+(z'(\tau))^{3/2}]+c_{19}
\al\int_{t_1-m_1(\tau)}^{t_2+m_2(\tau)}\big{(}g^{-(d+1)}_i(s)+g^{(1-p)d-1}(s)\big{)}ds$\quad for $\tau\in[0,\tau_1]$,}
\end{equation}
where
\begin{equation*}
z(\tau)=\int_{\Omega_i(\tau;t)}(|\nabla\mathbf{u}|^2+|\nabla\mathbf{u}|^p)\chi(x,\tau)dx,
\end{equation*}
and $\tau_1$ is a constant such that $\bar{t}=t_1-m_1(\tau_1)$.
In addition,
\begin{align*}
&1+\int_0^{\bar{t}}(g_i^{-(d+1)}(s)+g^{(1-p)d-1}_i(s))ds\\
&\leqslant 1+\int_0^{\bar{t}}(g_0^{-(d+1)}+g^{(1-p)d-1}_0)ds\\
&\leqslant (1+\bar{t}(g_0^{-(d+1)}+g^{(1-p)d-1}_0))\frac{(g_i(0)+\beta^{-1}\bar{t})^{d+1}+(g_i(0)
+\beta^{-1}\bar{t})^{(p-1)d+1}}{\bar{t}}\\
&\qquad \times \int_{\bar{t}}^{2\bar{t}}(g_i^{-(d+1)}(s)+g^{(1-p)d-1}_i(s))ds.
\end{align*}
This, together with (7.43) and $2\bar{t}<t<t_2+m_2(\tau)$, yields
\begin{equation}
z(\tau_1)\leqslant c_{20}\al\int_{t_1-m_1(\tau_1)}^{t_2+m_2(\tau_1)}\big{(}g_i^{-(d+1)}(s)+g^{(1-p)d-1}_i(s)\big{)}ds.
\end{equation}
Set
\begin{equation*}
\varphi(\tau)= c_{21}\al\int_{t_1-m_1(\tau)}^{t_2+m_2(\tau)}\big{(}g_i^{-(d+1)}(s)+g_i^{(1-p)d-1}(s)\big{)}ds
+c_{22}\al\big{(}g_i^{\f{4}{3}(1-d)}(t)+g_i^{\f{4}{3}-(p-\f{2}{3})d}(t)\big{)},
\end{equation*}
where $c_{21}=\max\{2c_{19},c_{20}\}$, and $c_{22}$ is a sufficiently large constant. We now prove
\begin{equation}
\varphi(\tau)\geqslant 2c_{18}(\varphi'(\tau)+(\varphi'(\tau))^{3/2}).
\end{equation}
At first, it is easy to get \begin{align*}
\varphi'(\tau)&=c_{21}\al\bigg(\frac{m'_2(\tau)}{g^{d+1}_i(t_2+m_2(\tau))}+\frac{m'_1(\tau)}{g^{d+1}_i(t_1-m_1(\tau))}\bigg)\\
&\qquad +c_{21}\al\bigg(\frac{m'_2(\tau)}{g^{(p-1)d+1}_i(t_2+m_2(\tau))}+\frac{m'_1(\tau)}{g^{(p-1)d+1}_i(t_1-m_1(\tau))}\bigg)\\
=&c_{21}\al\big{(}g_i^{\f{4}{3}(1-d)}(t_2+m_2(\tau))+g_i^{\f{4}{3}(1-d)}(t_1-m_1(\tau))\\
&\qquad +g_i^{\f{4}{3}-(p-\f{2}{3})d}(t_2+m_2(\tau))+g_i^{\f{4}{3}-(p-\f{2}{3})d}(t_1-m_1(\tau))\big{)}\\
\equiv &c_{21}\al A(\tau).
\end{align*}
Since $p\geqslant\f{4+2d}{3p}$, we have $A(\tau)\leqslant 2(g^{\f{4}{3}(1-d)}_0+g^{\f{4}{3}-(p-\f{2}{3})d}_0)$.
Thus $(\varphi'(\tau))^{3/2}\leqslant \sqrt{2}(c_{21}\al)^{\frac{3}{2}}(g^{\f{2}{3}(1-d)}_0
+g^{\f{2}{3}-(\f{p}{2}-\f{1}{3})d}_0) A(\tau)$.
On the other hand, by virtue of (7.50),
\begin{align*}
&A(\tau)=2g_i^{\f{4}{3}(1-d)}(t)-\int_{t_1-m_1(\tau)}^{t}\frac{d}{ds}g_i^{\f{4}{3}(1-d)}(s)ds
+\int_{t}^{t_1+m_2(\tau)}\frac{d}{ds}g_i^{\f{4}{3}(1-d)}(s)ds\\
&\qquad +2g_i^{\f{4}{3}-(p-\f{2}{3})d}(t)-\int_{t_1-m_1(\tau)}^{t}\frac{d}{ds}g_i^{\f{4}{3}-(p-\f{2}{3})d}(s)ds
+\int_{t}^{t_1+m_2(\tau)}\frac{d}{ds}g_i^{\f{4}{3}-(p-\f{2}{3})d}(s)ds\\
&\leqslant2\big{(}g_i^{\f{4}{3}(1-d)}(t)+g_i^{\f{4}{3}-(p-\f{2}{3})d}(t)\big{)}
+c_{22}\tilde{\gamma}\int_{t_1-m_1(\tau)}^{t_2+m_2(\tau)}\big{(}g_i^{-(d+1)}(s)+g_i^{(1-p)d-1}(s)\big{)}ds.
\end{align*}
Hence, if $\tilde{\gamma}$ is sufficient small, then (7.54) holds. This yields for $\tau\in[0,\tau_1]$,
 \begin{equation*}
\text{$z(\tau)\leqslant\varphi(\tau)$.}
\end{equation*}
Choosing $\tau=0$, we then have
\begin{align*}
&\beta\int_{t-\beta g(t)}^tdx_1\int_{\Sigma_i(x_1)}(|\nabla\mathbf{u}|^2+|\nabla\mathbf{u}|^p)dy\\
&\leqslant z(0)=\int_{\Omega_i(0;t)}(|\nabla\mathbf{u}|^2+|\nabla\mathbf{u}|^p)\chi(x,0)dx\\
&\leqslant \varphi(0)=c_{20}\al\int_{t_1}^{t_2}\big{(}g^{-(d+1)}_i(s)+g^{(1-p)d-1}_i(s)\big{)}ds+c_{21}\al\big{(}g_i^{\f{4}{3}(1-d)}(t)+g_i^{\f{4}{3}-(p-\f{2}{3})d}(t)\big{)}.
\end{align*}
Assume that $|g'_i(s)|\leqslant\gamma_1$ for $s\in [t_1,t_2]$. From this, we have
\begin{equation*}
\text{$t_2\leqslant t+\frac{\beta}{1-\beta\gamma_1}g_i(t)$ and $t_1\geqslant t-\beta(2+\beta\gamma_1)g_i(t)$.}
\end{equation*}
Hence, if $c_{23}=\beta[2+\beta\gamma_1+\frac{1}{1-\beta\gamma_1}]<\gamma_1^{-1}$, then
\begin{align*}
&\int_{t_1}^{t_2}\big{(}g^{-(d+1)}_i(s)+g^{(1-p)d-1}_i(s)\big{)}ds\\
&\leqslant \frac{c_{23}}{(1-c_{23}\gamma_1)^{d+1}}g_i^{-d}(t)+\frac{c_{23}}{(1-c_{23}\gamma_1)^{(p-1)d+1}}g_i^{(1-p)d}(t)\\
&\leqslant
c_{24}\big{(}g_i^{\f{4}{3}(1-d)}(t)+g_i^{\f{4}{3}-(p-\f{2}{3})d}(t)\big{)}.
\end{align*}
Therefore
\begin{equation*}
\int_{\hat{\omega}_i(t)}(|\nabla\mathbf{u}|^2+|\nabla\mathbf{u}|^p)dx
\leqslant c\al(g_i^{\f{4}{3}(1-d)}(t)+g_i^{\f{4}{3}-(p-\f{2}{3})d}(t)).
\end{equation*}
Since $\int_{\hat{\omega}_i(t)}(|\nabla\mathbf{a}|^2+|\nabla\mathbf{a}|^p)dx \leqslant c\al(g_i^{\f{4}{3}(1-d)}(t)+g_i^{\f{4}{3}-(p-\f{2}{3})d}(t))$, we have that
\begin{equation*}
\int_{\hat{\omega}_i(t)}(|\nabla\mathbf{v}|^2+|\nabla\mathbf{v}|^p)dx
\leqslant c_{17}\al(g_i^{\f{4}{3}(1-d)}(t)+g_i^{\f{4}{3}-(p-\f{2}{3})d}(t)).
\end{equation*}

\vskip 0.2 true cm

$\mathbf{Case \ II}$. {\bf $\mu_0=0$ , $2< p\leqslant 3-\f{2}{d}$}

\vskip 0.2 true cm

Choosing functions $m_1(\tau)$ and $m_2(\tau)$ such that
\begin{equation*}
\frac{dm_1(\tau)}{d\tau}=g^{\f{2p+3}{3}-(1-\f{p}{3})d}_i(t_1-m_1(\tau)), \qquad  \frac{dm_2(\tau)}{d\tau}=g^{\f{2p+3}{3}-(1-\f{p}{3})d}_i(t_2+m_2(\tau))
\end{equation*}
and $m_1(0)=m_2(0)=0$.
Similarly to the proof of Theorem 7.1, we multiply the equation in \eqref{problem 2} by $\mathbf{u}\chi$ and integrate over $\Omega_i(\tau;t)$
to get
\begin{equation}
\text{$z(\tau)\leqslant c_{25}[z'(\tau)+(z'(\tau))^{2/p}+(z'(\tau))^{3/p}]+c_{26}\al\int_{t_1-m_1(\tau)}^{t_2+m_2(\tau)}g^{(1-p)d-1}(s)ds$\quad for $\tau\in[0,\tau_1]$,}
\end{equation}
where
\begin{equation*}
z(\tau)=\int_{\Omega_i(\tau;t)}|\nabla\mathbf{u}|^p\chi(x,\tau)dx.
\end{equation*}
Note that
\begin{align*}
&1+\int_0^{\bar{t}}g^{(1-p)d-1}_i(s)ds\leqslant 1+\int_0^{\bar{t}}g^{(1-p)d-1}_0ds\\
&\leqslant (1+\bar{t}g^{(1-p)d-1}_0)\frac{(g_i(0)+\beta^{-1}\bar{t})^{(p-1)d+1}}{\bar{t}}\int_{\bar{t}}^{2\bar{t}}g^{(1-p)d-1}_i(s)ds.
\end{align*}
By this inequality and (7.43), we have
\begin{equation}
z(\tau_1)\leqslant c_{27}\al\int_{t_1-m_1(\tau_1)}^{t_2+m_2(\tau_1)}g^{(1-p)d-1}_i(s)ds.
\end{equation}
Set
\begin{equation*}
\varphi(\tau)= c_{29}\al\int_{t_1-m_1(\tau)}^{t_2+m_2(\tau)}g_i^{(1-p)d-1}(s)ds+c_{30}\al g_i^{\f{2p}{3}(1-d)}(t),
\end{equation*}
where $c_{29}=\max\{2c_{26},c_{27}\}$, and $c_{30}$ is a sufficiently large constant. Next we prove
\begin{equation}
\varphi(\tau)\geqslant 2c_{25}(\varphi'(\tau)+(\varphi'(\tau))^{2/p}+(\varphi'(\tau))^{3/p}).
\end{equation}
It follows from a direct computation that
\begin{align*}
\varphi'(\tau)&=c_{29}\al\bigg(\frac{m'_2(\tau)}{g^{(p-1)d+1}_i(t_2+m_2(\tau))}
+\frac{m'_1(\tau)}{g^{(p-1)d+1}_i(t_1-m_1(\tau))}\bigg)\\
=&c_{29}\al\big{(}g_i^{\f{2p}{3}(1-d)}(t_2+m_2(\tau))+g_i^{\f{2p}{3}(1-d)}(t_1-m_1(\tau))\big{)}\\
\equiv &c_{29}\al A(\tau).
\end{align*}
Together with $A(\tau)\leqslant 2g^{\f{2p}{3}(1-d)}_0$, this yields
$(\varphi'(\tau))^{3/p}\leqslant 2^{\f{3}{p}-1}(c_{29}\al)^{\frac{3}{p}}g^{\f{2(3-p)}{3p}(1-d)}_0 A(\tau)$.
Moreover, by (7.50)
\begin{align*}
&A(\tau)=2g_i^{\f{2p}{3}(1-d)}(t)-\int_{t_1-m_1(\tau)}^{t}\frac{d}{ds}g_i^{\f{2p}{3}(1-d)}(s)ds
+\int_t^{t_1+m_2(\tau)}\frac{d}{ds}g_i^{\f{2p}{3}(1-d)}(s)ds\\
&\leqslant2g_i^{\f{2p}{3}(1-d)}(t)+c_{31}\tilde{\gamma}\int_{t_1-m_1(\tau)}^{t_2+m_2(\tau)}g_i^{(1-p)d-1}(s)ds.
\end{align*}
Hence, if $\tilde{\gamma}$ is sufficient small, then (7.57) holds. Therefore,
\begin{equation*}
\text{$z(\tau)\leqslant\varphi(\tau)$,\quad for $\tau\in[0,\tau_1]$}.
\end{equation*}
Choosing $\tau=0$, we have
\begin{align*}
&\beta\int_{t-\beta g(t)}^tdx_1\int_{\Sigma_i(x_1)}|\nabla\mathbf{u}|^pdy\\
&\leqslant z(0)=\int_{\Omega_i(0;t)}|\nabla\mathbf{u}|^p\chi(x,0)dx\\
&\leqslant \varphi(0)=c_{29}\al\int_{t_1}^{t_2}g^{(1-p)d-1}_i(s)ds+c_{30}\al g_i^{\f{2p}{3}(1-d)}(t).
\end{align*}
Assume that $|g'_i(s)|\leqslant\gamma_1$ for $s\in [t_1,t_2]$. Then
\begin{equation*}
\text{$t_2\leqslant t+\frac{\beta}{1-\beta\gamma_1}g(t)$ and $t_1\geqslant t-\beta(2+\beta\gamma_1)g_i(t)$.}
\end{equation*}
Hence, if $c_{32}=\beta[2+\beta\gamma_1+\frac{1}{1-\beta\gamma_1}]\leqslant\gamma_1^{-1}$, then
\begin{equation*}
\int_{t_1}^{t_2}g^{(1-p)d-1}_i(s)ds\leqslant \frac{c_{32}}{(1-c_{32}\gamma_1)^{(p-1)d}}g_i^{(1-p)d}(t)\leqslant
c_{33} g_i^{\f{2p}{3}(1-d)}(t).
\end{equation*}
Therefore,
\begin{equation*}
\int_{\hat{\omega}_i(t)}|\nabla\mathbf{u}|^pdx\leqslant c\al g_i^{\f{2p}{3}(1-d)}(t).
\end{equation*}
It is easy to check that $\int_{\hat{\omega}_i(t)}|\nabla\mathbf{a}|^pdx \leqslant c\al g_i^{\f{4}{3}-(p-\f{2}{3})d}(t)$, hence
\begin{equation*}
\int_{\hat{\omega}_i(t)}|\nabla\mathbf{v}|^pdx\leqslant c_{17}\al g_i^{\f{2p}{3}(1-d)}(t).
\end{equation*}
Collecting all the estimates above, we complete the proof of Lemma 7.1. \qed
\vskip 0.3 true cm
Next, let $\mathbf{v}$ be a divergence free vector field in $W^{1,p}_{loc}(\bar{\Omega})$, and we assume that there exists a constant $\nu$ satisfying $\nu(p-2)\leqslant\f{p-2}{p-1}+(\f{2}{3}+\f{1}{p})d-\f{8}{3}$ such that
\begin{equation}
\text{$|\mcD(\mathbf{v})(x^i)|\geqslant c g_i^{-\nu}(x^i_1)|y^i|^{\f{1}{p-1}}$\quad for $x^i\in\Omega_i$,}
\end{equation}
where $y^i=(x^i_2,...,x^i_d)$ and $1\le i \le N$.

{\bf Theorem 7.4.}  {\it Let $d=3$, we assume that $g_i(t)$ satisfies (7.50) and $\alpha$ is sufficiently small, and
there exists a solution $\mathbf{v}$ satisfying (7.58) to problem (2.5).
Then the solution $\mathbf{v}$ of problem (2.5) is unique for $\mu_0=0$, $2< p\leqslant 3-\f{2}{d}$ .
When $\mu_0>0$, even if the assumption (7.58) on $\mathbf{u}$ is removed,
the solution $\mathbf{v}$ of problem (2.5) exists uniquely for $p\geqslant\f{16-d}{3d}$.\\
}

{\bf Remark 7.3.} {\it If $\mathbf{v}(x_1,y)$ is the Hagen-Poiseuille flow in pipes of
circular cross section with radius $g_i(t)$, then it follows from (7.15) of \cite{Rob} that
\begin{equation}
\text{$|\mcD(\mathbf{v})(x_1,y)|\geqslant c g_i^{-\f{3p-2}{p-1}}(t)|y|^{\f{1}{p-1}}$\quad for some positive number $c>0$}.
\end{equation}
This means that the assumption (7.58) is reasonable under some cases.
}

\vskip 0.1 true cm

{\bf Proof.}  By the assumptions in Theorem 7.4, from Lemma 7.1, since $d=3$, one has
\begin{equation}
\begin{split}
&\text{$\int_{\hat{\omega}_i(t)}|\nabla\mathbf{v}_2|^2dx\leqslant
c\al\big{(}g_i^{\f{4}{3}(1-d)}(t)+g_i^{\f{4}{3}-(p-\f{2}{3})d}\big{)}\leqslant c\al g_i^{d-4}(t)$, \quad if $\mu_0>0$, $p\geqslant\f{16-d}{3d}$;}\\
&\text{$\int_{\hat{\omega}_i(t)}|\nabla\mathbf{v}_2|^pdx\leqslant
c\al g^{\f{2p}{3}(1-d)}_i(t)$, \qquad\qquad\qquad\qquad\qquad\qquad\ \ if $\mu_0=0$, $2< p\leqslant3-\f{2}{d}$.}
\end{split}
\end{equation}
Set $\mathbf{w}=\mathbf{v}_1-\mathbf{v}_2$. Then
\begin{multline}
-\mu_0div(\mcD(\mathbf{w}))-\mu_1div\Big{(}|\mcD(\mathbf{w})+\mcD(\mathbf{v}_2)|^{p-2}\big{(}\mcD(\mathbf{w})+\mcD(\mathbf{v}_2)\big{)}
-|\mcD(\mathbf{v}_2)|^{p-2}\mcD(\mathbf{v}_2)\Big{)}\\
+\mathbf{w}\cdot\nabla\mathbf{w}+\mathbf{w}\cdot\nabla\mathbf{v}_2
+\mathbf{v}_2\cdot\nabla\mathbf{w}+\nabla(\pi_1-\pi_2)=0.
\end{multline}
From the proof of Theorem 7.1, we just need to estimate $(\mathbf{w}\cdot\nabla\mathbf{v}_2,\mathbf{w})_{\hat{\omega}_i(t)}$.
While, by $(7.59)_1$ and Lemma 3.1 if $\mu_0>0$,
\begin{equation}
\begin{split}
&\Big{|}\int_{\hat{\omega}_i(t)}\mathbf{w}\cdot\nabla\mathbf{v}_2\cdot\mathbf{w}dx\Big{|}\\
&\leqslant\|\mathbf{w}\|^2_{4,\hat{\omega}_i(t)}
\|\nabla\mathbf{v}_2\|_{2,\hat{\omega}_i(t)}\\
&\leqslant g_i^{2-\f{1}{2}d}(t)\|\nabla\mathbf{w}\|^2_{2,\hat{\omega}_i(t)}\|\nabla\mathbf{v}_2\|_{2,\hat{\omega}_i(t)}\\
&\leqslant c\al^{\f{1}{2}}\|\nabla\mathbf{w}\|^2_{2,\hat{\omega}_i(t)};
\end{split}
\end{equation}
if $\mu_0=0$, $2< p\leqslant 3-\f{2}{d}$,
\begin{equation}
\begin{split}
&\Big{|}\int_{\hat{\omega}_i(t)}\mathbf{w}\cdot\nabla\mathbf{v}_2\cdot\mathbf{w}dx\Big{|}\\
&\leqslant\|\mathbf{w}\|^2_{2p',\hat{\omega}_i(t)}
\|\nabla\mathbf{v}_2\|_{p,\hat{\omega}_i(t)}\\
&\leqslant g_i^{2+(\f{p-1}{p}-\f{2}{r})d}(t)\|\nabla\mathbf{w}\|^2_{r,\hat{\omega}_i(t)}\|\nabla\mathbf{v}_2\|_{p,\hat{\omega}_i(t)}.
\end{split}
\end{equation}
Meanwhile, if $1\leqslant r<\f{2(d-1)(p-1)}{d p-(d+1)}$, by (7.58)
\begin{equation}
\begin{split}
&\int_{\hat{\omega}_i(t)}|\nabla\mathbf{w}|^rdx\\
&\leqslant\int_{\hat{\omega}_i(t)}|\mcD(\mathbf{v}_2)|^{\f{(p-2)r}{2}}|\mcD(\mathbf{w})|^r|\mcD(\mathbf{v}_2)|^{-\f{(p-2)r}{2}}dx\\
&\leqslant\bigg(\int_{\hat{\omega}_i(t)}|\mcD(\mathbf{v}_2)|^{(p-2)}|\mcD(\mathbf{w})|^2dx\bigg)^{\f{r}{2}}
\bigg(\int_{\hat{\omega}_i(t)}|\mcD(\mathbf{v}_2)|^{-\f{(p-2)r}{2-r}}dx\bigg)^{\f{2-r}{2}}\\
&\leqslant cg_i^{\f{\nu r (p-2)}{2}+\f{(2-r)}{2}}(t)\||\mcD(\mathbf{v}_2)|^{p-2}\mcD(\mathbf{w})\|^r_{2,\hat{\omega}_i(t)}
\bigg(\int_0^{g_i(t)}s^{d-2-\f{(p-2)r}{(p-1)(2-r)}}ds\bigg)^{\f{2-r}{2}}\\
&\leqslant cg_i^{\f{d(2-r)}{2}-\f{(p-2)r}{2(p-1)}+\f{\nu r (p-2)}{2}}(t)\|
|\mcD(\mathbf{v}_2)|^{p-2}\mcD(\mathbf{w})\|^r_{2,\hat{\omega}_i(t)}.
\end{split}
\end{equation}
Hence, we have
\begin{equation}
\Big{|}\int_{\hat{\omega}_i(t)}\mathbf{w}\cdot\nabla\mathbf{v}_2\cdot\mathbf{w}dx\Big{|}\leqslant c\al^{\f{1}{p}} g_i^{\theta}(t)\||\mcD(\mathbf{v}_2)|^{p-2}\mcD(\mathbf{w})\|^2_{2,\hat{\omega}_i(t)},
\end{equation}
where $\theta=-\f{p-2}{p-1}-(\f{2}{3}+\f{1}{p})d+\f{8}{3}+\nu(p-2)\leqslant0$.
Therefore,
\begin{equation}
\Big{|}\int_{\hat{\omega}_i(t)}\mathbf{w}\cdot\nabla\mathbf{v}_2\cdot\mathbf{w}dx\Big{|}\leqslant c\al\||\mcD(\mathbf{v}_2)|^{p-2}\mcD(\mathbf{w})\|^2_{2,\hat{\omega}_i(t)}.
\end{equation}
From (7.61) and (7.65), similarly to the proof of Theorem 7.1, if $\al$ is sufficiently small, we can get
\begin{equation}
\begin{split}
&\text{$y(t)\leqslant c[y'(t)+(y'(t))^{\frac{3}{2}}]$,\qquad\qquad\qquad\quad if $\mu_0>0$, $p\geqslant\frac{16-d}{3d}$}, \\
&\text{$y(t)\leqslant c[y'(t)+(y'(t))^{\frac{2}{p}}+(y'(t))^{\frac{3}{p}}]$, \qquad\ if $\mu_0=0$, $2< p\leqslant 3-\f{2}{d}$,}
\end{split}
\end{equation}
where $y(t)=\int_{\Omega(t)}(|\nabla\mathbf{w}|^2+|\nabla\mathbf{w}|^p)dx$, if $\mu_0>0$; $y(t)=\int_{\Omega(t)}|\nabla\mathbf{w}|^pdx$, if $\mu_0=0$.
If $y(t)$ is not identically zero, it then follows from Lemma 3.3 (iii) that
\begin{equation*}
\text{$\liminf\limits_{t\rightarrow\infty}t^{-3}y(t)>0$, if $\mu_0>0$, $p\geqslant\frac{16-d}{3d}$};
\text{$\liminf\limits_{t\rightarrow\infty}t^{-\f{3}{3-p}}y(t)>0$, if $\mu_0=0$, $2< p\leqslant 3-\f{2}{d}$}.
\end{equation*}
These contradict with \eqref{problem 2}$_5$. Hence, $y(t)\equiv0$, and further $\mathbf{v}_1=\mathbf{v}_2$.
\qquad\qquad \qquad \qquad \qquad \qquad  $\square$


\end{document}